\documentclass{article}
\usepackage{amsfonts}
\usepackage{amsmath}
\usepackage{amssymb}
\usepackage{geometry}
\usepackage{setspace}
\usepackage{sectsty}
\usepackage{lineno}
\usepackage{theorem}
\usepackage{fancyhdr}
\usepackage{graphicx}

\newtheorem{theorem}{Theorem}

\newtheorem{corollary}[theorem]{Corollary}

\newtheorem{lemma}[theorem]{Lemma}

\DeclareMathOperator{\leb}{Leb}

\DeclareMathOperator{\var}{Var}

\input{tcilatex}

\begin{document}

\title{Optimal upper and lower bounds for the true and empirical excess
risks in heteroscedastic least-squares regression}
\author{A. Saumard\thanks{%
Research partly supported by the french Agence Nationale de la Recherche
(ANR 2011 BS01 010 01 projet Calibration).} \\
University Rennes 1, IRMAR\\
adrien.saumard@univ-rennes1.fr}
\date{February 24, 2012}
\maketitle

\begin{abstract}
We consider the estimation of a bounded regression function with
nonparametric heteroscedastic noise and random design. We study the true and
empirical excess risks of the least-squares estimator on finite-dimensional
vector spaces. We give upper and lower bounds on these quantities that are
nonasymptotic and optimal to first order, allowing the dimension to depend
on sample size. These bounds show the equivalence between the true and
empirical excess risks when, among other things, the least-squares estimator
is consistent in sup-norm with the projection of the regression function
onto the considered model. Consistency in the sup-norm is then proved for
suitable histogram models and more general models of piecewise polynomials
that are endowed with a localized basis structure.\bigskip

\noindent \textbf{keywords:} Least-squares regression, Heteroscedasticity,
Excess risk, Lower bounds, Sup-norm, Localized basis, Empirical process.
\end{abstract}

\section{Introduction\label{section_intro_fixed_reg}}

A few years ago, Birg\'{e} and Massart \cite{BirMas:07} introduced a
data-driven calibration method for penalized criteria in model selection,
called the \textit{Slope Heuristics}. Their algorithm is based on the
concept of the minimal penalty, under which a model selection procedure
fails. Given the shape of the ideal penalty, which in their Gaussian setting
is a known function of the dimension of the considered models, the algorithm
first provides a data-driven estimate of the minimal penalty. This is done
by taking advantage of a sudden change in the behavior of the model
selection procedure around this level of penalty. Then, the algorithm
selects a model by using a penalty that is twice the estimated minimal
penalty. Birg\'{e} and Massart prove\ in \cite{BirMas:07} that an
asymptotically optimal penalty is twice the minimal one, in the sense that
the associated selected model achieves a nonasymptotic oracle inequality
with leading constant converging to one when the sample size tends to
infinity.

The slope heuristics algorithm has been recently extended by Arlot and
Massart \cite{ArlotMassart:09} to the selection of M-estimators, whenever
the number of models is not more than polynomial in the sample size. Arlot
and Massart highlight that, in this context, the mean of the empirical
excess risk on each model should be a good, rather general candidate for the
- unknown - minimal penalty. In addition, they note that an optimal penalty
is roughly given by the sum of the true and the empirical excess risks on
each model. A key fact underlying the asymptotic optimality of the slope
heuristics algorithm is the equivalence - in the sense that the ratio tends
to one when the sample size tends to infinity - between the true and
empirical excess risk, for each model which is likely to be selected.
Generally, these models are of moderate\ dimension, typically between $%
\left( \log \left( n\right) \right) ^{c}$ and $n/\left( \log \left( n\right)
\right) ^{c}$, where $c$ is a positive constant and $n$ is the sample size.
This equivalence leads, quite straightforwardly, to the factor two between
the minimal penalty and the optimal one.

Arlot and Massart prove in \cite{ArlotMassart:09}, by considering the
selection of finite-dimensional models of histograms in heteroscedastic
regression with a random design, that the slope heuristics algorithm is
asymptotically optimal. The authors conjecture in \cite{ArlotMassart:09},
Section 1, that the restriction to histograms is \textquotedblleft mainly
technical\textquotedblright , and that the slope heuristics
\textquotedblleft remains valid at least in the general least squares
regression framework\textquotedblright .

The first motivation of the present paper is thus to tackle the challenging
mathematical problem raised by Arlot and Massart in \cite{ArlotMassart:09},
concerning the validity slope heuristics. More precisely, we isolate the
question of the equivalence, for a fixed model, between the true and
empirical excess risks. As emphasized in \cite{ArlotMassart:09}, this
constitutes the principal part of the conjecture, since other arguments
leading to model selection results are now well understood. We thus postpone
model selection issues to a forthcoming paper, and focus on the fixed model
case.

We consider least squares regression with heteroscedastic noise and random
design, using a finite-dimensional linear model. Our analysis is
nonasymptotic in the sense that our results are available for a fixed value
of the sample size. It is also worth noticing that the dimension of the
considered model is allowed to depend on the sample size and consequently is
not treated as a constant of the problem. In order to determine the possible
equivalence between the true and empirical excess risks, we investigate
upper and lower deviation bounds for each quantity. We obtain first order
optimal bounds, thus exhibiting the first part of the asymptotic expansion
of the excess risks. This requires to determine not only the right rates of
convergence, but also the optimal constant on the leading order term. We
give two examples of models that satisfy our conditions: models of
histograms and models of piecewise polynomials, whenever the partition
defining these models satisfy some regularity condition with respect to the
unknown distribution of data. Our results concerning histograms roughly
recover those derived for a fixed model by Arlot and Massart \cite%
{ArlotMassart:09}, but with different techniques. Moreover, the case of
piecewise polynomials strictly extend these results, and thus tends to
confirm Arlot and Massart conjecture on the validity of the slope heuristics.

We believe that our deviation bounds, especially those concerning the true
excess risk, are interesting by themselves. Indeed, the optimization of the
excess risk is, from a general perspective, at the core of many
nonparametric approaches, especially those related to statistical learning
theory. Hence, any sharp control of this quantity is likely to be useful in
many contexts.

In the general bounded M-estimation framework, rates of convergence and
upper bounds for the excess risk are now well understood, see \cite{Tsy:04a}%
, \cite{MassartNedelec:06}, \cite{Kolt:06}, \cite{BarMen:06}, \cite%
{GineKolt:06}. However, the values of the constants in these deviation
bounds are suboptimal - or even unknown -, due in particular to the use of
chaining techniques. Concerning lower deviation bounds, there is no
convincing contribution to our knowledge, except the work of Bartlett and
Mendelson \cite{BarMen:06}, where an additional assumption on the behavior
of underlying empirical process is used to derive such a result. However,
this assumption is in general hard to check.

More specific frameworks, such as least squares regression with a fixed
design on linear models (see for instance \cite{BirMas:07}, \cite%
{BarGirHuet:09} and \cite{Arl_Bac:2009}), least squares estimation of
density on linear models (see \cite{Lerasle:09a} and references therein), or
least squares regression on histograms as in \cite{ArlotMassart:09}, allow
for sharp, explicit computations that lead to optimal upper and lower bounds
for the excess risks. Hence, a natural question is: is there a framework,
between the general one and the special cases, that would allow to derive
deviation bounds that are optimal at the first order ? In other words, how
far could optimal results concerning deviation bounds been extended ? The
results presented in this article can be seen as a first attempt to answer
these questions.

The article is organized as follows. We present the statistical framework in
Section \ref{section_framework_reg}, where we show in particular the
existence of an expansion of the least squares regression contrast into the
sum of a linear and a quadratic part. In Section \ref%
{section_outline_approach}, we detail the main steps of our approach at a
heuristic level and give a summary of the results presented in the paper. We
then derive some general results in Section \ref{section_reg_fixed_model}.
These theorems are then applied to the case of histograms and piecewise
polynomials in Sections \ref{section_histogram_case_fixed} and \ref%
{section_piecewise_case_fixed} respectively, where in particular, explicit
rates of convergence in sup-norm are derived. Finally, the proofs are
postponed to the end of the article.

\section{Regression framework and notations\label{section_framework_reg}}

\subsection{least squares estimator}

Let $\left( \mathcal{X},\mathcal{T}_{\mathcal{X}}\right) $ be a measurable
space and set $\mathcal{Z}=\mathcal{X\times }\mathbb{R}$. We assume that $%
\xi _{i}=\left( X_{i,}Y_{i}\right) \in \mathcal{X\times }\mathbb{R}$, $i\in
\left\{ 1,...,n\right\} $, are $n$ i.i.d. observations with distribution $P$%
. The marginal law of $X_{i}$ is denoted by $P^{X}.$ We assume that the data
satisfy the following relation%
\begin{equation}
Y_{i}=s_{\ast }\left( X_{i}\right) +\sigma \left( X_{i}\right) \varepsilon
_{i}\text{ },  \label{regression_model}
\end{equation}%
where $s_{\ast }\in L_{2}\left( P^{X}\right) $, $\varepsilon _{i}$ are
i.i.d. random variables with mean $0$ and variance $1$ conditionally to $%
X_{i}$ and $\sigma :$ $\mathcal{X\longrightarrow }\mathbb{R}$ is a
heteroscedastic noise level. A generic random variable of law $P$,
independent of $\left( \xi _{1},...,\xi _{n}\right) $, is denoted by $\xi
=\left( X_{,}Y\right) .$

\noindent Hence, $s_{\ast }$ is the regression function of $Y$ with respect
to $X$, to be estimated. Given a finite dimensional linear vector space $M$,
that we will call a \textquotedblleft model\textquotedblright , we denote by 
$s_{M}$ the linear projection of $s_{\ast }$ onto $M$ in $L_{2}\left(
P^{X}\right) $ and by $D$ the linear dimension of $M$.

\noindent We consider on the model $M$ a least squares estimator $s_{n}$
(possibly non unique), defined as follows 
\begin{equation}
s_{n}\in \arg \min_{s\in M}\left\{ \frac{1}{n}\sum_{i=1}^{n}\left(
Y_{i}-s\left( X_{i}\right) \right) ^{2}\right\} \text{ }.
\label{def_estimator_bis}
\end{equation}%
So, if we denote by 
\begin{equation*}
P_{n}=\frac{1}{n}\sum_{i=1}^{n}\delta _{\left( X_{i},Y_{i}\right) }
\end{equation*}%
the empirical distribution of the data and by $K:L_{2}\left( P^{X}\right)
\longrightarrow L_{1}\left( P\right) $ the least squares contrast, defined
by 
\begin{equation*}
K\left( s\right) =\left( x,y\right) \in \mathcal{Z}\rightarrow \left(
y-s\left( x\right) \right) ^{2}\text{, \ \ }s\in L_{2}\left( P^{X}\right) 
\text{,\ }
\end{equation*}%
we then remark that $s_{n}$ belongs to the general class of M-estimators, as
it satisfies%
\begin{equation}
s_{n}\in \arg \min_{s\in M}\left\{ P_{n}\left( K\left( s\right) \right)
\right\} \text{ }.  \label{def_estimator_2_bis}
\end{equation}

\subsection{Excess risk and contrast\label{section_excess_risk}}

As defined in (\ref{def_estimator_2_bis}), $s_{n}$ is the empirical risk
minimizer of the least squares contrast. The regression function $s_{\ast }$
can be defined as the minimizer in $L_{2}\left( P^{X}\right) $ of the mean
of the contrast over the unknown law $P$,%
\begin{equation*}
s_{\ast }=\arg \min_{s\in L_{2}\left( P^{X}\right) }PK\left( s\right) \text{ 
},
\end{equation*}%
where 
\begin{equation*}
PK\left( s\right) =P\left( Ks\right) =PKs=\mathbb{E}\left[ K\left( s\right)
\left( X,Y\right) \right] =\mathbb{E}\left[ \left( Y-s\left( X\right)
\right) ^{2}\right]
\end{equation*}%
is called the risk of the function $s$. In particular we have $PKs_{\ast }=%
\mathbb{E}\left[ \sigma ^{2}\left( X\right) \right] $. We first notice that
for any $s\in L_{2}\left( P^{X}\right) $, if we denote by 
\begin{equation*}
\left\Vert s\right\Vert _{2}=\left( \int_{\mathcal{X}}s^{2}dP^{X}\right)
^{1/2}
\end{equation*}%
its quadratic norm, then we have, by (\ref{regression_model}) above, 
\begin{align*}
PKs-PKs_{\ast }& =P\left( Ks-Ks_{\ast }\right) \\
& =\mathbb{E}\left[ \left( Y-s\left( X\right) \right) ^{2}-\left( Y-s_{\ast
}\left( X\right) \right) ^{2}\right] \\
& =\mathbb{E}\left[ \left( s_{\ast }-s\right) ^{2}\left( X\right) \right] +2%
\mathbb{E}\left[ \left( s_{\ast }-s\right) \left( X\right) \underset{=0}{%
\underbrace{\mathbb{E}\left[ Y-s_{\ast }\left( X\right) \left\vert X\right. %
\right] }}\right] \\
& =\left\Vert s-s_{\ast }\right\Vert _{2}^{2}\geq 0\text{ }.
\end{align*}%
The quantity $PKs-PKs_{\ast }$ is called the excess risk of $s$. Now, if we
denote by $s_{M}$ the linear projection of $s_{\ast }$ onto $M$ in $%
L_{2}\left( P^{X}\right) $, we have%
\begin{equation}
PKs_{M}-PKs_{\ast }=\inf_{s\in M}\left\{ PKs-PKs_{\ast }\right\} \text{ },
\label{def_sM}
\end{equation}%
and for all $s\in M$%
\begin{equation}
P^{X}\left( s%
\cdot%
\left( s_{M}-s_{\ast }\right) \right) =0\text{ }.  \label{projectionM}
\end{equation}%
From (\ref{def_sM}), we deduce that%
\begin{equation*}
s_{M}=\arg \min_{s\in M}PK\left( s\right) \text{ }.
\end{equation*}

\noindent Our goal\ is to study the performance of the least squares
estimator, that we measure by its excess risk. So we are mainly interested
in the random quantity $P\left( Ks_{n}-Ks_{\ast }\right) .$ Moreover, as we
can write%
\begin{equation*}
P\left( Ks_{n}-Ks_{\ast }\right) =P\left( Ks_{n}-Ks_{M}\right) +P\left(
Ks_{M}-Ks_{\ast }\right)
\end{equation*}%
we naturally focus on the quantity%
\begin{equation*}
P\left( Ks_{n}-Ks_{M}\right) \geq 0
\end{equation*}%
that we want to bound from upper and from below, with high probability. We
will often call this last quantity the excess risk of the estimator on $M$
or the true excess risk of $s_{n}$, in opposition to the empirical excess
risk\textbf{\ }for which the expectation is taken over the empirical
measure, 
\begin{equation*}
P_{n}\left( Ks_{M}-Ks_{n}\right) \geq 0\text{ }.
\end{equation*}

\noindent The following lemma establishes the expansion of the regression
contrast around $s_{M}$ on $M$. This expansion exhibits a linear part and a
quadratic parts.

\begin{lemma}
\label{dev_centered_sM copy(1)_reg}We have, for every $z=\left( x,y\right)
\in\mathcal{Z},$%
\begin{equation}
\left( Ks\right) \left( z\right) -\left( Ks_{M}\right) \left( z\right)
=\psi_{1,M}\left( z\right) \left( s-s_{M}\right) \left( x\right)
+\psi_{2}\left( \left( s-s_{M}\right) \left( x\right) \right)
\label{dev_contrast}
\end{equation}
with $\psi_{1,M}\left( z\right) =-2\left( y-s_{M}\left( x\right) \right) $
and $\psi_{2}\left( t\right) =t^{2}$, for all $t\in\mathbb{R}.$ Moreover,
for all $s\in M$,%
\begin{equation}
P\left( \psi_{1,M}%
\cdot%
s\right) =0\text{ }.  \label{center_M}
\end{equation}
\end{lemma}

\noindent \textbf{Proof.}\ Start with%
\begin{align*}
& \left( Ks\right) \left( z\right) -\left( Ks_{M}\right) \left( z\right) \\
& =\left( y-s\left( x\right) \right) ^{2}-\left( y-s_{M}\left( x\right)
\right) ^{2} \\
& =\left( \left( s-s_{M}\right) \left( x\right) \right) \left( \left(
s-s_{M}\right) \left( x\right) -2\left( y-s_{M}\left( x\right) \right)
\right) \\
& =-2\left( y-s_{M}\left( x\right) \right) \left( \left( s-s_{M}\right)
\left( x\right) \right) +\left( \left( s-s_{M}\right) \left( x\right)
\right) ^{2}\text{ ,}
\end{align*}%
which gives (\ref{dev_contrast}). Moreover, observe that for any $s\in M$,%
\begin{equation}
P\left( \psi _{1,M}%
\cdot%
s\right) =-2\mathbb{E}\left[ \left( Y-s_{\ast }\left( X\right) \right)
s\left( X\right) \right] +2\mathbb{E}\left[ s\left( X\right) \left(
s_{M}-s_{\ast }\right) \left( X\right) \right] \text{ }.  \label{dev_reg_1}
\end{equation}%
We have%
\begin{equation}
\mathbb{E}\left[ \left( Y-s_{\ast }\left( X\right) \right) s\left( X\right) %
\right] =\mathbb{E}\left[ \underset{=0}{\underbrace{\mathbb{E}\left[ \left(
Y-s_{\ast }\left( X\right) \right) \left\vert X\right. \right] }}s\left(
X\right) \right] =0\text{ }.  \label{dev_reg_2}
\end{equation}%
and, by (\ref{projectionM}),%
\begin{equation}
\mathbb{E}\left[ s\left( X\right) \left( s_{M}-s_{\ast }\right) \left(
X\right) \right] =P^{X}\left( s%
\cdot%
\left( s_{M}-s_{\ast }\right) \right) =0\text{ }.  \label{dev_reg_3}
\end{equation}%
Combining (\ref{dev_reg_1}), (\ref{dev_reg_2}) and (\ref{dev_reg_3}) we get
that for any $s\in M$, $P\left( \psi _{1,M}%
\cdot%
s\right) =0$. This concludes the proof.\ 
$\blacksquare$%

\section{Outline of the approach\label{section_outline_approach}}

\noindent Having introduced the framework and notations in Section \ref%
{section_framework_reg} above, we are now able to explain more precisely the
major steps of our approach to the problem of deriving optimal upper and
lower bounds for the excess risks. As mentioned in the introduction, one of
our main motivations is to determine whether the true excess risk is
equivalent to the empirical one or not:%
\begin{equation}
P\left( Ks_{n}-Ks_{M}\right) \sim P_{n}\left( Ks_{M}-Ks_{n}\right) \text{ }?
\label{question_equi}
\end{equation}%
Indeed, such an equivalence is a keystone to justify the slope heuristics, a
data-driven calibration method first proposed by Birg\'{e} and Massart \cite%
{BirMas:07} in a Gaussian setting and then extended by Arlot and Massart 
\cite{ArlotMassart:09} to the selection of M-estimators.

\noindent The goal of this section is twofold. Firstly, it helps the reader
to understand the role of the assumptions made in the forthcoming sections.
Secondly, it provides an outline of the proof of our main result, Theorem %
\ref{principal_reg} below. We suggest the reader interested in our proofs to
read this section before entering the proofs.

\noindent We start by rewriting the lower and upper bound problems, for the
true and empirical excess risks. Let $C$ and $\alpha $ be two positive
numbers. The question of bounding the true excess risk from upper and with
high probability can be stated as follows: find, at a fixed $\alpha >0$, the
smallest $C>0$ such that 
\begin{equation*}
\mathbb{P}\left[ P\left( Ks_{n}-Ks_{M}\right) >C\right] \leq n^{-\alpha }%
\text{ .}
\end{equation*}%
We then write, by definition of the M-estimator $s_{n}$ as a minimizer of
the empirical excess risk over the model $M$,%
\begin{eqnarray}
&&\mathbb{P}\left[ P\left( Ks_{n}-Ks_{M}\right) >C\right]  \notag \\
&\leq &\mathbb{P}\left[ \inf_{s\in M_{C}}P_{n}\left( Ks-Ks_{M}\right) \geq
\inf_{s\in M_{>C}}P_{n}\left( Ks-Ks_{M}\right) \right]  \notag \\
&=&\mathbb{P}\left[ \sup_{s\in M_{C}}P_{n}\left( Ks_{M}-Ks\right) \leq
\sup_{s\in M_{>C}}P_{n}\left( Ks_{M}-Ks\right) \right] \text{ ,}
\label{reecriture_1}
\end{eqnarray}%
where%
\begin{equation*}
M_{C}:=\left\{ s\in M\text{ };\text{ }P\left( Ks-Ks_{M}\right) \leq C\right\}
\end{equation*}%
and%
\begin{equation*}
M_{>C}:=M\backslash M_{C}=\left\{ s\in M\text{ };\text{ }P\left(
Ks-Ks_{M}\right) >C\right\}
\end{equation*}%
are subsets of the model $M$, localized in terms of excess risk. As a matter
of fact, $M_{C}$ is the closed ball of radius $C$ in $\left( M,\left\Vert 
\cdot%
\right\Vert _{2}\right) $. In the same manner, the question of bounding the
true excess risk from below and with high probability is formalized as
follows: find the larger $C>0$ such that 
\begin{equation*}
\mathbb{P}\left[ P\left( Ks_{n}-Ks_{M}\right) \leq C\right] \leq n^{-\alpha }%
\text{ .}
\end{equation*}%
We then have, by definition of the M-estimator $s_{n}$,%
\begin{eqnarray}
&&\mathbb{P}\left[ P\left( Ks_{n}-Ks_{M}\right) \leq C\right]  \notag \\
&\leq &\mathbb{P}\left[ \inf_{s\in M_{C}}P_{n}\left( Ks-Ks_{M}\right) \leq
\inf_{s\in M_{>C}}P_{n}\left( Ks-Ks_{M}\right) \right]  \notag \\
&=&\mathbb{P}\left[ \sup_{s\in M_{C}}P_{n}\left( Ks_{M}-Ks\right) \geq
\sup_{s\in M_{>C}}P_{n}\left( Ks_{M}-Ks\right) \right] \text{ .}
\label{reecriture_2}
\end{eqnarray}%
Expressions obtained in (\ref{reecriture_1}) and (\ref{reecriture_2}) allow
to reduce both upper and lower bounds problems for the excess risk to the
comparison of two quantities of interest,%
\begin{equation*}
\sup_{s\in M_{C}}P_{n}\left( Ks_{M}-Ks\right) \text{ \ and \ }\sup_{s\in
M_{>C}}P_{n}\left( Ks_{M}-Ks\right) \text{ .}
\end{equation*}%
Moreover, by setting $\mathcal{D}_{L}=\left\{ s\in M\text{ };\text{ }P\left(
Ks_{n}-Ks_{M}\right) =L\right\} $, we get%
\begin{eqnarray}
\sup_{s\in M_{C}}P_{n}\left( Ks_{M}-Ks\right) &=&\sup_{0\leq L\leq C}\left\{
\sup_{s\in \mathcal{D}_{L}}P_{n}\left( Ks_{M}-Ks\right) \right\}  \notag \\
&=&\sup_{0\leq L\leq C}\left\{ \sup_{s\in \mathcal{D}_{L}}\left\{ \left(
P_{n}-P\right) \left( Ks_{M}-Ks\right) +P\left( Ks_{M}-Ks\right) \right\}
\right\}  \notag \\
&=&\sup_{0\leq L\leq C}\left\{ \sup_{s\in \mathcal{D}_{L}}\left\{ \left(
P_{n}-P\right) \left( Ks_{M}-Ks\right) \right\} -L\right\}
\label{decoupe_sphere_M_C}
\end{eqnarray}%
and also 
\begin{equation}
\sup_{s\in M_{>C}}P_{n}\left( Ks_{M}-Ks\right) =\sup_{L>C}\left\{ \sup_{s\in 
\mathcal{D}_{L}}\left\{ \left( P_{n}-P\right) \left( Ks_{M}-Ks\right)
\right\} -L\right\} \text{ .}  \label{decoupe_sphere_M>C}
\end{equation}%
The study of the excess risk thus reduces to the control of the following
suprema, on the spheres $\mathcal{D}_{L}$ of radius $L$ in $\left(
M,\left\Vert 
\cdot%
\right\Vert _{2}\right) $, of the empirical process indexed by contrasted
increments of functions in $M$ around the projection $s_{M}$ of the target,%
\begin{equation}
\sup_{s\in \mathcal{D}_{L}}\left\{ \left( P-P_{n}\right) \left(
Ks-Ks_{M}\right) \right\} \text{ , }L\geq 0\text{ .}  \label{slide_interet}
\end{equation}%
Similarly, the empirical excess risk can be written, by definition of the
M-estimator $s_{n}$,%
\begin{eqnarray}
P_{n}\left( Ks_{M}-Ks_{n}\right) &=&\sup_{s\in M}P_{n}\left( Ks_{M}-Ks\right)
\notag \\
&=&\sup_{L\geq 0}\left\{ \sup_{s\in \mathcal{D}_{L}}P_{n}\left(
Ks_{M}-Ks\right) \right\}  \notag \\
&=&\sup_{L\geq 0}\left\{ \sup_{s\in \mathcal{D}_{L}}\left\{ \left(
P_{n}-P\right) \left( Ks_{M}-Ks\right) \right\} -L\right\} \text{ .}
\label{formule_exces_empirique}
\end{eqnarray}%
Hence, the study of the empirical excess risk reduces again to the control
of the quantities given in (\ref{slide_interet}). As these quantities are
(local) suprema of an empirical process, we can handle, under the right
hypotheses, the deviations from their mean \textit{via} the use of
concentration inequalities - deviations from the right being described with
optimal constants by Bousquet inequality (Bousquet, \cite{Bousquet:02},
recalled in Section \ref{section_probabilistic_tools} at the end of the
present paper) and deviations from left being controlled with sharp
constants by Klein and Rio inequality (Klein and Rio \cite{Klein_Rio:05},
also recalled in Section \ref{section_probabilistic_tools}). We can thus
expect that, under standard assumptions, the deviations are negligible
compared to the means with large enough probability, at least for radii $L$
not too small,%
\begin{equation}
\sup_{s\in \mathcal{D}_{L}}\left\{ \left( P-P_{n}\right) \left(
Ks-Ks_{M}\right) \right\} \sim \mathbb{E}\left[ \sup_{s\in \mathcal{D}%
_{L}}\left\{ \left( P-P_{n}\right) \left( Ks-Ks_{M}\right) \right\} \right] 
\text{ .}  \label{esperance}
\end{equation}

\bigskip

\noindent \textbf{Remark 1} \textit{It is worth noting that the above
computations, which allow to investigate both upper and lower bound
problems, only rely on the definition of }$s_{n}$\textit{\ as a minimizer of
the empirical risk over the model }$M$\textit{, and not on the particular
structure of the least squares contrast. Thus, formula (\ref{reecriture_1}),
(\ref{reecriture_2}), (\ref{decoupe_sphere_M_C}), (\ref{decoupe_sphere_M>C})
and (\ref{formule_exces_empirique}) are general facts of M-estimation -
whenever the projection }$s_{M}$\textit{\ of the target onto the model }$M$%
\textit{\ exists. Moreover, although presented in a quite different manner,
our computations related to the control of the true excess risk are in
essence very similar to those developed by Bartlett and Mendelson in \cite%
{BarMen:06}, concerning what they call "a direct analysis of the empirical
minimization algorithm". Indeed, the authors highlight in Section 3 of \cite%
{BarMen:06} that, under rather mild hypotheses, the true excess risk is
essentially the maximizer of the function }$\mathcal{V}_{n}\left( L\right)
-L $\textit{, where we set}%
\begin{equation*}
\mathcal{V}_{n}\left( L\right) :=\mathbb{E}\left[ \sup_{s\in \mathcal{D}%
_{L}}\left\{ \left( P-P_{n}\right) \left( Ks-Ks_{M}\right) \right\} \right] 
\text{ .}
\end{equation*}%
\textit{Now, combining (\ref{reecriture_1}), (\ref{reecriture_2}), (\ref%
{decoupe_sphere_M_C}) and (\ref{decoupe_sphere_M>C}), it is easily seen that
in the case where }$s_{n}$\textit{\ is unique and where}%
\begin{equation*}
\forall C\geq 0,\text{ }\sup_{s\in \mathcal{D}_{C}}P_{n}\left(
Ks_{M}-Ks\right) \text{ \textit{is achieved} }\left( =\max_{s\in \mathcal{D}%
_{C}}P_{n}\left( Ks_{M}-Ks\right) \right) \text{ ,}
\end{equation*}%
\textit{\ we have in fact the following exact formula,}%
\begin{eqnarray}
P\left( Ks_{n}-Ks_{M}\right) &=&\arg \max_{L\geq 0}\left\{ \max_{s\in 
\mathcal{D}_{L}}P_{n}\left( Ks_{M}-Ks\right) \right\}  \notag \\
&=&\arg \max_{L\geq 0}\left\{ \max_{s\in \mathcal{D}_{L}}\left(
P-P_{n}\right) \left( Ks-Ks_{M}\right) -L\right\} \text{ .}
\label{formule_theo_intro_2}
\end{eqnarray}%
\textit{So, if (\ref{esperance}) is satisfied with high probability, we
recover Bartlett and Mendelson's observation, which is}%
\begin{equation}
P\left( Ks_{n}-Ks_{M}\right) \sim \arg \max_{L\geq 0}\left\{ \mathcal{V}%
_{n}\left( L\right) -L\right\} \text{ .}  \label{formule_BM}
\end{equation}%
\textit{In Theorem 3.1 of \cite{BarMen:06}, a precise sense is given to (\ref%
{formule_BM}), in a rather general framework. In particular, a lower bound
for the excess risk is given but only through an additional condition
controlling the supremum of the empirical process of interest itself over a
subset of functions of \textquotedblleft small\textquotedblright\ excess
risks. This additional condition remains the major restriction concerning
the related result of Bartlett and Mendelson. In the following, we show in
our more restricted framework how to take advantage of the linearity of the
model, as well as the existence of an expansion of the least squares
contrast around the projection }$s_{M}$\textit{\ of the} \textit{target, to
derive lower bounds without additional assumptions on the behavior of the
empirical process of interest. Moreover, our methodology allow to explicitly
calculate the first order of the quantity given at the right side of (\ref%
{formule_BM}), thus exhibiting a rather simple complexity term controlling
the rate of convergence of the excess risk in the regression setting and
relating some geometrical characteristics of the model }$M$\textit{\ to the
unknown law }$P$\textit{\ of data.\bigskip }

\noindent \textbf{Remark 2} \textit{Formula (\ref{formule_exces_empirique})
and (\ref{formule_theo_intro_2}) above show that the true and empirical
excess risks are of different nature, in the sense that the first one is
referred to the arguments\ of the function }%
\begin{equation*}
\Gamma _{n}:L\left( \geq 0\right) \mapsto \max_{s\in \mathcal{D}_{L}}\left(
P-P_{n}\right) \left( Ks-Ks_{M}\right) -L\text{ ,}
\end{equation*}%
\textit{whereas the second one is measured from the values of the function }$%
\Gamma _{n}$\textit{. Hence, the equivalence between the true and the
empirical excess risks, when satisfied, is in general not straightforward.
It is a consequence of the following \textquotedblleft fixed point
type\textquotedblright\ equation,}%
\begin{equation*}
\arg \max_{\mathbb{R}_{+}}\left\{ \Gamma _{n}\right\} \sim \max_{\mathbb{R}%
_{+}}\left\{ \Gamma _{n}\right\} \text{ .}
\end{equation*}

\bigskip

Considering that the approximation stated in (\ref{esperance}) is suitably
satisfied, it remains to get an asymptotic first order expansion of its
right-hand term. Such a control is obtained through the use of the least
squares contrast expansion given in (\ref{dev_contrast}). Indeed, using (\ref%
{dev_contrast}), we get%
\begin{eqnarray}
&&\mathbb{E}\left[ \sup_{s\in \mathcal{D}_{L}}\left\{ \left( P-P_{n}\right)
\left( Ks-Ks_{M}\right) \right\} \right]  \notag \\
&=&\underset{\text{principal part}}{\underbrace{\mathbb{E}\left[ \sup_{s\in 
\mathcal{D}_{L}}\left\{ \left( P-P_{n}\right) \left( \psi _{1,M}%
\cdot%
\left( s-s_{M}\right) \right) \right\} \right] }}+\underset{\text{residual
term}}{\underbrace{\mathbb{E}\left[ \sup_{s\in \mathcal{D}_{L}}\left\{
\left( P-P_{n}\right) \left( \left( s-s_{M}\right) ^{2}\right) \right\} %
\right] }}\text{ .}  \label{decompo_outlines}
\end{eqnarray}%
In order to show that the residual term is negligible compared with the
principal part, it is natural to use a contraction principle (see Theorem
4.12\ of \cite{LedouxTal:91}, also recalled in Section \ref%
{section_probabilistic_tools}). Indeed, arguments of the empirical process
appearing in the residual term are related to the square of the arguments
defining the empirical process in the principal part. Moreover, it appears
by using the contraction principle, that the ratio of the residual term over
the principal part is roughly given by the supremum norm of the indexes: $%
\sup_{s\in \mathcal{D}_{L}}\left\vert \left( s-s_{M}\right) \left( x\right)
\right\vert $ (see Lemma \ref{lower_rademacher copy(1)} in Section \ref%
{useful_lemmas} for more details). Now, using assumption (\textbf{H3}) of
Section \ref{section_main_assump_reg_fixed}, concerning the unit envelope of
the linear model $M$, we get that the last quantity is of order $\sqrt{DL}$.
Since the values $L$ of interest are typically of order $D/n$, the quantity
controlling the ratio is not sharp enough as it does not converge to zero as
soon as the dimension $D$ is of order at least $\sqrt{n}$.

We thus have to refine our analysis in order to be able to neglect the
residual term. The assumption of sup-norm consistency, of the least squares
estimator $s_{n}$ toward the projection $s_{M}$ of the target onto the model 
$M$, appears here to be essential. Indeed, if assumption (\textbf{H5}) of
Section \ref{section_main_assump_reg_fixed} is satisfied, then all the above
computations can be restricted with high probability to the subset where
belongs the estimator $s_{n}$, this subset being more precisely%
\begin{equation}
B_{L_{\infty }}\left( s_{M},R_{n,D,\alpha }\right) =\left\{ s\in M\text{ };%
\text{ }\left\Vert s-s_{M}\right\Vert _{\infty }\leq R_{n,D,\alpha }\right\}
\subset M\text{ ,}  \label{ball}
\end{equation}
$R_{n,D,\alpha }\ll 1$ being the rate of convergence in sup-norm of $s_{n}$
toward $s_{M}$, defined in (\textbf{H5}). In particular, the spheres of
interest $\mathcal{D}_{L}$ are now replaced in the calculations by their
intersection $\mathcal{\tilde{D}}_{L}$ with the ball of radius $%
R_{n,D,\alpha }$ in $\left( M,\left\Vert 
\cdot%
\right\Vert _{\infty }\right) $,%
\begin{equation*}
\mathcal{\tilde{D}}_{L}=\mathcal{D}_{L}\cap B_{L_{\infty }}\left(
s_{M},R_{M,n,\alpha }\right) \text{ .}
\end{equation*}%
The ratio between the consequently modified residual term and principal part
of (\ref{decompo_outlines}) is then roughly controlled by $R_{n,D,\alpha }$
(see again Lemma \ref{lower_rademacher copy(1)} in Section \ref%
{useful_lemmas}), a quantity indeed converging to zero as desired. Hence,
under the assumption (\textbf{H5}), we get 
\begin{equation}
\mathbb{E}\left[ \sup_{s\in \mathcal{\tilde{D}}_{L}}\left\{ \left(
P-P_{n}\right) \left( Ks-Ks_{M}\right) \right\} \right] \sim \mathbb{E}\left[
\sup_{s\in \mathcal{\tilde{D}}_{L}}\left\{ \left( P-P_{n}\right) \left( \psi
_{1,M}%
\cdot%
\left( s-s_{M}\right) \right) \right\} \right] \text{ .}
\label{partie_principale_tilde}
\end{equation}%
A legitimate and important question is: how restrictive is assumption (%
\textbf{H5}) of consistency in sup-norm of the least squares estimator ? We
prove in Lemma \ref{consistency_histograms} of Section \ref%
{section_histogram_case_fixed} that this assumption is satisfied for models
of histograms defined on a partition satisfying some regularity condition,
at\ a rate of convergence of order $\sqrt{D\ln \left( n\right) /n}$.
Moreover, in Lemma \ref{lemma_piecewise_polynomials_consistency}, Section %
\ref{section_piecewise_case_fixed}, we extend this result for models of
piecewise polynomials uniformly bounded in their degrees, again under some
lower-regularity assumption on the partition defining the model; the rate of
convergence being also preserved. A systematical study of consistency in
sup-norm of least squares estimators, on more general finite-dimensional
linear models, is also postponed to a forthcoming paper.

The control of the right-hand side of (\ref{partie_principale_tilde}), which
is needed to be sharp, is particularly technical, and is essentially
contained in Lemmas \ref{bound_moment_order_1}\ and \ref%
{control_lower_infinity}\ of Section \ref{useful_lemmas}.\ Let us shortly
describe the mathematical figures underlying this control. First, by
bounding the variance of the considered supremum of the empirical process -
by using a result due to Ledoux \cite{Ledoux:96}, see Theorem \ref%
{Theorem_Ledoux} and also Corollary \ref{cor_alt_hoff_2}\ in Section \ref%
{section_probabilistic_tools} -, we roughly get, for values of $L$ of
interest,%
\begin{equation}
\mathbb{E}\left[ \sup_{s\in \mathcal{\tilde{D}}_{L}}\left\{ \left(
P-P_{n}\right) \left( \psi _{1,M}%
\cdot%
\left( s-s_{M}\right) \right) \right\} \right] \sim \mathbb{E}^{1/2}\left[
\left( \sup_{s\in \mathcal{\tilde{D}}_{L}}\left\{ \left( P-P_{n}\right)
\left( \psi _{1,M}%
\cdot%
\left( s-s_{M}\right) \right) \right\} \right) ^{2}\right] \text{ .}
\label{mom1_mom2}
\end{equation}%
Then, by assuming that the model $M$ is fulfilled with a localized
orthonormal basis, as stated in assumption (\textbf{H4}) of Section \ref%
{section_main_assump_reg_fixed}, it can be shown that the localization on
the ball $B_{L_{\infty }}\left( s_{M},R_{M,n,\alpha }\right) $ can be
removed from the right-hand side of (\ref{mom1_mom2}), in the sense that 
\begin{equation}
\mathbb{E}^{1/2}\left[ \left( \sup_{s\in \mathcal{\tilde{D}}_{L}}\left\{
\left( P-P_{n}\right) \left( \psi _{1,M}%
\cdot%
\left( s-s_{M}\right) \right) \right\} \right) ^{2}\right] \sim \mathbb{E}%
^{1/2}\left[ \left( \sup_{s\in \mathcal{D}_{L}}\left\{ \left( P-P_{n}\right)
\left( \psi _{1,M}%
\cdot%
\left( s-s_{M}\right) \right) \right\} \right) ^{2}\right] \text{ .}
\label{without_loc}
\end{equation}%
The property of localized basis is standard in model selection theory (see
for instance Chapter 7 of \cite{Massart:07}) and was first introduced by Birg%
\'{e} and Massart in \cite{BirgeMassart:98}, also for deriving sharp
exponential bounds in a M-estimation context. We show in Lemmas \ref%
{lemma_histograms_localized}\ and \ref{lemma_piecewise_polynomials_localized}%
\ that this assumption is satisfied for models of histograms and piecewise
polynomials respectively, when they satisfy a certain regularity assumption
concerning the underlying partition.

Finally, as $\mathcal{D}_{L}$ is a sphere in $\left( M,\left\Vert 
\cdot%
\right\Vert _{2}\right) $, we simply get, by the use of Cauchy-Schwarz
inequality, that the right-hand side of (\ref{without_loc}) is equal to $%
\sqrt{\left( L/n\right) 
\cdot%
\sum_{k=1}^{D}%
\var%
\left( \psi _{1,M}%
\cdot%
\varphi _{k}\right) }$, where $\left( \varphi _{k}\right) _{k=1}^{D}$ is an
orthonormal basis of $M$. Gathering our arguments, we then obtain%
\begin{eqnarray}
P\left( Ks_{n}-Ks_{M}\right) &\sim &\arg \max_{L\geq 0}\left\{ \sup_{s\in 
\mathcal{D}_{L}}\mathbb{E}\left[ \left( P_{n}-P\right) \left(
Ks_{M}-Ks\right) \right] -L\right\}  \notag \\
&\sim &\arg \max_{L\geq 0}\left\{ \sqrt{\frac{L%
\cdot%
\sum_{k=1}^{D_{M}}%
\var%
\left( \psi _{1,M}%
\cdot%
\varphi _{k}\right) }{n}}-L\right\} =\frac{1}{4}\frac{D_{M}}{n}\mathcal{K}%
_{1,M}^{2}\text{ ,}  \label{bound_true}
\end{eqnarray}%
where $\mathcal{K}_{1,M}^{2}:=\frac{1}{D_{M}}\sum_{k=1}^{D_{M}}%
\var%
\left( \psi _{1,M}%
\cdot%
\varphi _{k}\right) $. As shown in Section \ref{section_complexity_model}\
below, the (normalized) complexity term $\mathcal{K}_{1,M}$ is independent
of the choice of the basis $\left( \varphi _{k}\right) _{k=1}^{D}$ and is,
under our assumptions, of the order of a constant. Concerning the empirical
excess risk, we have%
\begin{eqnarray}
P_{n}\left( Ks_{M}-Ks_{n}\right) &=&\max_{L\geq 0}\left\{ \sup_{s\in 
\mathcal{D}_{L}}\mathbb{E}\left[ \left( P_{n}-P\right) \left(
Ks_{M}-Ks\right) \right] -L\right\}  \notag \\
&\sim &\max_{L\geq 0}\left\{ \sqrt{\frac{L%
\cdot%
\sum_{k=1}^{D_{M}}%
\var%
\left( \psi _{1,M}%
\cdot%
\varphi _{k}\right) }{n}}-L\right\} =\frac{1}{4}\frac{D_{M}}{n}\mathcal{K}%
_{1,M}^{2}\text{ .}  \label{bound_emp}
\end{eqnarray}%
In particular, the equivalence%
\begin{equation*}
P\left( Ks_{n}-Ks_{M}\right) \sim P_{n}\left( Ks_{M}-Ks_{n}\right) \left(
\sim \frac{1}{4}\frac{D_{M}}{n}\mathcal{K}_{1,M}^{2}\right)
\end{equation*}%
is justified.

In Theorem \ref{principal_reg}\ below, a precise, non-asymptotic sense, is
given to equivalences described in (\ref{bound_true}) and (\ref{bound_emp}).
This is done under the structural constraints stated in conditions (\textbf{%
H4}) and (\textbf{H5}), for models of reasonable dimension. Moreover, we
give in Theorem \ref{theorem_upper_true_small_models} upper bounds for the
true and empirical excess risks, that are less precise than the bounds of
Theorem \ref{principal_reg}, but that are also valid for\ models of small
dimension. Corollaries of these theorems are given in the case of histograms
and piecewise polynomials, in Corollaries \ref%
{corollary_upper_lower_histograms} and \ref%
{corollary_upper_lower_piecewise_ploynomials}\ respectively. Indeed, we show
that in these particular cases, our general conditions (\textbf{H4}) and (%
\textbf{H5}) essentially reduce to a simple lower-regularity assumption on
the underlying partition.

\section{True and empirical excess risk bounds\label{section_reg_fixed_model}%
}

In this section, we derive under general constraints on the linear model $M$%
, upper and lower bounds for the true and empirical excess risk, that are
optimal - and equal - at the first order. In particular, we show that the
true excess risk is equivalent to the empirical one when the model is of
reasonable dimension. For smaller dimensions, we only achieve some upper
bounds.

\subsection{Main assumptions\label{section_main_assump_reg_fixed}}

We turn now to the statement of some assumptions that will be needed to
derive our results in Section \ref{section_theorems_fixed model copy(1)}.
These assumptions will be further discussed in Section \ref%
{section_complexity_model}.

\bigskip

\noindent \textbf{Boundedness assumptions:}

\bigskip

\begin{itemize}
\item (\textbf{H1}) The data and the linear projection of the target onto $M$
are bounded: a positive finite constant $A$ exists such that 
\begin{equation}
\left\vert Y_{i}\right\vert \leq A\text{ }a.s.  \label{bounded_data}
\end{equation}%
and%
\begin{equation}
\left\Vert s_{M}\right\Vert _{\infty }\leq A\text{ }.
\label{bounded_projection}
\end{equation}
\end{itemize}

\noindent Hence, from (\textbf{H1}) we deduce that 
\begin{equation}
\left\Vert s_{\ast}\right\Vert _{\infty}=\left\Vert \mathbb{E}\left[
Y\left\vert X=%
\cdot%
\right. \right] \right\Vert _{\infty}\leq A  \label{bounded_target}
\end{equation}
and that there exists a constant $\sigma_{\max}>0$ such that 
\begin{equation}
\sigma^{2}\left( X_{i}\right) \leq\sigma_{\max}^{2}\leq A^{2}\text{ \ }a.s.
\label{majo_sigma_max}
\end{equation}
Moreover, as $\psi_{1,M}\left( z\right) =-2\left( y-s_{M}\left( x\right)
\right) $ for all $z=\left( x,y\right) \in\mathcal{Z}$, we also deduce that%
\begin{equation}
\left\vert \psi_{1,M}\left( X_{i},Y_{i}\right) \right\vert \leq4A\text{ \ }%
a.s.  \label{control_phi_1M}
\end{equation}

\begin{itemize}
\item (\textbf{H2}) The heteroscedastic noise level $\sigma $ is uniformly
bounded from below: a positive finite constant $\sigma _{\min }$ exists such
that%
\begin{equation*}
0<\sigma _{\min }\leq \sigma \left( X_{i}\right) \text{ \ }a.s.
\end{equation*}
\end{itemize}

\bigskip

\noindent \textbf{Models with localized basis in }$L_{2}\left( P^{X}\right) $%
\textbf{:}

\bigskip

\noindent Let us define a function $\Psi_{M}$ on $\mathcal{X}$, that we call
the unit envelope of $M$, such that%
\begin{equation}
\Psi_{M}\left( x\right) =\frac{1}{\sqrt{D}}\sup_{s\in M,\left\Vert
s\right\Vert _{2}\leq1}\left\vert s\left( x\right) \right\vert \text{ }.
\label{def_enveloppe}
\end{equation}
As $M$ is a finite dimensional real vector space, the supremum in (\ref%
{def_enveloppe}) can also be taken over a countable subset of $M$, so $%
\Psi_{M}$ is a measurable function.

\begin{itemize}
\item (\textbf{H3}) The unit envelope of $M$ is uniformly bounded on $%
\mathcal{X}$: a positive constant $A_{3,M}$ exists such that 
\begin{equation*}
\left\Vert \Psi _{M}\right\Vert _{\infty }\leq A_{3,M}<\infty \text{ }.
\end{equation*}
\end{itemize}

\noindent The following assumption is stronger than (\textbf{H3}).

\begin{itemize}
\item (\textbf{H4}) Existence of a localized basis in $\left( M,\left\Vert 
\cdot%
\right\Vert _{2}\right) $: there exists an orthonormal basis $\varphi
=\left( \varphi _{k}\right) _{k=1}^{D}$ in $\left( M,\left\Vert 
\cdot%
\right\Vert _{2}\right) $ that satisfies, for a positive constant $%
r_{M}\left( \varphi \right) $ and all $\beta =\left( \beta _{k}\right)
_{k=1}^{D}\in \mathbb{R}^{D}$,%
\begin{equation*}
\left\Vert \sum_{k=1}^{D}\beta _{k}\varphi _{k}\right\Vert _{\infty }\leq
r_{M}\left( \varphi \right) \sqrt{D}\left\vert \beta \right\vert _{\infty }%
\text{ ,}
\end{equation*}%
where $\left\vert \beta \right\vert _{\infty }=\max \left\{ \left\vert \beta
_{k}\right\vert ;k\in \left\{ 1,...,D\right\} \right\} $ is the sup-norm of
the $D$-dimensional vector $\beta .$
\end{itemize}

\noindent \textbf{Remark 3} \textit{(\textbf{H4}) implies (\textbf{H3}) and
in that case }$A_{3,M}=r_{M}\left( \varphi \right) $\textit{\ is convenient.}

\bigskip

\noindent \textbf{The assumption of consistency in sup-norm:}

\bigskip

\noindent In order to handle second order terms in the expansion of the
contrast (\ref{dev_contrast}), we assume that the least squares estimator is
consistent for the sup-norm on the space $\mathcal{X}$. More precisely, this
requirement can be stated as follows.

\begin{itemize}
\item (\textbf{H5}) Assumption of consistency in sup-norm\textbf{: }for any $%
A_{+}>0$, if $M$ is a model of dimension $D$ satisfying%
\begin{equation*}
D\leq A_{+}\frac{n}{\left( \ln n\right) ^{2}}\text{ },
\end{equation*}%
then for every $\alpha >0$, we can find a positive integer $n_{1}$ and a
positive constant $A_{cons}$ satisfying the following property: there exists 
$R_{n,D,\alpha }>0$ depending on $D,$ $n$ and $\alpha $, such that%
\begin{equation}
R_{n,D,\alpha }\leq \frac{A_{cons}}{\sqrt{\ln n}}  \label{def_A_cons}
\end{equation}%
and by setting 
\begin{equation}
\Omega _{\infty ,\alpha }=\left\{ \left\Vert s_{n}-s_{M}\right\Vert _{\infty
}\leq R_{n,D,\alpha }\right\} \text{ },  \label{def_omega_infini}
\end{equation}%
it holds for all $n\geq n_{1}$, 
\begin{equation}
\mathbb{P}\left[ \Omega _{\infty ,\alpha }\right] \geq 1-n^{-\alpha }\text{ }%
.  \label{proba_omega_infini}
\end{equation}
\end{itemize}

\subsection{Theorems\label{section_theorems_fixed model copy(1)}}

We state here the general results of this article, that will be applied in
Section \ref{section_histogram_case_fixed} and \ref%
{section_piecewise_case_fixed} in the case of piecewise constant functions
and piecewise polynomials respectively.

\begin{theorem}
\label{principal_reg}Let $A_{+},A_{-},\alpha >0$ and let $M$ be a linear
model of finite dimension $D$. Assume that (\textbf{H1}), (\textbf{H2}), (%
\textbf{H4}) and (\textbf{H5}) hold and take $\varphi =\left( \varphi
_{k}\right) _{k=1}^{D}$ an orthonormal basis of $\left( M,\left\Vert 
\cdot%
\right\Vert _{2}\right) $ satisfying (\textbf{H4)}. If it holds 
\begin{equation}
A_{-}\left( \ln n\right) ^{2}\leq D\leq A_{+}\frac{n}{\left( \ln n\right)
^{2}}\text{ },  \label{hypo_dim_reg}
\end{equation}%
then a positive finite constant $A_{0}$ exists, only depending on $\alpha
,A_{-}$ and on the constants $A,\sigma _{\min },r_{M}\left( \varphi \right) $
defined in assumptions (\textbf{H1), (H2) }and \textbf{(H4) }respectively,
such that by setting%
\begin{equation}
\varepsilon _{n}=A_{0}\max \left\{ \left( \frac{\ln n}{D}\right) ^{1/4},%
\text{ }\left( \frac{D\ln n}{n}\right) ^{1/4},\text{ }\sqrt{R_{n,D,\alpha }}%
\right\} \text{ },  \label{def_A_0}
\end{equation}%
we have for all $n\geq n_{0}\left( A_{-},A_{+},A,A_{cons},r_{M}\left(
\varphi \right) ,\sigma _{\min },n_{1},\alpha \right) $,%
\begin{align}
\mathbb{P}\left[ P\left( Ks_{n}-Ks_{M}\right) \geq \left( 1-\varepsilon
_{n}\right) \frac{1}{4}\frac{D}{n}\mathcal{K}_{1,M}^{2}\right] & \geq
1-5n^{-\alpha }\text{ },  \label{lower_true_risk} \\
\mathbb{P}\left[ P\left( Ks_{n}-Ks_{M}\right) \leq \left( 1+\varepsilon
_{n}\right) \frac{1}{4}\frac{D}{n}\mathcal{K}_{1,M}^{2}\right] & \geq
1-5n^{-\alpha }\text{ },  \label{upper_true_risk} \\
\mathbb{P}\left[ P_{n}\left( Ks_{M}-Ks_{n}\right) \geq \left( 1-\varepsilon
_{n}^{2}\right) \frac{1}{4}\frac{D}{n}\mathcal{K}_{1,M}^{2}\right] & \geq
1-2n^{-\alpha }\text{ },  \label{lower_emp_risk} \\
\mathbb{P}\left[ P_{n}\left( Ks_{M}-Ks_{n}\right) \leq \left( 1+\varepsilon
_{n}^{2}\right) \frac{1}{4}\frac{D}{n}\mathcal{K}_{1,M}^{2}\right] & \geq
1-3n^{-\alpha }\text{ },  \label{upper_emp_risk}
\end{align}%
where $\mathcal{K}_{1,M}^{2}=\frac{1}{D}\sum_{k=1}^{D}%
\var%
\left( \psi _{1,M}%
\cdot%
\varphi _{k}\right) $. In addition, when (\textbf{H5}) does not hold, but (%
\textbf{H1}), (\textbf{H2}) and (\textbf{H4}) are satisfied, we still have
for all $n\geq n_{0}\left( A_{-},A_{+},A,r_{M}\left( \varphi \right) ,\sigma
_{\min },\alpha \right) $,%
\begin{equation}
\mathbb{P}\left( P_{n}\left( Ks_{M}-Ks_{n}\right) \geq \left( 1-A_{0}\max
\left\{ \sqrt{\frac{\ln n}{D}},\sqrt{\frac{D\ln n}{n}}\right\} \right) \frac{%
D}{4n}\mathcal{K}_{1,M}^{2}\right) \geq 1-2n^{-\alpha }\text{ }.
\label{lower_emp_risk_general}
\end{equation}
\end{theorem}

\noindent In Theorem \ref{principal_reg} above, we achieve sharp upper and
lower bounds for the true and empirical excess risks on $M$. They are
optimal at the first order since the leading constants are equal for upper
and lower bounds. Moreover, Theorem \ref{principal_reg} states the
equivalence with high probability of the true and empirical excess risks for
models of reasonable dimensions. We notice that second orders are smaller
for the empirical excess risk than for the true one. Indeed, when normalized
by the first order, the deviations of the empirical excess risk are square
of the deviations of the true one. Our bounds also give another evidence of
the concentration phenomenon of the empirical excess risk exhibited by
Boucheron and Massart \cite{BouMas:10} in the slightly different context of
M-estimation with bounded contrast where some margin condition hold. Notice
that considering the lower bound of the empirical excess risk given in (\ref%
{lower_emp_risk_general}), we do not need to assume the consistency of the
least squares estimator $s_{n}$ towards the linear projection $s_{M}$.

\noindent We turn now to upper bounds in probability for the true and
empirical excess risks on models with possibly small dimensions. In this
context, we do not achieve sharp or explicit constants in the rates of
convergence.

\begin{theorem}
\label{theorem_upper_true_small_models}Let $\alpha ,A_{+}>0$ be fixed and
let $M$ be a linear model of finite dimension 
\begin{equation*}
1\leq D\leq A_{+}\frac{n}{\left( \ln n\right) ^{2}}\text{ .}
\end{equation*}%
Assume that assumptions (\textbf{H1}), (\textbf{H3}) and (\textbf{H5}) hold.
Then a positive constant $A_{u}$ exists, only depending on $%
A,A_{cons},A_{3,M}$ and $\alpha $, such that for all $n\geq n_{0}\left(
A_{cons},n_{1}\right) $,%
\begin{equation}
\mathbb{P}\left[ P\left( Ks_{n}-Ks_{M}\right) \geq A_{u}\frac{D\vee \ln n}{n}%
\right] \leq 3n^{-\alpha }  \label{upper_true_small}
\end{equation}%
and%
\begin{equation}
\mathbb{P}\left[ P_{n}\left( Ks_{M}-Ks_{n}\right) \geq A_{u}\frac{D\vee \ln n%
}{n}\right] \leq 3n^{-\alpha }\text{ }.  \label{upper_emp_small}
\end{equation}
\end{theorem}

\noindent Notice that on contrary to the situation of Theorem \ref%
{principal_reg}, we do not assume that (\textbf{H2}) hold. This assumption
states that the noise level is uniformly bounded away from zero over the
space $\mathcal{X}$, and allows in Theorem \ref{principal_reg} to derive
lower bounds for the true and empirical excess risks, as well as to achieve
sharp constants in the deviation bounds for models of reasonable dimensions.
In Theorem \ref{theorem_upper_true_small_models}, we just derive upper
bounds and assumption (\textbf{H2}) is not needed. The price to pay is that
constants in the rates of convergence derived in (\ref{upper_true_small})
and (\ref{upper_emp_small}) are possibly larger than the corresponding ones
of Theorem \ref{principal_reg}, but our results still hold true for small
models. Moreover, in the case of models with reasonable dimensions, that is
dimensions satisfying assumption (\ref{hypo_dim_reg}) of Theorem \ref%
{principal_reg}, the rate of decay is preserved compared to Theorem \ref%
{principal_reg} and is proportional to $D/n$.

\noindent The proofs of the above theorems can be found in Section \ref%
{section_proof_fixed_model}.

\subsection{Some additional comments\label{section_complexity_model}}

Let us first comment on the assumptions given in Section \ref%
{section_main_assump_reg_fixed}. Assumptions (\ref{bounded_data}) and (%
\textbf{H2}) are rather mild and can also be found in the work of Arlot and
Massart \cite{ArlotMassart:09} related to the case of histograms, where they
are respectively denoted by (\textbf{Ab}) and (\textbf{An}). These
assumptions state respectively that the response variable $Y$ is uniformly
bounded and that the noise level is uniformly bounded away from zero. In 
\cite{ArlotMassart:09}, Arlot and Massart also notice that their results can
be extended to the unbounded case, where assumption (\textbf{Ab}) is
replaced by some condition on the moments of the noise, and where (\textbf{An%
}) is weakened into mild regularity conditions for the noise level. We
believe that moments conditions on the noise, in the spirit of assumptions
stated by Arlot and Massart, could also been taken into account in our study
in order to weaken (\ref{bounded_data}), but at the prize of many technical
efforts that are beyond the scope of the present paper. However, we explain
at the end of this section how condition (\textbf{H2}) can be relaxed - see
hypothesis (\textbf{H2bis}) below.

\noindent In assumption (\textbf{H4}) we require that the model $M$ is
provided with an orthonormal localized basis in $L_{2}\left( P^{X}\right) $.
This property is convenient when dealing with the $L_{\infty }$-structure on
the model, and this allows us to control the sup-norm of the functions in
the model by the sup-norm of the vector of their coordinates in the
localized basis. For examples of models with localized basis, and their use
in a model selection framework, we refer for instance to Section 7.4.2\ of
Massart \cite{Massart:07}, where it is shown that models of histograms,
piecewise polynomials and compactly supported wavelets are typical examples
of models with localized basis for the $L_{2}\left( 
\leb%
\right) $ structure, considering that $\mathcal{X\subset }\mathbb{R}^{k}$.
In Sections \ref{section_histogram_case_fixed} and \ref%
{section_piecewise_case_fixed}, we show that models of piecewise constant
and piecewise polynomials respectively can also have a localized basis for
the $L_{2}\left( P^{X}\right) $ structure, under rather mild assumptions on $%
P^{X}$. Assumption (\textbf{H4}) is needed in Theorem \ref{principal_reg},
whereas in Theorem \ref{theorem_upper_true_small_models} we only use the
weaker assumption (\textbf{H3}) on the unit envelope of the model $M$,
relating the $L_{2}$-structure of the model to the $L_{\infty }$-structure.
In fact, assumption (\textbf{H4}) allows us in the proof of Theorem \ref%
{principal_reg} to achieve sharp lower bounds for the quantities of
interest, whereas in Theorem \ref{theorem_upper_true_small_models} we only
give upper bounds in the case of small models.

\noindent We ask in assumption (\textbf{H5}) that the M-estimator is
consistent towards the linear projection $s_{M}$ of $s_{\ast}$ onto the
model $M$, at a rate at least better than $\left( \ln n\right) ^{-1/2}$ .
This can be considered as a rather strong assumption, but it is essential
for our methodology. Moreover, we show in Sections \ref%
{section_histogram_case_fixed} and \ref{section_piecewise_case_fixed} that
this assumption is satisfied under mild conditions for histogram models and
models of piecewise polynomials respectively, both at the rate 
\begin{equation*}
R_{n,D,\alpha}\propto\sqrt{\frac{D\ln n}{n}}\text{ }.
\end{equation*}

\noindent Secondly, let us comment on the rates of convergence given in
Theorem \ref{principal_reg} for models of reasonable dimensions. As we can
see in Theorem \ref{principal_reg}, the rate of estimation in a fixed model $%
M$ of reasonable dimension is determined at the first order by a key
quantity that relates the structure of the model to the unknown law $P$ of
data. We call this quantity the \textit{complexity} of the model\textbf{\ }$%
M $ and we denote it by $\mathcal{C}_{M}.$ More precisely, let us define%
\begin{equation*}
\mathcal{C}_{M}=\frac{1}{4}D\times \mathcal{K}_{1,M}^{2}
\end{equation*}%
where 
\begin{equation*}
\mathcal{K}_{1,M}=\sqrt{\frac{1}{D}\sum_{k=1}^{D}%
\var%
\left( \psi _{1,M}%
\cdot%
\varphi _{k}\right) }
\end{equation*}%
for a localized orthonormal basis $\left( \varphi _{k}\right) _{k=1}^{D}$ of 
$\left( M,\left\Vert 
\cdot%
\right\Vert _{2}\right) .$ Notice that $\mathcal{K}_{1,M}$ is well defined
as it does not depend on the choice of the basis $\left( \varphi _{k}\right)
_{k=1}^{D}.$ Indeed, since we have $P\left( \psi _{1,M}%
\cdot%
\varphi _{k}\right) =0$, we deduce that 
\begin{equation*}
\mathcal{K}_{1,M}^{2}=P\left( \psi _{1,M}^{2}%
\cdot%
\left( \frac{1}{D}\sum_{k=1}^{D}\varphi _{k}^{2}\right) \right) \text{ }.
\end{equation*}%
Now observe that, by using Cauchy-Schwarz inequality in Definition (\ref%
{def_enveloppe}), as pointed out by Birg\'{e} and Massart \cite%
{BirgeMassart:98}, we get%
\begin{equation}
\Psi _{M}^{2}=\frac{1}{D}\sum_{k=1}^{D}\varphi _{k}^{2}  \label{phi_M}
\end{equation}%
and so%
\begin{align}
\mathcal{K}_{1,M}^{2}& =P\left( \psi _{1,M}^{2}\Psi _{M}^{2}\right)  \notag
\\
& =4\mathbb{E}\left[ \mathbb{E}\left[ \left( Y-s_{M}\left( X\right) \right)
^{2}\left\vert X\right. \right] \Psi _{M}^{2}\left( X\right) \right]  \notag
\\
& =4\left( \mathbb{E}\left[ \sigma ^{2}\left( X\right) \Psi _{M}^{2}\left(
X\right) \right] +\mathbb{E}\left[ \left( s_{M}-s_{\ast }\right) ^{2}\left(
X\right) \Psi _{M}^{2}\left( X\right) \right] \right) \text{ }.
\label{formule_K1M_bis}
\end{align}%
On the one hand, if we assume (\textbf{H1}) then we obtain by elementary
computations%
\begin{equation}
\mathcal{K}_{1,M}\leq 2\sigma _{\max }+4A\leq 6A\text{ }.  \label{upper_K1M}
\end{equation}%
On the other hand, (\textbf{H2}) implies%
\begin{equation}
\mathcal{K}_{1,M}\geq 2\sigma _{\min }>0\text{ }.  \label{lower_K1M}
\end{equation}%
To fix ideas, let us explicitly compute $\mathcal{K}_{1,M}^{2}$ in a simple
case. Consider homoscedastic regression on a histogram model $M$, in which
the homoscedastic noise level $\sigma $ is such that 
\begin{equation*}
\sigma ^{2}\left( X\right) =\sigma ^{2}\text{ \ \ }a.s.\text{ ,}
\end{equation*}%
so we have%
\begin{equation*}
\mathbb{E}\left[ \sigma ^{2}\left( X\right) \Psi _{M}^{2}\left( X\right) %
\right] =\sigma ^{2}\mathbb{E}\left[ \Psi _{M}^{2}\left( X\right) \right]
=\sigma ^{2}\text{ }.
\end{equation*}%
Now, under notations of Lemma \ref{lemma_histograms_localized} below, 
\begin{equation*}
s_{M}=\sum_{I\in \mathcal{P}}\mathbb{E}\left[ Y\varphi _{I}\left( X\right) %
\right] \varphi _{I}=\sum_{I\in \mathcal{P}}\mathbb{E}\left[ Y\left\vert
X\in I\right. \right] \mathbf{1}_{I}\text{ },
\end{equation*}%
thus we deduce, by (\ref{phi_M}) and the previous equality, that%
\begin{align*}
\mathbb{E}\left[ \left( s_{M}-s_{\ast }\right) ^{2}\left( X\right) \Psi
_{M}^{2}\left( X\right) \right] & =\frac{1}{\left\vert \mathcal{P}%
\right\vert }\sum_{I\in \mathcal{P}}\mathbb{E}\left[ \left( s_{M}-s_{\ast
}\right) ^{2}\left( X\right) \varphi _{I}^{2}\left( X\right) \right] \\
& =\frac{1}{\left\vert \mathcal{P}\right\vert }\sum_{I\in \mathcal{P}}%
\mathbb{E}\left[ \left( \mathbb{E}\left[ Y\left\vert X\in I\right. \right] -%
\mathbb{E}\left[ Y\left\vert X\right. \right] \right) ^{2}\frac{\mathbf{1}%
_{X\in I}}{P^{X}\left( I\right) }\right] \\
& =\frac{1}{\left\vert \mathcal{P}\right\vert }\sum_{I\in \mathcal{P}}%
\mathbb{E}\left[ \left( \mathbb{E}\left[ Y\left\vert X\in I\right. \right] -%
\mathbb{E}\left[ Y\left\vert X\right. \right] \right) ^{2}\left\vert X\in
I\right. \right] \\
& =\frac{1}{\left\vert \mathcal{P}\right\vert }\sum_{I\in \mathcal{P}}%
\mathbb{V}\left[ \mathbb{E}\left[ Y\left\vert X\right. \right] \left\vert
X\in I\right. \right] \text{ },
\end{align*}%
where the conditional variance $\mathbb{V}\left[ U\left\vert \mathcal{A}%
\right. \right] $ of a variable $U$ with respect to the event $\mathcal{A}$
is defined to be 
\begin{equation*}
\mathbb{V}\left[ U\left\vert \mathcal{A}\right. \right] :=\mathbb{E}\left[
\left( U-\mathbb{E}\left[ U\left\vert \mathcal{A}\right. \right] \right)
^{2}\left\vert \mathcal{A}\right. \right] =\mathbb{E}\left[ U^{2}\left\vert 
\mathcal{A}\right. \right] -\left( \mathbb{E}\left[ U\left\vert \mathcal{A}%
\right. \right] \right) ^{2}\text{ .}
\end{equation*}%
By (\ref{formule_K1M_bis}), we explicitly get%
\begin{equation}
\mathcal{K}_{1,M}^{2}=4\left( \sigma ^{2}+\frac{1}{\left\vert \mathcal{P}%
\right\vert }\sum_{I\in \mathcal{P}}\mathbb{V}\left[ \mathbb{E}\left[
Y\left\vert X\right. \right] \left\vert X\in I\right. \right] \right) \text{ 
}.  \label{formule_K1M_histo}
\end{equation}%
A careful look at the proof of Theorem \ref{principal_reg} given in Section %
\ref{section_proof_fixed_model} show that condition (\textbf{H2}) is only
used through the lower bound (\ref{lower_K1M}), and thus (\textbf{H2}) can
be replaced by the following slightly more general assumption :

\begin{description}
\item[(\textbf{H2bis})] Lower bound on the normalized complexity $\mathcal{K}%
_{1,M}$ : a positive constant $A_{\min}$ exists such that%
\begin{equation*}
\mathcal{K}_{1,M}\geq A_{\min}>0\text{ .}
\end{equation*}
\end{description}

\noindent When (\textbf{H2}) holds, we see from Inequality \ref{lower_K1M}
that (\textbf{H2bis}) is satisfied with $A_{\min}=2\sigma_{\min}$. For
suitable models we can have for a positive constant $A_{\Psi}^{-}$ and for
all $x\in\mathcal{X}$, 
\begin{equation}
\Psi_{M}\left( x\right) \geq A_{\Psi}^{-}>0\text{ ,}  \label{lower_phi_M}
\end{equation}
and this allows to consider vanishing noise level, as we then have by (\ref%
{formule_K1M_bis}),%
\begin{equation*}
\mathcal{K}_{1,M}\geq2A_{\Psi}^{-}\sqrt{\mathbb{E}\left[ \sigma^{2}\left(
X\right) \right] }=2A_{\Psi}^{-}\left\Vert \sigma\right\Vert _{2}>0\text{ .}
\end{equation*}
As we will see in Sections \ref{section_histogram_case_fixed} and \ref%
{section_piecewise_case_fixed}, Inequality (\ref{lower_phi_M}) can be
satisfied for histogram and piecewise polynomial models on a partition
achieving some upper regularity assumption with respect to the law $P^{X}$ .

\section{The histogram case\label{section_histogram_case_fixed}}

In this section, we particularize the results stated in Section \ref%
{section_reg_fixed_model} to the case of piecewise constant functions. We
show that under a lower regularity assumption on the considered partition,
the assumption (\textbf{H4}) of existence of a localized basis in $%
L_{2}\left( P^{X}\right) $ and (\textbf{H5}) of consistency in sup-norm of
the M-estimator towards the linear projection $s_{M}$ are satisfied.

\subsection{Existence of a localized basis\label%
{section_histogram_case_fixed_localized}}

The following lemma states the existence of an orthonormal localized basis
for piecewise constant functions in $L_{2}\left( P^{X}\right) $, on a
partition which is lower-regular for the law $P^{X}$.

\begin{lemma}
\label{lemma_histograms_localized}Let consider a linear model $M$ of
histograms defined on a finite partition $\mathcal{P}$ on $\mathcal{X}$, and
write $\left\vert \mathcal{P}\right\vert =D$ the dimension of $M$. Moreover,
assume that for a positive finite constant $c_{M,P}$, 
\begin{equation}
\sqrt{\left\vert \mathcal{P}\right\vert \inf_{I\in \mathcal{P}}P^{X}\left(
I\right) }\geq c_{M,P}>0\text{ }.  \label{lower_reg_part_hist}
\end{equation}%
Set, for $I\in \mathcal{P}$, 
\begin{equation*}
\varphi _{I}=\left( P^{X}\left( I\right) \right) ^{-1/2}\mathbf{1}_{I}\text{
.\ }
\end{equation*}%
Then the family $\left( \varphi _{I}\right) _{I\in \Lambda _{M}}$ is an
orthonormal basis in $L_{2}\left( P^{X}\right) $ and we have,%
\begin{equation}
\text{for all }\beta =\left( \beta _{I}\right) _{I\in \mathcal{P}}\in 
\mathbb{R}^{D},\text{ \ }\left\Vert \sum_{I\in \mathcal{P}}\beta _{I}\varphi
_{I}\right\Vert _{\infty }\leq c_{M,P}^{-1}\sqrt{D}\left\vert \beta
\right\vert _{\infty }\text{ }.  \label{inequality_localized_histograms}
\end{equation}
\end{lemma}

\noindent Condition (\ref{lower_reg_part_hist}) can also be found in Arlot
and Massart \cite{ArlotMassart:09} and is named lower regularity of the
partition $\mathcal{P}$ for the law $P^{X}$. It is easy to see that the
lower regularity of the partition is equivalent to the property of localized
basis in the case of histograms, i.e. (\ref{lower_reg_part_hist}) is
equivalent to (\ref{inequality_localized_histograms}). The proof of Lemma %
\ref{lemma_histograms_localized} is straightforward and can be found in
Section \ref{section_proof_histogram}.

\subsection{Rates of convergence in sup-norm\label%
{section_histogram_case_fixed_consistency}}

The following lemma allows to derive property (\textbf{H5}) for histogram
models.

\begin{lemma}
\label{consistency_histograms}Consider a linear model $M$ of histograms
defined on a finite partition $\mathcal{P}$ of $\mathcal{X}$, and denote by $%
\left\vert \mathcal{P}\right\vert =D$ the dimension of $M$. Assume that
Inequality (\ref{bounded_data}) holds, that is, a positive constant $A$
exists such that $\left\vert Y\right\vert \leq A$ $a.s.$ Moreover, assume
that for some positive finite constant $c_{M,P}$, 
\begin{equation}
\sqrt{\left\vert \mathcal{P}\right\vert \inf_{I\in \mathcal{P}}P^{X}\left(
I\right) }\geq c_{M,P}>0  \label{assump_lower_part_reg}
\end{equation}%
and that $D\leq A_{+}n\left( \ln n\right) ^{-2}\leq n$ for some positive
finite constant $A_{+}.$ Then, for any $\alpha >0$ and for all $n\geq
n_{0}\left( \alpha ,c_{M,P},A_{+}\right) $, there exists an event of
probability at least $1-n^{-\alpha }$ on which $s_{n}$ exists, is unique and
it holds,%
\begin{equation}
\left\Vert s_{n}-s_{M}\right\Vert _{\infty }\leq L_{A_{+},A,c_{M,P},\alpha }%
\sqrt{\frac{D\ln n}{n}}\text{ }.  \label{rate_sup_norm_histo}
\end{equation}
\end{lemma}

\noindent In Lemma \ref{consistency_histograms} we thus achieve the
convergence in sup-norm of the regressogram $s_{n}$ towards the linear
projection $s_{M}$ at the rate $\sqrt{D\ln\left( n\right) /n}$ . It is worth
noticing that for a model of histograms satisfying the assumptions of Lemma %
\ref{consistency_histograms}, if we set 
\begin{equation*}
A_{cons}=L_{A,c_{M,P},\alpha}\sqrt{A_{+}}\text{ , }n_{1}=n_{0}\left(
\alpha,c_{M,P},A_{+}\right) \text{\ and\ }R_{n,D,\alpha}=L_{A_{+},A,c_{M,P},%
\alpha}\sqrt{\frac{D\ln n}{n}},
\end{equation*}
then Assumption \textbf{(H5) }is satisfied. To derive Inequality (\ref%
{rate_sup_norm_histo}), we need to assume that the response variable $Y$ is
almost surely bounded and that the considered partition is lower-regular for
the law $P^{X}$. Hence, we fit again with the framework of \cite%
{ArlotMassart:09} and we can thus view the general set of assumptions
exposed in Section \ref{section_main_assump_reg_fixed} as a natural
generalization for linear models of the framework developed in \cite%
{ArlotMassart:09} in the case of histograms. The proof of Lemma \ref%
{consistency_histograms} can be found in Section \ref%
{section_proof_histogram}.

\subsection{Bounds for the excess risks\label%
{section_histogram_case_fixed_bounds}}

The next results is a straightforward application of Lemmas \ref%
{lemma_histograms_localized}, \ref{consistency_histograms} and Theorems \ref%
{principal_reg}, \ref{theorem_upper_true_small_models}.

\begin{corollary}
\label{corollary_upper_lower_histograms}Given $A_{+},A_{-},\alpha >0$,
consider a linear model $M$ of histograms defined on a finite partition $%
\mathcal{P}$ of $\mathcal{X}$, and write $\left\vert \mathcal{P}\right\vert
=D$ the dimension of $M$. Assume that for some positive finite constant $%
c_{M,P}$, it holds 
\begin{equation}
\sqrt{\left\vert \mathcal{P}\right\vert \inf_{I\in \mathcal{P}}P^{X}\left(
I\right) }\geq c_{M,P}>0\text{ .}  \label{lower_reg_part_hist_theo}
\end{equation}%
If (\textbf{H1}) and (\textbf{H2}) of Section \ref%
{section_main_assump_reg_fixed} are satisfied and if 
\begin{equation*}
A_{-}\left( \ln n\right) ^{2}\leq D\leq A_{+}\frac{n}{\left( \ln n\right)
^{2}}\text{ },
\end{equation*}%
then there exists a positive finite constant $A_{0}$, only depending on $%
\alpha ,A,\sigma _{\min },A_{-},A_{+},c_{M,P}$ such that, by setting%
\begin{equation*}
\varepsilon _{n}=A_{0}\max \left\{ \left( \frac{\ln n}{D}\right)
^{1/4},\left( \frac{D\ln n}{n}\right) ^{1/4}\right\}
\end{equation*}%
we have, for all $n\geq n_{0}\left( A_{-},A_{+},A,,\sigma _{\min
},c_{M,P},\alpha \right) $,%
\begin{equation}
\mathbb{P}\left[ \left( 1+\varepsilon _{n}\right) \frac{1}{4}\frac{D}{n}%
\mathcal{K}_{1,M}^{2}\geq P\left( Ks_{n}-Ks_{M}\right) \geq \left(
1-\varepsilon _{n}\right) \frac{1}{4}\frac{D}{n}\mathcal{K}_{1,M}^{2}\right]
\geq 1-10n^{-\alpha }  \label{p1_histo}
\end{equation}%
and%
\begin{equation}
\mathbb{P}\left[ \left( 1+\varepsilon _{n}^{2}\right) \frac{1}{4}\frac{D}{n}%
\mathcal{K}_{1,M}^{2}\geq P_{n}\left( Ks_{M}-Ks_{n}\right) \geq \left(
1-\varepsilon _{n}^{2}\right) \frac{1}{4}\frac{D}{n}\mathcal{K}_{1,M}^{2}%
\right] \geq 1-5n^{-\alpha }\text{ }.  \label{p2_histo}
\end{equation}%
If (\ref{lower_reg_part_hist_theo}) holds together with (\textbf{H1}) and if
we assume that 
\begin{equation*}
1\leq D\leq A_{+}\frac{n}{\left( \ln n\right) ^{2}}\text{ ,}
\end{equation*}%
then a positive constant $A_{u}$ exists, only depending on $A,c_{M,P},A_{+}$
and $\alpha $, such that for all $n\geq n_{0}\left( A,c_{M,P},A_{+},\alpha
\right) $,%
\begin{equation*}
\mathbb{P}\left[ P\left( Ks_{n}-Ks_{M}\right) \geq A_{u}\frac{D\vee \ln n}{n}%
\right] \leq 3n^{-\alpha }
\end{equation*}%
and%
\begin{equation*}
\mathbb{P}\left[ P_{n}\left( Ks_{M}-Ks_{n}\right) \geq A_{u}\frac{D\vee \ln n%
}{n}\right] \leq 3n^{-\alpha }\text{ }.
\end{equation*}
\end{corollary}

\noindent We recover in Corollary \ref{corollary_upper_lower_histograms} the
general results of Section \ref{section_theorems_fixed model copy(1)} for
the case of histograms on a lower-regular partition. Moreover, in the case
of histograms, assumption (\ref{bounded_projection}) which is part of (%
\textbf{H1}) is a straightforward consequence of (\ref{bounded_data}).
Indeed, we easily see that the projection $s_{M}$ of the regression function 
$s_{\ast }$ onto the model of piecewise constant functions with respect to $%
\mathcal{P}$ can be written%
\begin{equation}
s_{M}=\sum_{I\in \mathcal{P}}\mathbb{E}\left[ Y\left\vert X\in I\right. %
\right] \mathbf{1}_{I}\text{ }.  \label{formule_projection}
\end{equation}%
Under (\ref{bounded_data}), we have $\left\vert \mathbb{E}\left[ Y\left\vert
X\in I\right. \right] \right\vert \leq \left\Vert Y\right\Vert _{\infty
}\leq A$ for every $I\in \mathcal{P}$ and we deduce by (\ref%
{formule_projection}) that $\left\Vert s_{M}\right\Vert _{\infty }\leq A.$

\subsection{Comments\label{section_histo_comments}}

\noindent Our bounds in Corollary \ref{corollary_upper_lower_histograms} are
obtained by following a general methodology that consists, among other
things, in expanding the contrast and to take advantage of explicit
computations that can be derived on the linear part of the contrast - for
more details, see the proofs in Section \ref{section_proof_fixed_model}\
below. It is then instructive to compare them to the best available results
in this special case. Let us compare them to the bounds obtained by Arlot
and Massart in \cite{ArlotMassart:09}, in the case of a fixed model. Such
results can be found in Propositions 10, 11 and 12 of \cite{ArlotMassart:09}.

\noindent The strategy adopted by the authors in this case is as follows.
They first notice that the mean of the empirical excess risk on histograms
is given by%
\begin{equation*}
\mathbb{E}\left[ P_{n}\left( Ks_{M}-Ks_{n}\right) \right] =\frac{D}{4n}%
\mathcal{K}_{1,M}^{2}\text{ .}
\end{equation*}%
Then they derive concentration inequalities for the true excess risk and its
empirical counterpart around their mean. Finally, the authors compare the
mean of the true excess risk to the mean of the empirical excess risk.

\noindent More precisely, using our notations, inequality (34) of
Proposition 10 in \cite{ArlotMassart:09} states that for every $x\geq 0$
there exists an event of probability at least $1-e^{1-x}$ on which,%
\begin{align}
& \left\vert P_{n}\left( Ks_{M}-Ks_{n}\right) -\mathbb{E}\left[ P_{n}\left(
Ks_{M}-Ks_{n}\right) \right] \right\vert  \notag \\
& \leq \frac{L}{\sqrt{D_{M}}}\left[ P(Ks_{M}-Ks_{\ast })+\frac{A^{2}\mathbb{E%
}\left[ P_{n}\left( Ks_{M}-Ks_{n}\right) \right] }{\sigma _{\min }^{2}}%
\left( \sqrt{x}+x\right) \right] \text{ ,}  \label{arlMas}
\end{align}%
for some absolute constant $L$. One can notice that inequality (\ref{arlMas}%
), which is a special case of general concentration inequalities given by
Boucheron and Massart \cite{BouMas:10}, involves the bias of the model $%
P(Ks_{M}-Ks_{\ast })$. By pointing out that the bias term arises from the
use of some margin conditions that are satisfied for bounded regression, we
believe that it can be removed from Proposition 10 of \cite{ArlotMassart:09}%
, since in the case of histograms models for bounded regression, some
margin-like conditions hold, that are directly pointed at the linear
projection $s_{M}$. Apart for the bias term, the deviations of the empirical
excess risk are then of the order 
\begin{equation*}
\frac{\ln \left( n\right) \sqrt{D_{M}}}{n}\text{ },
\end{equation*}%
considering the same probability of event as ours, inequality (\ref{arlMas})
becomes significantly better than inequality (\ref{p2_histo}) for large
models.

\noindent Concentration inequalities for the true excess risk given in
Proposition 11 of \cite{ArlotMassart:09} give a magnitude of deviations that
is again smaller than ours for sufficiently large models and that is in fact
closer to $\varepsilon _{n}^{2}$ than $\varepsilon _{n}$, where $\varepsilon
_{n}$ is defined in Corollary \ref{corollary_upper_lower_histograms}. But
the mean of the true excess risk has to be compared to the mean of the
empirical excess risk and it is remarkable that in Proposition 12 of \cite%
{ArlotMassart:09} where such a result is given in a way that seems very
sharp, there is a term lower bounded by 
\begin{equation*}
\left( n\times \inf_{I\in \mathcal{P}}P^{X}\left( I\right) \right)
^{-1/4}\propto \left( \frac{D}{n}\right) ^{1/4}\text{ ,}
\end{equation*}%
due to the lower regularity assumption on the partition. This tends to
indicate that, up to a logarithmic factor, the term proportional to $\left( 
\frac{D\ln n}{n}\right) ^{1/4}$ appearing in $\varepsilon _{n}$ is not
improvable in general, and that the empirical excess risk concentrates
better around its mean than the true excess risk.

\noindent We conclude that the bounds given in Proposition 10, 11 and 12 of 
\cite{ArlotMassart:09} are essentially more accurate than ours, apart for
the bias term involved in concentration inequalities of Proposition 10, but
this term could be removed as explained above. Furthermore, concentration
inequalities for the empirical excess risk are significantly sharper than
ours for large models.

\noindent Arlot and Massart \cite{ArlotMassart:09} also propose
generalizations in the case of unbounded noise and when the noise level
vanishes. The unbounded case seems to be beyond the reach of our strategy,
due to our repeated use of Bousquet and Klein-Rio's inequalities along the
proofs. However, we recover the case of vanishing noise level for histogram
models, when the partition is upper regular with respect to the law $P^{X}$,
a condition also needed in \cite{ArlotMassart:09} in this case. Indeed, we
have noticed in Section \ref{section_complexity_model} that assumption (%
\textbf{H2}) can be weakened into (\textbf{H2bis}), where we assume that 
\begin{equation*}
\mathcal{K}_{1,M}\geq A_{\min }>0
\end{equation*}%
for some positive constant $A_{\min }$. So, it suffices to bound from below
the normalized complexity. We have from identity (\ref{formule_K1M_bis}),%
\begin{equation*}
\mathcal{K}_{1,M}^{2}\geq 4\mathbb{E}\left[ \sigma ^{2}\left( X\right) \Psi
_{M}^{2}\left( X\right) \right] \text{ .}
\end{equation*}%
Moreover, from identity (\ref{phi_M}), we have in the case of histograms,%
\begin{equation*}
\Psi _{M}^{2}\left( x\right) =\frac{1}{\left\vert \mathcal{P}\right\vert }%
\sum_{I\in \mathcal{P}}\frac{\mathbf{1}_{x\in I}}{P^{X}\left( I\right) }%
\text{ , \ for all }x\in \mathcal{X}\text{ .}
\end{equation*}%
Now, if we assume the upper regularity of the partition $\mathcal{P}$ with
respect to $P^{X}$, that is%
\begin{equation}
\left\vert \mathcal{P}\right\vert \sup_{I\in \mathcal{P}}P^{X}\left(
I\right) \leq c_{M,P}^{+}<+\infty  \label{upper_reg}
\end{equation}%
for some positive constant $c_{M,P}^{+}$, we then have%
\begin{equation*}
\Psi _{M}^{2}\left( x\right) \geq \left( c_{M,P}^{+}\right) ^{-1}>0\text{ ,
\ for all }x\in \mathcal{X}\text{ ,}
\end{equation*}%
and so $A_{\min }=2\left( c_{M,P}^{+}\right) ^{-1/2}\left\Vert \sigma
\right\Vert _{2}>0$ is convenient in (\textbf{H2bis}).

\section{The case of piecewise polynomials\label%
{section_piecewise_case_fixed}}

In this Section, we generalize the results given in Section \ref%
{section_histogram_case_fixed} for models of piecewise constant functions to
models of piecewise polynomials uniformly bounded in their degrees.

\subsection{Existence of a localized basis\label%
{section_piecewise_case_fixed_localized}}

The following lemma states the existence of a localized orthonormal basis in 
$\left( M,\left\Vert 
\cdot%
\right\Vert _{2}\right) $, where $M$ is a model of piecewise polynomials and 
$\mathcal{X=}\left[ 0,1\right] $ is the unit interval.

\begin{lemma}
\label{lemma_piecewise_polynomials_localized}Let $%
\leb%
$ denote the Lebesgue measure on $\left[ 0,1\right] $. Let assume that $%
\mathcal{X=}\left[ 0,1\right] $ and that $P^{X}$ has a density $f$ with
respect to $%
\leb%
$ satisfying, for a positive constant $c_{\min }$,%
\begin{equation*}
f\left( x\right) \geq c_{\min }>0,\text{ \ }x\in \left[ 0,1\right] \text{ }.
\end{equation*}%
Consider a linear model $M$ of piecewise polynomials on $\left[ 0,1\right] $
with degree $r$ or smaller, defined on a finite partition $\mathcal{P}$ made
of intervals. Then there exists an orthonormal basis $\left\{ \varphi _{I,j},%
\text{ }I\in \mathcal{P},\text{ }j\in \left\{ 0,...,r\right\} \right\} $ of $%
\left( M,\left\Vert 
\cdot%
\right\Vert _{2}\right) $ such that,%
\begin{equation*}
\text{for all }j\in \left\{ 0,...,r\right\} \text{, \ \ \ \ \ }\varphi _{I,j}%
\text{ is supported by the element }I\text{ of }\mathcal{P}\text{,}
\end{equation*}%
and a constant $L_{r,c_{\min }}$ depending only on $r,c_{\min }$ exists,
satisfying for all $I\in \mathcal{P},$%
\begin{equation}
\max_{j\in \left\{ 0,...,r\right\} }\left\Vert \varphi _{I,j}\right\Vert
_{\infty }\leq L_{r,c_{\min }}\frac{1}{\sqrt{%
\leb%
\left( I\right) }}\text{ .}  \label{inequality_localized_leb_I_piecewise}
\end{equation}%
As a consequence, if it holds%
\begin{equation}
\sqrt{\left\vert \mathcal{P}\right\vert \inf_{I\in \mathcal{P}}%
\leb%
\left( I\right) }\geq c_{M,%
\leb%
}  \label{lower_part_leb}
\end{equation}%
a constant $L_{r,c_{\min },c_{M,%
\leb%
}}$ depending only on $r,c_{\min }$ and $c_{M,%
\leb%
}$ exists, such that for all $\beta =\left( \beta _{I,j}\right) _{I\in 
\mathcal{P},j\in \left\{ 0,...,r\right\} }\in \mathbb{R}^{D}$, 
\begin{equation}
\left\Vert \sum_{I,j}\beta _{I,j}\varphi _{I,j}\right\Vert _{\infty }\leq
L_{r,c_{\min },c_{M,%
\leb%
}}\sqrt{D}\left\vert \beta \right\vert _{\infty }
\label{inequality_localized_piecewise}
\end{equation}%
where $D=\left( r+1\right) \left\vert \mathcal{P}\right\vert $ is the
dimension of $M$.
\end{lemma}

\noindent Lemma \ref{lemma_piecewise_polynomials_localized} states that if $%
\mathcal{X=}\left[ 0,1\right] $ is the unit interval and if $P^{X}$ has a
density with respect to the Lebesgue measure $%
\leb%
$ on $\mathcal{X}$, which is uniformly bounded away form zero, then there
exists an orthonormal basis in $\left( M,\left\Vert 
\cdot%
\right\Vert _{2}\right) $ satisfying good enough properties in terms of the
sup-norm of its elements. Moreover, if we assume the lower regularity of the
partition with respect to $%
\leb%
$, then the orthonormal basis is localized and the constant of localization
given in (\ref{inequality_localized_piecewise}) depend on the maximal degree 
$r$. We notice that in the case of piecewise constant functions we do not
need to assume the existence of a density for $P^{X}$ or to restrict
ourselves to the unit interval. The proof of Lemma \ref%
{lemma_piecewise_polynomials_localized} can be found in Section \ref%
{section_proof_piecewise}.

\subsection{Rates of convergence in sup-norm\label%
{section_piecewise_case_fixed_consistency}}

The following lemma allows to derive property (\textbf{H5}) for piecewise
polynomials.

\begin{lemma}
\label{lemma_piecewise_polynomials_consistency}Assume that Inequality (\ref%
{bounded_data}) holds, that is a positive constant $A$ exists such that $%
\left\vert Y\right\vert \leq A$ $a.s.$ Denote by $%
\leb%
$ the Lebesgue measure on $\left[ 0,1\right] $. Assume that $\mathcal{X=}%
\left[ 0,1\right] $ and that $P^{X}$ has a density $f$ with respect to $%
\leb%
$, satisfying for positive constants $c_{\min }$ and $c_{\max }$,%
\begin{equation}
0<c_{\min }\leq f\left( x\right) \leq c_{\max }<+\infty ,\text{ }x\in \left[
0,1\right] \text{ }.  \label{density_piecwise_consistency}
\end{equation}%
Consider a linear model $M$ of piecewise polynomials on $\left[ 0,1\right] $
with degree less than $r$, defined on a finite partition $\mathcal{P}$ made
of intervals, that satisfies for some finite positive constants $c_{M,%
\leb%
}$%
\begin{equation}
\sqrt{\left\vert \mathcal{P}\right\vert \inf_{I\in \mathcal{P}}%
\leb%
\left( I\right) }\geq c_{M,%
\leb%
}>0\text{ }.  \label{partition_piecewise_consistency}
\end{equation}%
Assume moreover that $D\leq A_{+}n\left( \ln n\right) ^{-2}$ for a positive
finite constant $A_{+}.$ Then, for any $\alpha >0$, there exists an event of
probability at least $1-n^{-\alpha }$ such that $s_{n}$ exists, is unique on
this event and it holds, for all $n\geq n_{0}\left( r,A_{+},c_{\min },c_{M,%
\leb%
},\alpha \right) $,%
\begin{equation}
\left\Vert s_{n}-s_{M}\right\Vert _{\infty }\leq L_{A,r,A_{+},c_{\min
},c_{\max },c_{M,%
\leb%
},\alpha }\sqrt{\frac{D\ln n}{n}}\text{ }.  \label{rate_sup_norm_localized}
\end{equation}
\end{lemma}

\noindent In Lemma \ref{consistency_histograms}, we thus obtain the
convergence in sup-norm of the M-estimator $s_{n}$ toward the linear
projection $s_{M}$ at the rate $\sqrt{D\ln \left( n\right) /n}$ . It is
worth noting that, for a model of piecewise polynomials satisfying the
assumptions of Lemma \ref{consistency_histograms}, if we set 
\begin{gather*}
A_{cons}=L_{A,r,A_{+},c_{\min },c_{\max },c_{M,%
\leb%
},\alpha }\sqrt{A_{+}}\text{ \ },\text{ \ }R_{n,D,\alpha
}=L_{A,r,A_{+},c_{\min },c_{\max },c_{M,%
\leb%
},\alpha }\sqrt{\frac{D\ln n}{n}}\text{ }, \\
n_{1}=n_{0}\left( r,A_{+},c_{\min },c_{M,%
\leb%
},\alpha \right) \text{ },
\end{gather*}%
then Assumption \textbf{(H5) }is satisfied. The proof of Lemma \ref%
{lemma_piecewise_polynomials_consistency} can be found in Section \ref%
{section_proof_piecewise}.

\subsection{Bounds for the excess risks\label%
{section_piecewise_case_fixed_bounds}}

The forthcoming result is a straightforward application of Lemmas \ref%
{lemma_piecewise_polynomials_localized}, \ref%
{lemma_piecewise_polynomials_consistency} and Theorems \ref{principal_reg}, %
\ref{theorem_upper_true_small_models}.

\begin{corollary}
\label{corollary_upper_lower_piecewise_ploynomials}Denote by $%
\leb%
$ the Lebesgue measure on $\left[ 0,1\right] $ and fix some positive finite
constant $\alpha $. Assume that $\mathcal{X=}\left[ 0,1\right] $ and that $%
P^{X}$ has a density $f$ with respect to $%
\leb%
$ satisfying, for some positive finite constants $c_{\min }$ and $c_{\max }$,%
\begin{equation}
0<c_{\min }\leq f\left( x\right) \leq c_{\max }<+\infty ,\text{ }x\in \left[
0,1\right] \text{ }.  \label{control_density}
\end{equation}%
Consider a linear model $M$ of piecewise polynomials on $\left[ 0,1\right] $
with degree less than $r$, defined on a finite partition $\mathcal{P}$ made
of intervals, that satisfy for a finite constant $c_{M,%
\leb%
}$, 
\begin{equation}
\sqrt{\left\vert \mathcal{P}\right\vert \inf_{I\in \mathcal{P}}%
\leb%
\left( I\right) }\geq c_{M,%
\leb%
}>0\text{ }.  \label{lower_reg_part_pp_theo}
\end{equation}%
Assume that (\textbf{H1}) and (\textbf{H2}) hold. Then, if there exist some
positive finite constants $A_{-}$ and $A_{+}$ such that%
\begin{equation*}
A_{-}\left( \ln n\right) ^{2}\leq D\leq A_{+}\frac{n}{\left( \ln n\right)
^{2}}\text{ },
\end{equation*}%
then there exists a positive finite constant $A_{0}$, depending on $\alpha
,A,\sigma _{\min },A_{-},A_{+},r,c_{M,%
\leb%
},c_{\min }$ and $c_{\max }$ such that, by setting%
\begin{equation*}
\varepsilon _{n}=A_{0}\max \left\{ \left( \frac{\ln n}{D}\right)
^{1/4},\left( \frac{D\ln n}{n}\right) ^{1/4}\right\}
\end{equation*}%
we have, for all $n\geq n_{0}\left( A_{-},A_{+},A,r,\sigma _{\min },c_{M,%
\leb%
},c_{\min },c_{\max },\alpha \right) $,%
\begin{equation*}
\mathbb{P}\left[ \left( 1+\varepsilon _{n}\right) \frac{1}{4}\frac{D}{n}%
\mathcal{K}_{1,M}^{2}\geq P\left( Ks_{n}-Ks_{M}\right) \geq \left(
1-\varepsilon _{n}\right) \frac{1}{4}\frac{D}{n}\mathcal{K}_{1,M}^{2}\right]
\geq 1-10n^{-\alpha }
\end{equation*}%
and%
\begin{equation*}
\mathbb{P}\left[ \left( 1+\varepsilon _{n}^{2}\right) \frac{1}{4}\frac{D}{n}%
\mathcal{K}_{1,M}^{2}\geq P_{n}\left( Ks_{M}-Ks_{n}\right) \geq \left(
1-\varepsilon _{n}^{2}\right) \frac{1}{4}\frac{D}{n}\mathcal{K}_{1,M}^{2}%
\right] \geq 1-5n^{-\alpha }\text{ }.
\end{equation*}%
Moreover, if (\ref{control_density}) and (\ref{lower_reg_part_pp_theo}) hold
together with (\textbf{H1}) and if we assume that 
\begin{equation*}
1\leq D\leq A_{+}\frac{n}{\left( \ln n\right) ^{2}}\text{ ,}
\end{equation*}%
then a positive constant $A_{u}$ exists, only depending on $A_{+},A,r,c_{M,%
\leb%
},c_{\min }$ and $\alpha $, such that for all $n\geq n_{0}\left(
A_{+},A,r,c_{\min },c_{\max },c_{M,%
\leb%
},\alpha \right) $,%
\begin{equation*}
\mathbb{P}\left[ P\left( Ks_{n}-Ks_{M}\right) \geq A_{u}\frac{D\vee \ln n}{n}%
\right] \leq 3n^{-\alpha }
\end{equation*}%
and%
\begin{equation*}
\mathbb{P}\left[ P_{n}\left( Ks_{M}-Ks_{n}\right) \geq A_{u}\frac{D\vee \ln n%
}{n}\right] \leq 3n^{-\alpha }\text{ }.
\end{equation*}
\end{corollary}

\noindent We derive in Corollary \ref%
{corollary_upper_lower_piecewise_ploynomials} optimal upper and lower bounds
for the excess risk and its empirical counterpart in the case of models of
piecewise polynomials uniformly bounded in their degree, with reasonable
dimension. We give also upper bounds for models of possibly small dimension,
without assumption (\textbf{H2}). Notice that we need stronger assumptions
than in the case of histograms. Namely, we require the existence of a
density uniformly bounded from above and from below for the unknown law $%
P^{X}$, with respect to the Lebesgue measure on the unit interval. However,
we recover essentially the bounds of Corollary \ref%
{corollary_upper_lower_histograms}, since by Lemma \ref%
{lemma_piecewise_polynomials_consistency}, we still have $R_{n,D,\alpha
}\propto \sqrt{D\ln \left( n\right) /n}$.

\noindent Moreover, as in the case of histograms, assumption (\ref%
{bounded_projection}) which is part of (\textbf{H1}), is a straightforward
consequence of (\ref{bounded_data}). Indeed, we easily see that the
projection $s_{M}$ of the regression function $s_{\ast }$ onto the model of
piecewise polynomials with respect to $\mathcal{P}$ can be written 
\begin{equation*}
s_{M}=\sum_{\left( I,j\right) \in \mathcal{P}\times \left\{ 0,...,r\right\}
}P\left( Y\varphi _{I,j}\right) \varphi _{I,j}\text{ ,}
\end{equation*}%
where $\varphi _{I,j}$ is the orthonormal basis given in Lemma \ref%
{lemma_piecewise_polynomials_localized}. It is then easy to show, using (\ref%
{inequality_localized_leb_I_piecewise}) and (\ref{bounded_data}), that $%
\left\Vert s_{M}\right\Vert _{\infty }\leq L_{A,r,c_{\min },c_{\max }}.$

\noindent Again, we can consider vanishing noise at the prize to ask that
the partition is upper regular with respect to $%
\leb%
$. By (\textbf{H2bis}) of Section \ref{section_complexity_model}, if we show
that 
\begin{equation*}
\mathcal{K}_{1,M}\geq A_{\min }>0
\end{equation*}%
for a positive constant $A_{\min }$ instead of (\textbf{H2}), then the
conclusions of Corollary \ref{corollary_upper_lower_piecewise_ploynomials}
still hold. Now, from identity (\ref{formule_K1M_bis}) we have%
\begin{equation*}
\mathcal{K}_{1,M}^{2}\geq 4\mathbb{E}\left[ \sigma ^{2}\left( X\right) \Psi
_{M}^{2}\left( X\right) \right] \text{ .}
\end{equation*}%
Moreover, from identity (\ref{phi_M}), it holds in the case of piecewise
polynomials, for all $x\in \mathcal{X}$,%
\begin{equation}
\Psi _{M}^{2}\left( x\right) =\frac{1}{\left( r+1\right) \left\vert \mathcal{%
P}\right\vert }\sum_{\left( I,j\right) \in \mathcal{P}\times \left\{
0,...,r\right\} }\varphi _{I,j}^{2}\geq \frac{1}{\left( r+1\right)
\left\vert \mathcal{P}\right\vert }\sum_{I\in \mathcal{P}}\frac{\mathbf{1}%
_{x\in I}}{P^{X}\left( I\right) }\text{ .}  \label{lower_phi_m}
\end{equation}%
Furthermore, if we ask that%
\begin{equation}
\left\vert \mathcal{P}\right\vert \sup_{I\in \mathcal{P}}%
\leb%
\left( I\right) \leq c_{M,P}^{+}<+\infty  \label{upper_part}
\end{equation}%
for some positive constant $c_{M,P}^{+}$, then\ by using (\ref%
{control_density}), (\ref{lower_phi_m}) and (\ref{upper_part}), we obtain
for all $x\in \mathcal{X}$,%
\begin{equation*}
\Psi _{M}^{2}\left( x\right) \geq \left( c_{\max }\times c_{M,P}^{+}\times
\left( r+1\right) \right) ^{-1}>0\text{ ,}
\end{equation*}%
and so $A_{\min }=2\left( c_{\max }\times c_{M,P}^{+}\times \left(
r+1\right) \right) ^{-1/2}\sqrt{\mathbb{E}\left[ \sigma ^{2}\left( X\right) %
\right] }>0$ is convenient in (\textbf{H2bis}).

\section{Proofs}

We begin with the simpler proofs of Sections \ref%
{section_histogram_case_fixed} and \ref{section_piecewise_case_fixed}, in
Sections \ref{section_proof_histogram} and \ref{section_proof_piecewise}
respectively. The proofs of Theorems \ref{principal_reg} and \ref%
{theorem_upper_true_small_models} of Section \ref{section_theorems_fixed
model copy(1)} can be found in Section \ref{section_proof_fixed_model}.

\subsection{Proofs of Section \protect\ref{section_histogram_case_fixed} 
\label{section_proof_histogram}}

\subsection*{}%
\textbf{Proof of Lemma \ref{lemma_histograms_localized}. }It suffices to
observe that%
\begin{align*}
\left\Vert \sum_{I\in\mathcal{P}}\beta_{I}\varphi_{I}\right\Vert _{\infty} &
\leq\left\vert \beta\right\vert _{\infty}\sup_{I\in\mathcal{P}}\left\Vert
\varphi_{I}\right\Vert _{\infty} \\
& =\left\vert \beta\right\vert _{\infty}\sup_{I\in\mathcal{P}}\left(
P^{X}\left( I\right) \right) ^{-1/2} \\
& \leq c_{M,P}^{-1}\sqrt{D}\left\vert \beta\right\vert _{\infty}\text{ }.
\end{align*}
\ \ \ \ \ \ \ \ \ \ \ \ \ \ \ \ \ \ \ \ \ \ \ \ \ \ \ \ \ \ \ \ \ \ \ \ \ \
\ \ \ \ \ \ \ \ \ \ \ \ \ \ \ \ \ \ \ \ \ \ \ \ \ \ \ \ \ \ \ \ \ \ \ \ \ \
\ \ \ \ \ \ \ \ \ \ \ \ \ \ \ \ \ \ \ \ \ \ \ \ \ \ \ \ \ \ \ \ \ \ \ \ \ \
\ \ \ \ \ \ \ \ \ \ \ \ \ \ \ \ \ \ \ \ \ \ \ \ \ \ \ \ \ \ \ \ \ \ \ \ \ \
\ \ \ \ 

\ \ \ \ \ \ \ \ \ \ \ \ \ \ \ \ \ \ \ \ \ \ \ \ \ \ \ \ \ \ \ \ \ \ \ \ \ \
\ \ \ \ \ \ \ \ \ \ \ \ \ \ \ \ \ \ \ \ \ \ \ \ \ \ \ \ \ \ \ \ \ \ \ \ \ \
\ \ \ \ \ \ \ \ \ \ \ \ \ \ \ \ \ \ \ \ \ \ \ \ \ \ \ \ \ \ \ \ \ \ \ \ \ \
\ \ \ \ \ \ \ \ \ \ \ \ \ \ \ \ \ \ \ \ \ \ \ \ \ \ \ \ \ \ \ 
$\blacksquare$%

We now intend to prove (\ref{rate_sup_norm_histo}) under the assumptions of
Lemma \ref{consistency_histograms}.

\subsection*{}%
\textbf{Proof of Lemma \ref{consistency_histograms}. }Along the proof, we
denote by abusing the notation, for any $I\in \mathcal{P}$,%
\begin{equation*}
P\left( I\right) :=P\left( I\times \mathbb{R}\right) =P^{X}\left( I\right) 
\text{ and }P_{n}\left( I\right) :=P_{n}\left( I\times \mathbb{R}\right) 
\text{ .}
\end{equation*}%
Let $\alpha >0$ be fixed and let $\beta >0$ to be chosen later. We first
show that, since we have $D\leq A_{+}n\left( \ln n\right) ^{-2}$, it holds
with large probability and for all $n$ sufficiently large,%
\begin{equation*}
\inf_{I\in \mathcal{P}}P_{n}\left( I\right) >0\text{ }.
\end{equation*}%
Since 
\begin{equation*}
\left\Vert \mathbf{1}_{I}\right\Vert _{\infty }\leq 1\text{ \ \ \ and \ \ \ }%
\mathbb{E}\left[ \mathbf{1}_{I}^{2}\right] =P\left( I\right) \text{ ,}
\end{equation*}%
we get by Bernstein's inequality (\ref{bernstein_ineq_reg}), for any $x>0$
and $I\in \mathcal{P}$,%
\begin{equation}
\mathbb{P}\left[ \left\vert \left( P_{n}-P\right) \left( I\right)
\right\vert \geq \sqrt{\frac{2P\left( I\right) x}{n}}+\frac{x}{3n}\right]
\leq 2\exp \left( -x\right) \text{ }.  \label{bern_hist}
\end{equation}%
Further note that by (\ref{assump_lower_part_reg}), $D\geq
c_{M,P}^{2}P\left( I\right) ^{-1}>0$ for any $I\in \mathcal{P}$, and thus by
taking $x=\beta \ln n$, we easily deduce from inequality (\ref{bern_hist})
that there exists a positive constant $L_{\beta ,c_{M,P,A_{+}}}^{(1)}$ only
depending on $c_{M,P}$ and $\beta $ such that, for any $I\in \mathcal{P}$,%
\begin{equation}
\mathbb{P}\left[ \frac{\left\vert \left( P_{n}-P\right) \left( I\right)
\right\vert }{P\left( I\right) }\geq L_{\beta ,c_{M,P,A_{+}}}^{(1)}\sqrt{%
\frac{D\ln n}{n}}\right] \leq 2n^{-\beta }\text{ }.
\label{bernstein_I_proof}
\end{equation}%
Now, as $D\leq A_{+}n\left( \ln n\right) ^{-2}$ for some positive constant $%
A_{+}$, a positive integer $n_{0}\left( \beta ,c_{M,P},A_{+}\right) $ exists
such that 
\begin{equation}
L_{\beta ,c_{M,P,A_{+}}}^{(1)}\sqrt{\frac{D\ln n}{n}}\leq \frac{1}{2}\text{,
for all }n\geq n_{0}\left( \beta ,c_{M,P},A_{+}\right) \text{ }.
\label{def_n0_hist_proof}
\end{equation}%
Therefore we get, for all $n\geq n_{0}\left( \beta ,c_{M,P},A_{+}\right) $,%
\begin{align*}
& \mathbb{P}\left[ \forall I\in \mathcal{P},\text{ }P_{n}\left( I\right) >0%
\right] \\
& \geq \mathbb{P}\left[ \forall I\in \mathcal{P},\text{ }\frac{P\left(
I\right) }{2}>\left\vert \left( P_{n}-P\right) \left( I\right) \right\vert %
\right] \\
& \geq \mathbb{P}\left[ \forall I\in \mathcal{P},\text{ }\frac{\left\vert
\left( P_{n}-P\right) \left( I\right) \right\vert }{P\left( I\right) }%
<L_{\beta ,c_{M,P,A_{+}}}^{(1)}\sqrt{\frac{D\ln n}{n}}\right] \text{ \ by (%
\ref{def_n0_hist_proof})} \\
& \geq 1-2Dn^{-\beta }\text{ }.
\end{align*}%
Introduce the event 
\begin{equation*}
\Omega _{+}=\left\{ \forall I\in \mathcal{P},\text{ }P_{n}\left( I\right)
>0\right\} \text{ }.
\end{equation*}%
We have shown that 
\begin{equation}
\mathbb{P}\left[ \Omega _{+}\right] \geq 1-2Dn^{-\beta }\text{ }.
\label{proba_omega_plus}
\end{equation}%
Moreover, on the event $\Omega _{+}$, the least squares estimator $s_{n}$
exists, is unique and it holds%
\begin{equation*}
s_{n}=\sum_{I\in \mathcal{P}}\frac{P_{n}\left( y\mathbf{1}_{x\in I}\right) }{%
P_{n}\left( I\right) }\mathbf{1}_{I}\text{ }.
\end{equation*}%
We also have 
\begin{equation*}
s_{M}=\sum_{I\in \mathcal{P}}\frac{P\left( y\mathbf{1}_{x\in I}\right) }{%
P\left( I\right) }\mathbf{1}_{I}\text{ }.
\end{equation*}%
Hence it holds on $\Omega _{+},$%
\begin{align}
\left\Vert s_{n}-s_{M}\right\Vert _{\infty }& =\sup_{I\in \mathcal{P}%
}\left\vert \frac{P_{n}\left( y\mathbf{1}_{x\in I}\right) }{P_{n}\left(
I\right) }-\frac{P\left( y\mathbf{1}_{x\in I}\right) }{P\left( I\right) }%
\right\vert  \notag \\
& =\sup_{I\in \mathcal{P}}\left\vert \frac{P_{n}\left( y\mathbf{1}_{x\in
I}\right) }{P\left( I\right) \left( 1+\frac{\left( P_{n}-P\right) \left(
I\right) }{P\left( I\right) }\right) }-\frac{P\left( y\mathbf{1}_{x\in
I}\right) }{P\left( I\right) }\right\vert  \notag \\
& \leq \sup_{I\in \mathcal{P}}\left\vert \frac{\left( P_{n}-P\right) \left( y%
\mathbf{1}_{x\in I}\right) }{P\left( I\right) \left( 1+\frac{\left(
P_{n}-P\right) \left( I\right) }{P\left( I\right) }\right) }\right\vert 
\notag \\
& +\sup_{I\in \mathcal{P}}\left\vert \frac{P\left( y\mathbf{1}_{x\in
I}\right) }{P\left( I\right) }\right\vert \times \sup_{I\in \mathcal{P}%
}\left\vert 1-\frac{1}{1+\frac{\left( P_{n}-P\right) \left( I\right) }{%
P\left( I\right) }}\right\vert \text{ .}
\label{control_norm_infini_proof_hist_bis}
\end{align}%
Moreover, by Bernstein's inequality (\ref{bernstein_ineq_reg}), as%
\begin{equation*}
\left\Vert y\mathbf{1}_{x\in I}\right\Vert _{\infty }\leq A\text{ \ \ \ \
and \ \ \ \ }\mathbb{E}\left[ \left( Y\mathbf{1}_{X\in I}\right) ^{2}\right]
\leq A^{2}P\left( I\right) \text{ ,}
\end{equation*}%
we get for all $I\in \mathcal{P}$,%
\begin{equation*}
\mathbb{P}\left[ \left\vert \left( P_{n}-P\right) \left( y\mathbf{1}_{x\in
I}\right) \right\vert \geq \sqrt{\frac{2A^{2}P\left( I\right) x}{n}}+\frac{Ax%
}{3n}\right] \leq 2\exp \left( -x\right) \text{ }.
\end{equation*}%
By putting $x=\beta \ln n$ in the latter inequality and using the fact that $%
D\geq c_{M,P}^{2}P\left( I\right) ^{-1}$ it follows that there exists a
positive constant $L_{A,c_{M,P},\beta ,A_{+}}^{(2)}$ only depending on $A,$ $%
c_{M,P}$ and $\beta $ such that%
\begin{equation}
\mathbb{P}\left[ \frac{\left\vert \left( P_{n}-P\right) \left( y\mathbf{1}%
_{x\in I}\right) \right\vert }{P\left( I\right) }\geq L_{A,c_{M,P},\beta
,A_{+}}^{(2)}\sqrt{\frac{D\ln n}{n}}\right] \leq 2n^{-\beta }\text{ }.
\label{dev_y_I}
\end{equation}%
Now define%
\begin{equation*}
\Omega _{1,2}=\bigcap\limits_{I\in \mathcal{P}}\left\{ \left\{ \frac{%
\left\vert \left( P_{n}-P\right) \left( I\right) \right\vert }{P\left(
I\right) }<L_{\beta ,c_{M,P,A_{+}}}^{(1)}\sqrt{\frac{D\ln n}{n}}\right\}
\bigcap \left\{ \frac{\left\vert \left( P_{n}-P\right) \left( y\mathbf{1}%
_{x\in I}\right) \right\vert }{P\left( I\right) }<L_{A,c_{M,P},\beta
,A_{+}}^{(2)}\sqrt{\frac{D\ln n}{n}}\right\} \right\} \text{ }.
\end{equation*}%
Clearly, since $D\leq n$ we have, by (\ref{bernstein_I_proof}) and (\ref%
{dev_y_I}), 
\begin{equation}
\mathbb{P}\left[ \Omega _{1,2}^{c}\right] \leq 4n^{-\beta +1}\text{ }.
\label{proba_omega_1_2}
\end{equation}%
Moreover, for all $n\geq n_{0}\left( \beta ,c_{M,P},A_{+}\right) $, we get
by (\ref{def_n0_hist_proof}) that%
\begin{equation*}
\frac{\left\vert \left( P_{n}-P\right) \left( I\right) \right\vert }{P\left(
I\right) }<\frac{1}{2}
\end{equation*}%
on the event $\Omega _{1,2}$, and so, for all $n\geq n_{0}\left( \beta
,c_{M,P},A_{+}\right) $, $\Omega _{1,2}\subset \Omega _{+}$. Hence, we get
that%
\begin{align}
& \sup_{I\in \mathcal{P}}\left\vert \frac{\left( P_{n}-P\right) \left( y%
\mathbf{1}_{x\in I}\right) }{P\left( I\right) \left( 1+\frac{\left(
P_{n}-P\right) \left( I\right) }{P\left( I\right) }\right) }\right\vert
+\sup_{I\in \mathcal{P}}\left\vert \frac{P\left( y\mathbf{1}_{x\in I}\right) 
}{P\left( I\right) }\right\vert \times \sup_{I\in \mathcal{P}}\left\vert 1-%
\frac{1}{1+\frac{\left( P_{n}-P\right) \left( I\right) }{P\left( I\right) }}%
\right\vert  \notag \\
& \leq 2\sup_{I\in \mathcal{P}}\left\vert \frac{\left( P_{n}-P\right) \left(
y\mathbf{1}_{x\in I}\right) }{P\left( I\right) }\right\vert +2\sup_{I\in 
\mathcal{P}}\left\vert \frac{P\left( y\mathbf{1}_{x\in I}\right) }{P\left(
I\right) }\right\vert \times \sup_{I\in \mathcal{P}}\left\vert \frac{\left(
P_{n}-P\right) \left( I\right) }{P\left( I\right) }\right\vert  \notag \\
& \leq 2L_{A,c_{M,P},\beta ,A_{+}}^{\left( 2\right) }\sqrt{\frac{D\ln n}{n}}%
+2L_{\beta ,c_{M,P,A_{+}}}^{\left( 1\right) }\sqrt{\frac{D\ln n}{n}}\times
\sup_{I\in \mathcal{P}}\left\vert \frac{P\left( y\mathbf{1}_{x\in I}\right) 
}{P\left( I\right) }\right\vert \text{ }.
\label{control_norm_infini_proof_hist_2}
\end{align}%
Finally we have, for any $I\in \mathcal{P}$,%
\begin{equation}
\left\vert P\left( y\mathbf{1}_{x\in I}\right) \right\vert \leq P\left(
\left\vert y\right\vert \mathbf{1}_{x\in I}\right) \leq AP\left( I\right) ,
\label{last_cons_histo}
\end{equation}%
so by (\ref{control_norm_infini_proof_hist_bis}), (\ref%
{control_norm_infini_proof_hist_2}) and (\ref{last_cons_histo}) we finally
get, on the event $\Omega _{1,2}$ and for all $n\geq n_{0}\left( \beta
,c_{M,P},A_{+}\right) $,%
\begin{equation*}
\left\Vert s_{n}-s_{M}\right\Vert _{\infty }\leq \left( 2L_{A,c_{M,P},\beta
,A_{+}}^{\left( 2\right) }+2AL_{\beta ,c_{M,P,A_{+}}}^{\left( 1\right)
}\right) \sqrt{\frac{D\ln n}{n}}\text{ .}
\end{equation*}%
Taking $\beta =\alpha +3$, we get by (\ref{proba_omega_1_2}) for all $n\geq
2 $, $\mathbb{P}\left[ \Omega _{1,2}^{c}\right] \leq n^{-\alpha }$ which
implies (\ref{rate_sup_norm_histo}).

\ \ \ \ \ \ \ \ \ \ \ \ \ \ \ \ \ \ \ \ \ \ \ \ \ \ \ \ \ \ \ \ \ \ \ \ \ \
\ \ \ \ \ \ \ \ \ \ \ \ \ \ \ \ \ \ \ \ \ \ \ \ \ \ \ \ \ \ \ \ \ \ \ \ \ \
\ \ \ \ \ \ \ \ \ \ \ \ \ \ \ \ \ \ \ \ \ \ \ \ \ \ \ \ \ \ \ \ \ \ \ \ \ \
\ \ \ \ \ \ \ \ \ \ \ \ \ \ \ \ \ \ \ \ \ \ \ \ \ \ \ \ \ \ \ 
$\blacksquare$%

\subsection{Proofs of Section \protect\ref{section_piecewise_case_fixed} 
\label{section_proof_piecewise}}

Under the assumptions of Lemma \ref{lemma_piecewise_polynomials_localized},
we intend to establish (\ref{inequality_localized_piecewise}).

\subsection*{}%
\textbf{Proof of Lemma \ref{lemma_piecewise_polynomials_localized}.} Let $I$
be any interval of $\left[ 0,1\right] $ and $w$ a positive measurable
function on $I$. Denote by $L_{2}\left( I,%
\leb%
\right) $ the space of square integrable functions on $I$ with respect to
the Lebesgue measure $%
\leb%
$ and set 
\begin{equation*}
L_{2}\left( I,w\right) =\left\{ g:I\longrightarrow \mathbb{R}\text{ ; }g%
\sqrt{w}\in L_{2}\left( I,%
\leb%
\right) \right\} \text{ }.
\end{equation*}%
This space is equipped with the natural inner product%
\begin{equation*}
\left\langle g,h\right\rangle _{I,w}=\int_{x\in I}g\left( x\right) h\left(
x\right) w\left( x\right) dx\text{ }.
\end{equation*}%
Write $\left\Vert .\right\Vert _{I,w}$ its associated norm.

\noindent Now, consider an interval $I$ of $\mathcal{P}$ with bounds $a$ and 
$b$, $a<b$. Also denote by $f_{\mid I}:x\in I\longmapsto f\left( x\right) $
the restriction of the density $f$ to the interval $I.$ We readily have for $%
g,h\in L_{2}\left( I,f_{\mid I}\right) ,$%
\begin{align}
& \int\limits_{x\in I}g\left( x\right) h\left( x\right) f_{\mid I}\left(
x\right) \frac{dx}{%
\leb%
\left( I\right) }  \notag \\
& =\int\limits_{y\in \left[ 0,1\right] }g\left( \left( b-a\right) y+a\right)
h\left( \left( b-a\right) y+a\right) f_{|I}\left( \left( b-a\right)
y+a\right) dy\text{ .}  \label{transport_L2_end}
\end{align}%
Define the function $f^{I}$ from $\left[ 0,1\right] $ to $\mathbb{R}_{+}$ by 
\begin{equation*}
\text{\ }f^{I}\left( y\right) =f_{\mid I}\left( \left( b-a\right) y+a\right)
,\text{ \ }y\in \left[ 0,1\right] \text{ }.
\end{equation*}%
If $\left( p_{I,0},p_{I,1},...p_{I,r}\right) $ is an orthonormal family of
polynomials in $L_{2}\left( \left[ 0,1\right] ,f^{I}\right) $ then by
setting, for all $x\in I$, $j\in \left\{ 0,...,r\right\} $,%
\begin{equation*}
\tilde{\varphi}_{I,j}\left( x\right) =p_{I,j}\left( \frac{x-a}{b-a}\right) 
\frac{1}{\sqrt{%
\leb%
\left( I\right) }}\text{ },
\end{equation*}%
we deduce from equality (\ref{transport_L2_end}) that $\left( \tilde{\varphi}%
_{I,j}\right) _{j=0}^{r}$ is an orthonormal family of polynomials in $%
L_{2}\left( I,f_{\mid I}\right) $ such that $\deg \left( \tilde{\varphi}%
_{I,j}\right) =\deg \left( p_{I,j}\right) $.

\noindent Now, it is a classical fact of orthogonal polynomials theory (see
for example Theorems 1.11 and 1.12 of \cite{CrouMignot:84}) that there
exists a unique family $\left( q_{I,0},q_{I,1},...q_{I,r}\right) $ of
orthogonal polynomials on $\left[ 0,1\right] $ such that $\deg \left(
q_{I,j}\right) =j$ and the coefficient of the highest monomial $x^{j}$ of $%
q_{I,j}$ is equal to $1$. Moreover, each $q_{I,j}$ has $j$ distinct real
roots belonging to $\left] 0,1\right[ $. Thus, we can write%
\begin{equation}
q_{I,j}\left( x\right) =\prod\limits_{k=1}^{j}\left( x-\alpha
_{I,j}^{k}\right) ,\text{ \ \ }\alpha _{I,j}^{k}\in \left] 0,1\right[ \text{
and\ }\alpha _{I,j}^{k}\neq \alpha _{I,j}^{l}\text{\ for }k\neq l\text{ .}
\label{qIj_scinde_end}
\end{equation}%
Clearly, $\left\Vert q_{I,j}\right\Vert _{\infty }\leq 1.$ Moreover,%
\begin{align*}
\left\Vert q_{I,j}\right\Vert _{\left[ 0,1\right] ,f^{I}}^{2}& =\int\limits_{%
\left[ 0,1\right] }\left( q_{I,j}\right) ^{2}f^{I}dx \\
& \geq c_{\min }\int\limits_{\left[ 0,1\right] }\left( q_{I,j}\right) ^{2}dx%
\text{ }.
\end{align*}%
Now we set $B\left( \alpha ,r\right) =\left] \alpha -r,\alpha +r\right[ $
for $\alpha \in \mathbb{R}$, so that by (\ref{qIj_scinde_end}) we get 
\begin{equation*}
\forall x\in \left[ 0,1\right] \setminus \cup _{k=1}^{j}B\left( \alpha
_{I,j}^{k},\left( 4j\right) ^{-1}\right) ,\text{ \ }\left\vert q_{I,j}\left(
x\right) \right\vert \geq \left( 4j\right) ^{-j}\text{ },
\end{equation*}%
and%
\begin{equation*}
\leb%
\left( \left[ 0,1\right] \setminus \cup _{k=1}^{j}B\left( \alpha
_{I,j}^{k},\left( 4j\right) ^{-1}\right) \right) \geq \frac{1}{2}\text{ }.
\end{equation*}%
Therefore,%
\begin{align*}
\left\Vert q_{I,j}\right\Vert _{\left[ 0,1\right] ,f^{I}}^{2}& \geq c_{\min
}\int\limits_{\left[ 0,1\right] }\left( q_{I,j}\right) ^{2}dx \\
& \geq c_{\min }\int\limits_{\left[ 0,1\right] \setminus \cup
_{k=1}^{j}B\left( \alpha _{I,j}^{k},\left( 4j\right) ^{-1}\right) }\left(
q_{I,j}\right) ^{2}dx \\
& \geq \frac{c_{\min }}{2}\left( 4j\right) ^{-2j}\text{ }.
\end{align*}%
Finally, introduce $p_{I,j}=\left\Vert q_{I,j}\right\Vert _{\left[ 0,1\right]
,f^{I}}^{-1}q_{I,j}$ and denote by $\varphi _{I,j}$ its associated
orthonormal family of $L_{2}\left( I,f_{\mid I}\right) .$ Then, by
considering the extension $\varphi _{I,j}$ of $\tilde{\varphi}_{I,j}$ to $%
\left[ 0,1\right] $ by adding null values, it is readily checked that the
family 
\begin{equation*}
\left\{ \varphi _{I,j},\text{ }I\in \mathcal{P},\text{ }j\in \left\{
0,...,r\right\} \right\}
\end{equation*}%
is an orthonormal basis of $\left( M,\left\Vert 
\cdot%
\right\Vert _{2}\right) .$ In addition,%
\begin{align}
\left\Vert \varphi _{I,j}\right\Vert _{\infty }& =\left\Vert \tilde{\varphi}%
_{I,j}\right\Vert _{\infty }  \notag \\
& =\left\Vert q_{I,j}\right\Vert _{\left[ 0,1\right] ,f^{I}}^{-1}\left\Vert
q_{I,j}\right\Vert _{\infty }%
\leb%
\left( I\right) ^{-1/2}  \notag \\
& \leq \sqrt{2}c_{\min }^{-1/2}\left( 4r\right) ^{r}%
\leb%
\left( I\right) ^{-1/2}  \label{phi_I_j_norme_infinie} \\
& \leq \sqrt{2}c_{M,%
\leb%
}^{-1}c_{\min }^{-1/2}\left( 4r\right) ^{r}\left( r+1\right) ^{-1/2}\sqrt{D}
\label{phi_I_j_norme_infinie_end}
\end{align}%
where in the last inequality we used the fact that 
\begin{equation*}
\sqrt{\left\vert \mathcal{P}\right\vert \inf_{I\in \mathcal{P}}%
\leb%
\left( I\right) }\geq c_{M,%
\leb%
}\text{ \ and \ }D=\left( r+1\right) \left\vert \mathcal{P}\right\vert \text{
}.
\end{equation*}%
For all $j\in \left\{ 0,...,r\right\} $, $\varphi _{I,j}$ is supported by
the element $I$ of $\mathcal{P}$, hence we deduce from (\ref%
{phi_I_j_norme_infinie}) that the orthonormal basis $\left\{ \varphi _{I,j},%
\text{ }I\in \mathcal{P},\text{ }j\in \left\{ 0,...,r\right\} \right\} $ of $%
\left( M,\left\Vert 
\cdot%
\right\Vert _{2}\right) $ satisfies (\ref%
{inequality_localized_leb_I_piecewise}) with 
\begin{equation*}
L_{r,c_{\min }}=\sqrt{2}c_{\min }^{-1/2}\left( 4r\right) ^{r}\text{ }.
\end{equation*}%
To conclude, observe that%
\begin{align*}
\left\Vert \sum_{I,j}\beta _{I,j}\varphi _{I,j}\right\Vert _{\infty }&
=\max_{I\in \mathcal{P}}\left\{ \left\Vert \sum_{j=0}^{r}\beta _{I,j}\varphi
_{I,j}\right\Vert _{\infty }\right\} \\
& \leq \left\vert \beta \right\vert _{\infty }\max_{I\in \mathcal{P}}\left\{
\sum_{j=0}^{r}\left\Vert \varphi _{I,j}\right\Vert _{\infty }\right\} \\
& \leq \left( r+1\right) \left\vert \beta \right\vert _{\infty }\max_{I\in 
\mathcal{P}}\max_{j\in \left\{ 0,...,r\right\} }\left\{ \left\Vert \varphi
_{I,j}\right\Vert _{\infty }\right\}
\end{align*}%
and thus, by plugging (\ref{phi_I_j_norme_infinie_end}) into the right-hand
side of the last inequality, we finally obtain that the value 
\begin{equation*}
L_{r,c_{\min },c_{M,%
\leb%
}}=\sqrt{2}c_{M,%
\leb%
}^{-1}c_{\min }^{-1/2}\left( 4r\right) ^{r}\left( r+1\right) ^{1/2}
\end{equation*}%
gives the desired bound (\ref{inequality_localized_piecewise}).\ 
$\blacksquare$%

\noindent We now turn to the proof of (\ref{rate_sup_norm_localized}) under
the assumptions of Lemma \ref{lemma_piecewise_polynomials_consistency}. The
proof is based on concentration inequalities recalled in Section \ref%
{section_probabilistic_tools} and on inequality (\ref%
{inequality_localized_leb_I_piecewise}) of Lemma \ref%
{lemma_piecewise_polynomials_localized}, that allows us to control the
sup-norm of elements of an orthonormal basis for a model of piecewise
polynomials.

\subsection*{}%
\textbf{Proof of Lemma \ref{lemma_piecewise_polynomials_consistency}.} Let $%
\alpha >0$ be fixed and $\gamma >0$ to be chosen later. The partition $%
\mathcal{P}$ associated to $M$ will be denoted by%
\begin{equation*}
\mathcal{P=}\left\{ I_{0},...,I_{m-1}\right\} \text{ ,}
\end{equation*}%
so that $\left\vert \mathcal{P}\right\vert =m$ and $D=\left( r+1\right) m$
where $D$ is the dimension of the model $M$. By (\ref%
{inequality_localized_leb_I_piecewise}) of Lemma \ref%
{lemma_piecewise_polynomials_localized} there exist an orthonormal basis $%
\left\{ \varphi _{I_{k},j};\text{ }k\in \left\{ 0,...,m-1\right\} ,\text{ }%
j\in \left\{ 0,...,r\right\} \right\} $ of $\left( M,L_{2}\left(
P^{X}\right) \right) $ such that,%
\begin{equation*}
\varphi _{I_{k},j}\text{ is supported by the element }I_{k}\text{ of }%
\mathcal{P}\text{, \ \ for all }j\in \left\{ 0,...,r\right\}
\end{equation*}%
and a constant $L_{r,c_{\min }}$ depending only on $r,c_{\min }$ and
satisfying 
\begin{equation}
\max_{j\in \left\{ 0,...,r\right\} }\left\Vert \varphi _{I_{k},j}\right\Vert
_{\infty }\leq L_{r,c_{\min }}\frac{1}{\sqrt{%
\leb%
\left( I_{k}\right) }}\text{, \ for all }k\in \left\{ 0,...,m-1\right\} 
\text{.}  \label{localized_piecewise_proof}
\end{equation}%
In order to avoid cumbersome notation, we define a total ordering $\preceq $
on the set 
\begin{equation*}
\mathcal{I=}\left\{ \left( I_{k},j\right) ;k\in \left\{ 0,...,m-1\right\} ,%
\text{ }j\in \left\{ 0,...,r\right\} \right\} \text{ ,}
\end{equation*}%
as follows. Let $\prec $ be a binary relation on $\mathcal{I\times I}$ such
that%
\begin{equation*}
\left( I_{k},j\right) \prec \left( I_{l},i\right) \text{ \ if \ }\left( k<l%
\text{ \ or \ }\left( k=l\text{ \ and \ }j<i\right) \right) \text{, }
\end{equation*}%
and consider the total ordering $\preceq $ defined to be%
\begin{equation*}
\left( I_{k},j\right) \preceq \left( I_{l},i\right) \text{ \ if \ }\left(
\left( I_{k},j\right) =\left( I_{l},i\right) \text{ \ or \ }\left(
I_{k},j\right) \prec \left( I_{l},i\right) \right) \text{ .}
\end{equation*}%
So, from the definition of $\preceq $, the vector $\beta =\left( \beta
_{I_{k},j}\right) _{\left( I_{k},j\right) \in \mathcal{I}}\in \mathbb{R}^{D}$
has coordinate $\beta _{I_{k},j}$ at position $\left( r+1\right) k+j+1$ and
the matrix 
\begin{equation*}
A=\left( A_{\left( I_{k},j\right) ,\left( I_{l},i\right) }\right) _{\left(
I_{k},j\right) ,\left( I_{l},i\right) \in \mathcal{I\times I}}\in \mathbb{R}%
^{D\times D}\text{ },
\end{equation*}%
has coefficient $A_{\left( I_{k},j\right) ,\left( I_{l},i\right) }$ at line $%
\left( r+1\right) k+j+1$ and column $\left( r+1\right) l+i+1$.

\noindent Now, for some $s=\sum_{\left( I_{k},j\right) \in \mathcal{I}}\beta
_{I_{k},j}\varphi _{I_{k},j}\in M$, we have%
\begin{align*}
& P_{n}\left( K\left( s\right) \right) =P_{n}\left[ \left( y-\left(
\sum_{\left( I_{k},j\right) \in \mathcal{I}}\beta _{I_{k},j}\varphi
_{I_{k},j}\left( x\right) \right) \right) ^{2}\right] \\
& =P_{n}y^{2}-2\sum_{\left( I_{k},j\right) \in \mathcal{I}}\beta
_{I_{k},j}P_{n}\left( y\varphi _{I_{k},j}\left( x\right) \right)
+\sum_{\left( I_{k},j\right) ,\left( I_{l},i\right) \in \mathcal{I\times I}%
}\beta _{I_{k},j}\beta _{I_{l},i}P_{n}\left( \varphi _{I_{k},j}\varphi
_{I_{l},i}\right) \text{ }.
\end{align*}%
Hence, by taking the derivative with respect to $\beta _{I_{k},j}$ in the
last quantity, 
\begin{align}
& \frac{1}{2}\frac{\partial }{\partial \beta _{I_{k},j}}P_{n}\left[ \left(
y-\left( \sum_{\left( I_{k},j\right) \in \mathcal{I}}\beta _{I_{k},j}\varphi
_{I_{k},j}\left( x\right) \right) \right) ^{2}\right]  \notag \\
& =-P_{n}\left( y\varphi _{I_{k},j}\left( x\right) \right) +\sum_{\left(
I_{l},i\right) \in \mathcal{I}}\beta _{I_{l},i}P_{n}\left( \varphi
_{I_{k},j}\varphi _{I_{l},i}\right) \text{ .}  \label{derivative}
\end{align}%
We see that if $\beta ^{\left( n\right) }=\left( \beta _{I_{k},j}^{\left(
n\right) }\right) _{\left( I_{k},j\right) \in \mathcal{I}}\in \mathbb{R}^{D}$
is a critical point of 
\begin{equation*}
P_{n}\left[ \left( y-\left( \sum_{\left( I_{k},j\right) \in \mathcal{I}%
}\beta _{I_{k},j}\varphi _{I_{k},j}\left( x\right) \right) \right) ^{2}%
\right] \text{ ,}
\end{equation*}%
it holds%
\begin{equation*}
\left( \frac{\partial }{\partial \beta _{I_{k},j}}P_{n}\left[ \left(
y-\left( \sum_{\left( I_{k},j\right) \in \mathcal{I}}\beta _{I_{k},j}\varphi
_{I_{k},j}\left( x\right) \right) \right) ^{2}\right] \right) \left( \beta
^{\left( n\right) }\right) =0
\end{equation*}%
and by combining (\ref{derivative}) with the fact that 
\begin{equation*}
P\left( \varphi _{I_{k},j}\right) ^{2}=1\text{ },\text{ for all }\left(
I_{k},j\right) \in \mathcal{I}\text{ \ \ \ and \ \ \ }P\left( \varphi
_{I_{k},j}\varphi _{I_{l},i}\right) =0\text{ \ if \ }\left( I_{k},j\right)
\neq \left( I_{l},i\right) \text{ ,}
\end{equation*}%
we deduce that $\beta ^{\left( n\right) }$ satisfies the following random
linear system,%
\begin{equation}
\left( I_{D}+L_{n,D}\right) \beta ^{\left( n\right) }=X_{y,n}
\label{system_piecewise}
\end{equation}%
where $X_{y,n}=\left( P_{n}\left( y\varphi _{I_{k},j}\left( x\right) \right)
\right) _{\left( I_{k},j\right) \in \mathcal{I}}\in \mathbb{R}^{D},$ $I_{D}$
is the identity matrix of dimension $D$ and

\noindent $L_{n,D}=\left( \left( L_{n,D}\right) _{\left( I_{k},j\right)
,\left( I_{l},i\right) }\right) _{\left( I_{k},j\right) ,\left(
I_{l},i\right) \in \mathcal{I\times I}}$ is a $D\times D$ matrix satisfying%
\begin{equation*}
\left( L_{n,D}\right) _{\left( I_{k},j\right) ,\left( I_{l},i\right)
}=\left( P_{n}-P\right) \left( \varphi _{I_{k},j}\varphi _{I_{l},i}\right) 
\text{ }.
\end{equation*}%
Now, by inequality (\ref{lemma_sub_2}) in Lemma \ref%
{lemma_sub_consistency_piecewise} below, one can find a positive integer $%
n_{0}\left( r,A_{+},c_{\min },c_{M,%
\leb%
},\gamma \right) $ such that for all $n\geq n_{0}$, we have on an event $%
\Omega _{n}$ of probability at least $1-3Dn^{-\gamma }$, 
\begin{equation}
\left\Vert L_{n,D}\right\Vert \leq \frac{1}{2}\text{ },
\label{control_L_n_D}
\end{equation}%
where for a $D\times D$ matrix $L$, the operator norm $\left\Vert 
\cdot%
\right\Vert $ associated to the sup-norm on vectors is%
\begin{equation*}
\left\Vert L\right\Vert =\sup_{x\neq 0}\frac{\left\vert Lx\right\vert
_{\infty }}{\left\vert x\right\vert _{\infty }}\text{ }.
\end{equation*}%
Then we deduce from (\ref{control_L_n_D}) that $\left( I_{D}+L_{n,D}\right) $
is a non-singular $D\times D$ matrix and, as a consequence, that the linear
system (\ref{system_piecewise}) admits a unique solution $\beta ^{\left(
n\right) }$ on $\Omega _{n}$ for any $n\geq n_{0}\left( r,A_{+},c_{\min
},c_{M,%
\leb%
},\gamma \right) $. Moreover, since $P_{n}\left( y-\left( \sum_{\left(
I_{k},j\right) \in \mathcal{I}}\beta _{I_{k},j}\varphi _{I_{k},j}\left(
x\right) \right) \right) ^{2}$ is a nonnegative quadratic functional with
respect to $\left( \beta _{I_{k},j}\right) _{\left( I_{k},j\right) \in 
\mathcal{I}}\in \mathbb{R}^{D}$ we can easily deduce that on $\Omega _{n}$, $%
\beta ^{\left( n\right) }$ achieves the unique minimum of $P_{n}\left(
y-\left( \sum_{\left( I_{k},j\right) \in \mathcal{I}}\beta _{I_{k},j}\varphi
_{I_{k},j}\left( x\right) \right) \right) ^{2}$ on $\mathbb{R}^{D}$. In
other words, 
\begin{equation*}
s_{n}=\sum_{\left( I_{k},j\right) \in \mathcal{I}}\beta _{I_{k},j}^{\left(
n\right) }\varphi _{I_{k},j}
\end{equation*}%
is the unique least squares estimator on $M$, and by (\ref{system_piecewise}%
) it holds,%
\begin{equation}
\beta _{I_{k},j}^{\left( n\right) }\left( 1+\sum_{\left( I_{l},i\right) \in 
\mathcal{I}}\left( P_{n}-P\right) \left( \varphi _{I_{k},j}\varphi
_{I_{l},i}\right) \right) =P_{n}\left( y\varphi _{I_{k},j}\left( x\right)
\right) \text{ , for all }\left( I_{k},j\right) \in \mathcal{I}\text{.}
\label{equation_beta_n_1}
\end{equation}%
Now, as $\varphi _{I_{k},j}$ and $\varphi _{I_{l},i}$ have disjoint supports
when $k\neq l$, it holds $\varphi _{I_{k},j}\varphi _{I_{l},i}=0$ whenever $%
k\neq l$, and so equation (\ref{equation_beta_n_1}) reduces to%
\begin{equation}
\beta _{I_{k},j}^{\left( n\right) }\times \left( 1+\sum_{i=0}^{r}\left(
P_{n}-P\right) \left( \varphi _{I_{k},j}\varphi _{I_{k},i}\right) \right)
=P_{n}\left( y\varphi _{I_{k},j}\left( x\right) \right) \text{ , \ for all }%
\left( I_{k},j\right) \in \mathcal{I}\text{ .}  \label{equation_beta_n_2}
\end{equation}%
Moreover, recalling that $s_{M}=\sum_{\left( I_{k},j\right) \in \mathcal{I}%
}P\left( y\varphi _{I_{k},j}\left( x\right) \right) \varphi _{I_{k},j}$ , it
holds%
\begin{align}
\left\Vert s_{n}-s_{M}\right\Vert _{\infty }& =\left\Vert \sum_{\left(
I_{k},j\right) \in \mathcal{I}}\left( \beta _{I_{k},j}^{\left( n\right)
}-P\left( y\varphi _{I_{k},j}\left( x\right) \right) \right) \varphi
_{I_{k},j}\right\Vert _{\infty }  \notag \\
& \leq \max_{k\in \left\{ 0,...,m-1\right\} }\left\Vert \sum_{j=0}^{r}\left(
\beta _{I_{k},j}^{\left( n\right) }-P\left( y\varphi _{I_{k},j}\left(
x\right) \right) \right) \varphi _{I_{k},j}\right\Vert _{\infty }  \notag \\
& \leq \left( r+1\right) \max_{k\in \left\{ 0,...,m-1\right\} }\left\{
\left( \max_{j\in \left\{ 0,...,r\right\} }\left\vert \beta
_{I_{k},j}^{\left( n\right) }-P\left( y\varphi _{I_{k},j}\left( x\right)
\right) \right\vert \right) \right.  \notag \\
& \text{ \ \ \ \ \ \ \ \ \ \ \ \ \ \ \ \ \ \ \ \ \ \ \ \ \ \ \ \ \ \ \ \ \ \
\ \ \ \ \ \ }\left. \times \max_{j\in \left\{ 0,...,r\right\} }\left\Vert
\varphi _{I_{k},j}\right\Vert _{\infty }\right\}
\label{control_norme_infinie_piecewise}
\end{align}%
where the first inequality comes from the fact that $\varphi _{I_{k},j}$ and 
$\varphi _{I_{l},i}$ have disjoint supports when $k\neq l$. We next turn to
the control of the right-hand side of (\ref{control_norme_infinie_piecewise}%
). Let the index $\left( I_{k},j\right) $ be fixed. By subtracting the
quantity $\left( 1+\sum_{i=0}^{r}\left( P_{n}-P\right) \left( \varphi
_{I_{k},j}\varphi _{I_{k},i}\right) \right) \times P\left( y\varphi
_{I_{k},j}\left( x\right) \right) $\ in each side of equation (\ref%
{equation_beta_n_2}), we get%
\begin{align}
& \left( \beta _{I_{k},j}^{\left( n\right) }-P\left( y\varphi
_{I_{k},j}\left( x\right) \right) \right) \times \left(
1+\sum_{i=0}^{r}\left( P_{n}-P\right) \left( \varphi _{I_{k},j}\varphi
_{I_{k},i}\right) \right)  \notag \\
& =\left( P_{n}-P\right) \left( y\varphi _{I_{k},j}\left( x\right) \right)
-\left( \sum_{i=0}^{r}\left( P_{n}-P\right) \left( \varphi _{I_{k},j}\varphi
_{I_{k},i}\right) \right) \times P\left( y\varphi _{I_{k},j}\left( x\right)
\right) \text{ }.  \label{equation_beta_py}
\end{align}%
Moreover, by Inequality (\ref{lemma_sub_1}) of Lemma \ref%
{lemma_sub_consistency_piecewise}, we have for all $n\geq n_{0}\left(
r,A_{+},c_{\min },c_{M,%
\leb%
},\gamma \right) $,%
\begin{equation}
\sum_{i=0}^{r}\left\vert \left( P_{n}-P\right) \left( \varphi
_{I_{k},j}\varphi _{I_{k},i}\right) \right\vert \leq L_{r,A_{+},c_{\min
},c_{M,%
\leb%
},\gamma }\sqrt{\frac{\ln n}{n%
\leb%
\left( I_{k}\right) }}\leq \frac{1}{2}  \label{inequality_indice_fixed}
\end{equation}%
on the event $\Omega _{n}$. We thus deduce that%
\begin{equation}
\left\vert \left( \beta _{I_{k},j}^{\left( n\right) }-P\left( y\varphi
_{I_{k},j}\left( x\right) \right) \right) \times \left(
1+\sum_{i=0}^{r}\left( P_{n}-P\right) \left( \varphi _{I_{k},j}\varphi
_{I_{k},i}\right) \right) \right\vert \geq \frac{1}{2}\left\vert \beta
_{I_{k},j}^{\left( n\right) }-P\left( y\varphi _{I_{k},j}\left( x\right)
\right) \right\vert  \label{control_beta_py_left}
\end{equation}%
and 
\begin{equation}
\left\vert \left( \sum_{i=0}^{r}\left( P_{n}-P\right) \left( \varphi
_{I_{k},j}\varphi _{I_{k},i}\right) \right) \times P\left( y\varphi
_{I_{k},j}\left( x\right) \right) \right\vert \leq L_{r,A_{+},c_{\min },c_{M,%
\leb%
},\gamma }\sqrt{\frac{\ln n}{n%
\leb%
\left( I_{k}\right) }}\times \left\vert P\left( y\varphi _{I_{k},j}\left(
x\right) \right) \right\vert \text{ .}  \label{control_beta_py_right_2}
\end{equation}%
Moreover, by (\ref{bounded_data}), (\ref{density_piecwise_consistency}) and (%
\ref{localized_piecewise_proof}) we have%
\begin{align}
\left\vert P\left( y\varphi _{I_{k},j}\left( x\right) \right) \right\vert &
\leq A\left\Vert \varphi _{I_{k},j}\right\Vert _{\infty }P\left( I_{k}\right)
\notag \\
& \leq Ac_{\max }\left\Vert \varphi _{I_{k},j}\right\Vert _{\infty }%
\leb%
\left( I_{k}\right)  \notag \\
& \leq Ac_{\max }L_{r,c_{\min }}\sqrt{%
\leb%
\left( I_{k}\right) }  \notag \\
& \leq L_{A,r,c_{\min },c_{\max }}\sqrt{%
\leb%
\left( I_{k}\right) }\text{ }.  \label{py}
\end{align}%
Putting inequality (\ref{py}) in (\ref{control_beta_py_right_2}) we obtain%
\begin{equation}
\left\vert \left( \sum_{i=0}^{r}\left( P_{n}-P\right) \left( \varphi
_{I_{k},j}\varphi _{I_{k},i}\right) \right) \times P\left( y\varphi
_{I_{k},j}\left( x\right) \right) \right\vert \leq L_{r,A_{+},c_{\min
},c_{\max },c_{M,%
\leb%
},\gamma }\sqrt{\frac{\ln n}{n}}\text{ .}
\label{control_beta_py_right_2_end}
\end{equation}%
Hence, using inequalities (\ref{control_beta_py_left}), (\ref%
{control_beta_py_right_2_end}) and inequality (\ref{lemma_sub_3}) of Lemma %
\ref{lemma_sub_consistency_piecewise} in equation (\ref{equation_beta_py}),
we obtain that%
\begin{equation*}
\left\vert \beta _{I_{k},j}^{\left( n\right) }-P\left( y\varphi
_{I_{k},j}\left( x\right) \right) \right\vert \leq L_{A,r,A_{+},c_{\min
},c_{\max },c_{M,%
\leb%
},\gamma }\sqrt{\frac{\ln n}{n}}
\end{equation*}%
on $\Omega _{n}$. Since the constant $L_{A,r,A_{+},c_{\min },c_{\max },c_{M,%
\leb%
},\gamma }$ does not depend on the index $\left( I_{k},j\right) $ we deduce
by (\ref{localized_piecewise_proof}) that%
\begin{align}
& \left( \max_{j\in \left\{ 0,...,r\right\} }\left\vert \beta
_{I_{k},j}^{\left( n\right) }-P\left( y\varphi _{I_{k},j}\left( x\right)
\right) \right\vert \right) \times \max_{j\in \left\{ 0,...,r\right\}
}\left\Vert \varphi _{I_{k},j}\right\Vert _{\infty }  \notag \\
& \leq L_{A,r,A_{+},c_{\min },c_{\max },c_{M,%
\leb%
},\gamma }\sqrt{\frac{\ln n}{n}}\times \max_{j\in \left\{ 0,...,r\right\}
}\left\Vert \varphi _{I_{k},j}\right\Vert _{\infty }  \notag \\
& \leq L_{A,r,A_{+},c_{\min },c_{\max },c_{M,%
\leb%
},\gamma }\sqrt{\frac{\ln n}{n%
\leb%
\left( I_{k}\right) }}\text{ .}  \label{controle_last_piecewise}
\end{align}%
Finally, by using (\ref{partition_piecewise_consistency}) and (\ref%
{controle_last_piecewise}) in (\ref{control_norme_infinie_piecewise}), we
get for all $n\geq n_{0}\left( r,A_{+},c_{\min },c_{M,%
\leb%
},\gamma \right) $, on the event $\Omega _{n}$ of probability at least $%
1-3Dn^{-\gamma }$,%
\begin{align*}
\left\Vert s_{n}-s_{M}\right\Vert _{\infty }& \leq \left( r+1\right)
\max_{k\in \left\{ 0,...,m-1\right\} }\left\{ \left( \max_{j\in \left\{
0,...,r\right\} }\left\vert \beta _{I_{k},j}^{\left( n\right) }-P\left(
y\varphi _{I_{k},j}\left( x\right) \right) \right\vert \right) \times
\max_{j\in \left\{ 0,...,r\right\} }\left\Vert \varphi _{I_{k},j}\right\Vert
_{\infty }\right\} \\
& \leq L_{A,r,A_{+},c_{\min },c_{M,%
\leb%
},\gamma }\sqrt{\frac{\ln n}{n}}\max_{k\in \left\{ 0,...,m-1\right\} }\frac{1%
}{\sqrt{%
\leb%
\left( I_{k}\right) }} \\
& \leq L_{A,r,A_{+},c_{\min },c_{M,%
\leb%
},\gamma }\sqrt{\frac{\left\vert \mathcal{P}\right\vert \ln n}{n}} \\
& \leq L_{A,r,A_{+},c_{\min },c_{M,%
\leb%
},\gamma }\sqrt{\frac{D\ln n}{n}}\text{ }.
\end{align*}%
To conclude, simply take $\gamma =\frac{\ln 3}{\ln 2}+\alpha +1$, so that it
holds for $n\geq 2$, $\mathbb{P}\left[ \Omega _{n}^{c}\right] \leq
n^{-\alpha }$ which implies (\ref{rate_sup_norm_localized}).

\noindent It remains to prove the following lemma that has been used all
along the proof.

\begin{lemma}
\label{lemma_sub_consistency_piecewise}Recall that $L_{n,D}=\left( \left(
L_{n,D}\right) _{\left( I_{k},j\right) ,\left( I_{l},i\right) }\right)
_{\left( I_{k},j\right) ,\left( I_{l},i\right) \in \mathcal{I\times I}}$ is
a $D\times D$ matrix such that for all $\left( k,l\right) \in \left\{
0,...,m-1\right\} ^{2}$ , $\left( j,i\right) \in \left\{ 0,...,r\right\}
^{2} $ ,%
\begin{equation*}
\left( L_{n,D}\right) _{\left( I_{k},j\right) ,\left( I_{l},i\right)
}=\left( P_{n}-P\right) \left( \varphi _{I_{k},j}\varphi _{I_{l},i}\right) 
\text{ .}
\end{equation*}%
Also recall that for a $D\times D$ matrix $L$, the operator norm $\left\Vert 
\cdot%
\right\Vert $ associated to the sup-norm on the vectors is%
\begin{equation*}
\left\Vert L\right\Vert =\sup_{x\neq 0}\frac{\left\vert Lx\right\vert
_{\infty }}{\left\vert x\right\vert _{\infty }}\text{ }.
\end{equation*}%
Then, under the assumptions of Lemma \ref%
{lemma_piecewise_polynomials_consistency}, a positive integer $n_{0}\left(
r,A_{+},c_{\min },c_{M,%
\leb%
},\gamma \right) $ exists such that, for all $n\geq n_{0}\left(
r,A_{+},c_{\min },c_{M,%
\leb%
},\gamma \right) $, the following inequalities hold on an event $\Omega _{n}$
of probability at least $1-3Dn^{-\gamma }$,%
\begin{equation}
\left\Vert L_{n,D}\right\Vert \leq L_{r,A_{+},c_{\min },c_{M,%
\leb%
},\gamma }\sqrt{\frac{D\ln n}{n}}\leq \frac{1}{2}  \label{lemma_sub_2}
\end{equation}%
and for all $k\in \left\{ 0,...,m-1\right\} ,$ 
\begin{equation}
\max_{j\in \left\{ 0,...,r\right\} }\left\{ \sum_{i=0}^{r}\left\vert \left(
P_{n}-P\right) \left( \varphi _{I_{k},j}\varphi _{I_{k},i}\right)
\right\vert \right\} \leq L_{r,A_{+},c_{\min },c_{M,%
\leb%
},\gamma }\sqrt{\frac{\ln n}{n%
\leb%
\left( I_{k}\right) }}\leq \frac{1}{2}\text{ ,}  \label{lemma_sub_1}
\end{equation}%
\begin{equation}
\max_{j\in \left\{ 0,...,r\right\} }\left\vert \left( P_{n}-P\right) \left(
y\varphi _{I_{k},j}\left( x\right) \right) \right\vert \leq
L_{A,A_{+},r,c_{\min },c_{M,%
\leb%
},\gamma }\sqrt{\frac{\ln n}{n}}\text{ .}  \label{lemma_sub_3}
\end{equation}
\end{lemma}

\noindent \textbf{Proof of Lemma \ref{lemma_sub_consistency_piecewise}. }Let
us begin with the proof of inequality (\ref{lemma_sub_3}). Let the index $%
\left( I_{k},j\right) \in \mathcal{I}$ be fixed. By using Bernstein's
inequality (\ref{bernstein_ineq_reg}) and observing that, by (\ref%
{bounded_data}), 
\begin{equation*}
\var%
\left( y\varphi _{I_{k},j}\left( x\right) \right) \leq P\left[ \left(
y\varphi _{I_{k},j}\left( x\right) \right) ^{2}\right] \leq \left\Vert
Y\right\Vert _{\infty }^{2}\leq A^{2}
\end{equation*}%
and, by (\ref{bounded_data}), (\ref{localized_piecewise_proof}) and (\ref%
{partition_piecewise_consistency}),%
\begin{align*}
\left\Vert Y\varphi _{I_{k},j}\left( X\right) \right\Vert _{\infty }& \leq
A\left\Vert \varphi _{I_{k},j}\left( X\right) \right\Vert _{\infty } \\
& \leq AL_{r,c_{\min }}\frac{1}{\sqrt{%
\leb%
\left( I_{k}\right) }} \\
& \leq L_{A,r,c_{\min },c_{M,%
\leb%
}}\sqrt{\left\vert \mathcal{P}\right\vert } \\
& \leq L_{A,r,c_{\min },c_{M,%
\leb%
}}\sqrt{D}\text{ ,}
\end{align*}%
we get%
\begin{equation}
\mathbb{P}\left[ \left\vert \left( P_{n}-P\right) \left( y\varphi
_{I_{k},j}\left( x\right) \right) \right\vert \geq \sqrt{2A^{2}\frac{x}{n}}+%
\frac{L_{A,r,c_{\min },c_{M,%
\leb%
}}\sqrt{D}}{3n}x\right] \leq 2\exp \left( -x\right) \text{ }.  \label{bern}
\end{equation}%
By taking $x=\gamma \ln n$ in inequality (\ref{bern}), we obtain that%
\begin{equation}
\mathbb{P}\left[ \left\vert \left( P_{n}-P\right) \left( y\varphi
_{I_{k},j}\left( x\right) \right) \right\vert \geq \sqrt{2A^{2}\gamma \frac{%
\ln n}{n}}+\frac{L_{A,r,c_{\min },c_{M,%
\leb%
}}\sqrt{D}\gamma \ln n}{3n}\right] \leq 2n^{-\gamma }\text{ }.
\label{ineq_bernstein_py}
\end{equation}%
Now, as $D\leq A_{+}n\left( \ln n\right) ^{-2}$, we deduce from (\ref%
{ineq_bernstein_py}) that for some well chosen positive constant $%
L_{A,A_{+},r,c_{\min },c_{M,%
\leb%
},\gamma }$, we have%
\begin{equation*}
\mathbb{P}\left[ \left\vert \left( P_{n}-P\right) \left( y\varphi
_{I_{k},j}\left( x\right) \right) \right\vert \geq L_{A,A_{+},r,c_{\min
},c_{M,%
\leb%
},\gamma }\sqrt{\frac{\ln n}{n}}\right] \leq 2n^{-\gamma }\text{ }
\end{equation*}%
and by setting 
\begin{equation*}
\Omega _{n}^{\left( 1\right) }=\bigcap\limits_{\left( I_{k},j\right) \in 
\mathcal{I}}\left\{ \left\vert \left( P_{n}-P\right) \left( y\varphi
_{I_{k},j}\left( x\right) \right) \right\vert \leq L_{A,A_{+},r,c_{\min
},c_{M,%
\leb%
},\gamma }\sqrt{\frac{\ln n}{n}}\right\}
\end{equation*}%
we deduce that 
\begin{equation}
\mathbb{P}\left( \Omega _{n}^{\left( 1\right) }\right) \geq 1-2Dn^{-\gamma }%
\text{ .}  \label{proba_omega_1_n}
\end{equation}%
Hence the expected bound (\ref{lemma_sub_3}) holds on $\Omega _{n}^{\left(
1\right) }$, for all $n\geq 1$.

\noindent We turn now to the proof of inequality (\ref{lemma_sub_1}). Let
the index $\left( I_{k},j\right) \in \mathcal{I}$ be fixed. By
Cauchy-Schwarz inequality, we have%
\begin{equation}
\sum_{i=0}^{r}\left\vert \left( P_{n}-P\right) \left( \varphi
_{I_{k},j}\varphi _{I_{k},i}\right) \right\vert \leq \sqrt{r+1}\sqrt{%
\sum_{i=0}^{r}\left( \left( P_{n}-P\right) \left( \varphi _{I_{k},j}\varphi
_{I_{k},i}\right) \right) ^{2}}\text{ .}  \label{CS_khi}
\end{equation}%
Let write 
\begin{equation*}
\chi _{I_{k},j}=\sqrt{\sum_{i=0}^{r}\left( \left( P_{n}-P\right) \left(
\varphi _{I_{k},j}\varphi _{I_{k},i}\right) \right) ^{2}}\text{ \ and \ }%
B_{I_{k}}=\left\{ \sum_{i=0}^{r}\beta _{I_{k},i}\varphi _{I_{k},i}\text{ };%
\text{ }\left( \beta _{I_{k},i}\right) _{i=0}^{r}\in \mathbb{R}^{r+1}\text{
and }\sum_{i=0}^{r}\beta _{I_{k},i}^{2}\leq 1\right\} \text{ .}
\end{equation*}%
By Cauchy-Schwarz inequality again, it holds%
\begin{equation*}
\chi _{I_{k},j}=\sup_{s\in B_{I_{k}}}\left\vert \left( P_{n}-P\right) \left(
\varphi _{I_{k},j}s\right) \right\vert \text{ .}
\end{equation*}%
Then, Bousquet's inequality (\ref{bousquet}), applied with $\varepsilon =1$
and $\mathcal{F=}B_{I_{k}}$, implies that%
\begin{equation}
\mathbb{P}\left[ \chi _{I_{k},j}-\mathbb{E}\left[ \chi _{I_{k},j}\right]
\geq \sqrt{2\sigma _{I_{k},j}^{2}\frac{x}{n}}+\mathbb{E}\left[ \chi
_{I_{k},j}\right] +\frac{4}{3}\frac{b_{I_{k},j}x}{n}\right] \leq \exp \left(
-x\right)  \label{concentration_khi_ik_j}
\end{equation}%
where, by (\ref{localized_piecewise_proof}), 
\begin{equation}
\sigma _{I_{k},j}^{2}=\sup_{s\in B_{I_{k}}}%
\var%
\left( \varphi _{I_{k},j}s\right) \leq \left\Vert \varphi
_{I_{k},j}\right\Vert _{\infty }^{2}\leq \frac{L_{r,c_{\min }}}{%
\leb%
\left( I_{k}\right) }  \label{variance_khi_Ik}
\end{equation}%
and 
\begin{equation}
b_{I_{k},j}\leq 2\sup_{s\in B_{I_{k}}}\left\Vert \varphi
_{I_{k},j}s\right\Vert _{\infty }\leq 2\left\Vert \varphi
_{I_{k},j}\right\Vert _{\infty }\sup_{s\in B_{I_{k}}}\left\Vert s\right\Vert
_{\infty }\text{ }.  \label{b_Ik_j_1}
\end{equation}%
Moreover, for $s=\sum_{i=0}^{r}\beta _{I_{k},i}\varphi _{I_{k},i}\in
B_{I_{k}}$, we have $\max_{i}\left\vert \beta _{I_{k},i}\right\vert \leq 
\sqrt{\sum_{i=0}^{r}\beta _{I_{k},i}^{2}}\leq 1$, so by (\ref%
{localized_piecewise_proof}),%
\begin{equation*}
\sup_{s\in B_{I_{k}}}\left\Vert s\right\Vert _{\infty }\leq
\sum_{i=0}^{r}\left\Vert \varphi _{I_{k},i}\right\Vert _{\infty }\leq \frac{%
L_{r,c_{\min }}}{\sqrt{%
\leb%
\left( I_{k}\right) }}
\end{equation*}%
and injecting the last bound in (\ref{b_Ik_j_1}) we get%
\begin{equation}
b_{I_{k},j}\leq \left\Vert \varphi _{I_{k},j}\right\Vert _{\infty }\frac{%
L_{r,c_{\min }}}{\sqrt{%
\leb%
\left( I_{k}\right) }}\leq \frac{L_{r,c_{\min }}}{%
\leb%
\left( I_{k}\right) }\text{ }.  \label{b_Ik_j_1_end}
\end{equation}%
In addition, we have 
\begin{align}
\mathbb{E}\left[ \chi _{I_{k},j}\right] & \leq \sqrt{\mathbb{E}\left[ \chi
_{I_{k},j}^{2}\right] }=\sqrt{\frac{\sum_{i=0}^{r}%
\var%
\left( \varphi _{I_{k},j}\varphi _{I_{k},i}\right) }{n}}  \notag \\
& \leq \left\Vert \varphi _{I_{k},j}\right\Vert _{\infty }\sqrt{\frac{%
\sum_{i=0}^{r}P\left( \varphi _{I_{k},i}^{2}\right) }{n}}  \notag \\
& =\left\Vert \varphi _{I_{k},j}\right\Vert _{\infty }\sqrt{\frac{r+1}{n}} 
\notag \\
& \leq L_{r,c_{\min }}\sqrt{\frac{1}{n%
\leb%
\left( I_{k}\right) }}\text{ }.  \label{mean_khi}
\end{align}%
Therefore, combining (\ref{variance_khi_Ik}), (\ref{b_Ik_j_1_end}), (\ref%
{mean_khi}) and (\ref{concentration_khi_ik_j}) while taking $x=\gamma \ln n$%
, we get%
\begin{equation}
\mathbb{P}\left[ \chi _{I_{k},j}\geq L_{r,c_{\min },\gamma }\left( \sqrt{%
\frac{1}{n%
\leb%
\left( I_{k}\right) }}+\sqrt{\frac{\ln n}{n%
\leb%
\left( I_{k}\right) }}+\frac{\ln n}{n%
\leb%
\left( I_{k}\right) }\right) \right] \leq n^{-\gamma }\text{ .}
\label{concentration_khi_ik_j_2}
\end{equation}%
Now, since by (\ref{partition_piecewise_consistency}) and the fact that $%
D\leq A_{+}n\left( \ln n\right) ^{-2}$ we have%
\begin{equation*}
\frac{1}{%
\leb%
\left( I_{k}\right) }\leq c_{M,%
\leb%
}^{-2}D\leq c_{M,%
\leb%
}^{-2}A_{+}\frac{n}{\left( \ln n\right) ^{2}}\text{ ,}
\end{equation*}%
we obtain from (\ref{concentration_khi_ik_j_2}) that a positive constant $%
L_{r,A_{+},c_{\min },c_{M,%
\leb%
},\gamma }$ exists, depending only on $\gamma ,r,A_{+},c_{\min }$ and $c_{M,%
\leb%
}$ such that%
\begin{equation}
\mathbb{P}\left[ \chi _{I_{k},j}\geq L_{r,A_{+},c_{\min },c_{M,%
\leb%
},\gamma }\sqrt{\frac{\ln n}{n%
\leb%
\left( I_{k}\right) }}\right] \leq n^{-\gamma }\text{ .}
\label{concentration_khi_ik_j_3}
\end{equation}%
Finally, define%
\begin{equation*}
\Omega _{n}^{\left( 2\right) }=\bigcap\limits_{\left( I_{k},j\right) \in 
\mathcal{I}}\left\{ \chi _{I_{k},j}\leq L_{r,A_{+},c_{\min },c_{M,%
\leb%
},\gamma }\sqrt{\frac{\ln n}{n%
\leb%
\left( I_{k}\right) }}\right\} \text{ .}
\end{equation*}%
For all $n\geq n_{0}\left( r,A_{+},c_{\min },c_{M,%
\leb%
},\gamma \right) $, we have%
\begin{align}
& \sqrt{r+1}\times L_{r,A_{+},c_{\min },c_{M,%
\leb%
},\gamma }\sqrt{\frac{\ln n}{n%
\leb%
\left( I_{k}\right) }}  \notag \\
& \leq L_{r,A_{+},c_{\min },c_{M,%
\leb%
},\gamma }\sqrt{\frac{D\ln n}{n}}  \notag \\
& \leq L_{r,A_{+},c_{\min },c_{M,%
\leb%
},\gamma }\frac{1}{\sqrt{\ln n}}\leq \frac{1}{2}\text{ .}  \label{def_n0}
\end{align}%
Moreover by (\ref{concentration_khi_ik_j_3}) it holds 
\begin{equation}
\mathbb{P}\left( \Omega _{n}^{\left( 2\right) }\right) \geq 1-Dn^{-\gamma }
\label{omega_2}
\end{equation}%
and, by (\ref{CS_khi}), the expected bound (\ref{lemma_sub_1}) holds on $%
\Omega _{n}^{\left( 2\right) }$, for all $n\geq n_{0}\left( r,A_{+},c_{\min
},c_{M,%
\leb%
},\gamma \right) $.

\noindent Next, notice that for a $D\times D$ matrix $L=\left( L_{\left(
I_{k},j\right) ,\left( I_{l},i\right) }\right) _{\left( I_{k},j\right)
,\left( I_{l},i\right) \in\mathcal{I\times I}}$ we have the following
classical formula,%
\begin{equation*}
\left\Vert L\right\Vert =\max_{\left( I_{k},j\right) \in\mathcal{I}%
}\sum_{\left( I_{l},i\right) \in\mathcal{I}}\left\vert L_{\left(
I_{k},j\right) ,\left( I_{l},i\right) }\right\vert \text{ }.
\end{equation*}
Applied to the matrix of interest $L_{n,D}$ , this gives%
\begin{align}
\left\Vert L_{n,D}\right\Vert & =\max_{\left( I_{k},j\right) \in \mathcal{I}%
}\sum_{\left( I_{l},i\right) \in\mathcal{I}}\left\vert \left( P_{n}-P\right)
\left( \varphi_{I_{k},j}\varphi_{I_{l},i}\right) \right\vert \text{ }  \notag
\\
& =\max_{k\in\left\{ 0,...,m-1\right\} }\max_{j\in\left\{ 0,...,r\right\}
}\left\{ \sum_{\left( I_{l},i\right) \in\mathcal{I}}\left\vert \left(
P_{n}-P\right) \left( \varphi_{I_{k},j}\varphi_{I_{l},i}\right) \right\vert
\right\} \text{ }.  \label{formula_matrix}
\end{align}
Thus, using formula (\ref{formula_matrix}), inequalities (\ref{lemma_sub_1}%
), (\ref{partition_piecewise_consistency}) and (\ref{def_n0}) give that for
all $n\geq n_{0}\left( r,A_{+},c_{\min},c_{M,%
\leb%
},\gamma\right) $, we have on $\Omega_{n}^{\left( 2\right) }$,%
\begin{equation*}
\left\Vert L_{n,D}\right\Vert \leq L_{r,A_{+},c_{\min},c_{M,%
\leb%
},\gamma}\sqrt{\frac{D\ln n}{n}}\leq\frac{1}{2}\text{ }.
\end{equation*}
Finally, by setting $\Omega_{n}=\Omega_{n}^{\left( 1\right) }\bigcap
\Omega_{n}^{\left( 2\right) }$, we have $\mathbb{P}\left( \Omega _{n}\right)
\geq1-3Dn^{-\gamma}$, and inequalities (\ref{lemma_sub_1}), (\ref%
{lemma_sub_2}) and (\ref{lemma_sub_3}) are satisfied on $\Omega_{n}$ for all 
$n\geq n_{0}\left( r,A_{+},c_{\min},c_{M,%
\leb%
},\gamma\right) $, which completes the proof of Lemma \ref%
{lemma_sub_consistency_piecewise}.\ 
$\blacksquare$%

\subsection{Proofs of Section \protect\ref{section_reg_fixed_model} \label%
{section_proof_fixed_model}}

In order to express the quantities of interest in the proofs of Theorems \ref%
{principal_reg} and \ref{theorem_upper_true_small_models}, we need
preliminary definitions. Let $\alpha >0$ be fixed and for $R_{n,D,\alpha }$
defined in (\textbf{H5}), see Section \ref{section_main_assump_reg_fixed},
we set%
\begin{equation}
\tilde{R}_{n,D,\alpha }=\max \left\{ R_{n,D,\alpha }\text{ };\text{ }%
A_{\infty }\sqrt{\frac{D\ln n}{n}}\right\}  \label{def_R_n_D_tild}
\end{equation}%
where $A_{\infty }$ is a positive constant to be chosen later. Moreover, we
set 
\begin{equation}
\nu _{n}=\max \left\{ \sqrt{\frac{\ln n}{D}}\text{ };\text{ }\sqrt{\frac{%
D\ln n}{n}}\text{ };\text{ }R_{n,D,\alpha }\right\} \text{ }.
\label{nu_n_proof}
\end{equation}%
Thanks to the assumption of consistency in sup-norm (\textbf{H5}), our
analysis will be localized in the subset 
\begin{equation*}
B_{\left( M,L_{\infty }\right) }\left( s_{M},\tilde{R}_{n,D,\alpha }\right)
=\left\{ s\in M,\left\Vert s-s_{M}\right\Vert _{\infty }\leq \tilde{R}%
_{n,D,\alpha }\right\}
\end{equation*}%
of $M$.

\noindent Let us define several slices of excess risk on the model $M$ : for
any $C\geq 0$,%
\begin{align*}
\mathcal{F}_{C}& =\left\{ s\in M,P\left( Ks-Ks_{M}\right) \leq C\right\}
\bigcap B_{\left( M,L_{\infty }\right) }\left( s_{M},\tilde{R}_{n,D,\alpha
}\right) \\
\mathcal{F}_{>C}& =\left\{ s\in M,P\left( Ks-Ks_{M}\right) >C\right\}
\bigcap B_{\left( M,L_{\infty }\right) }\left( s_{M},\tilde{R}_{n,D,\alpha
}\right)
\end{align*}%
and for any interval $J\subset \mathbb{R},$%
\begin{equation*}
\mathcal{F}_{J}=\left\{ s\in M,P\left( Ks-Ks_{M}\right) \in J\right\}
\bigcap B_{\left( M,L_{\infty }\right) }\left( s_{M},\tilde{R}_{n,D,\alpha
}\right) \text{ }.
\end{equation*}%
We also define, for all $L\geq 0$,%
\begin{equation*}
D_{L}=\left\{ s\in M,P\left( Ks-Ks_{M}\right) =L\right\} \bigcap B_{\left(
M,L_{\infty }\right) }\left( s_{M},\tilde{R}_{n,D,\alpha }\right) \text{ }.
\end{equation*}%
Recall that, by Lemma \ref{dev_centered_sM copy(1)_reg} of Section \ref%
{section_excess_risk}, the contrasted functions satisfy, for every $s\in M$
and $z=\left( x,y\right) \in \mathcal{X\times }\mathbb{R},$%
\begin{equation*}
\left( Ks\right) \left( z\right) -\left( Ks_{M}\right) \left( z\right) =\psi
_{1,M}\left( z\right) \left( s-s_{M}\right) \left( x\right) +\psi _{2}\left(
\left( s-s_{M}\right) \left( x\right) \right)
\end{equation*}%
where $\psi _{1,M}\left( z\right) =-2\left( y-s_{M}\left( x\right) \right) $
and $\psi _{2}\left( t\right) =t^{2}$, for all $t\in \mathbb{R}.$ For
convenience, we will use the following notation, for any $s\in M$,%
\begin{equation*}
\psi _{2}\circ \left( s-s_{M}\right) :x\in \mathcal{X\longmapsto }\psi
_{2}\left( \left( s-s_{M}\right) \left( x\right) \right) \text{ .}
\end{equation*}%
Note that, for all $s\in M$,%
\begin{equation}
P\left( \psi _{1,M}%
\cdot%
s\right) =0  \label{center_proof}
\end{equation}%
and by (\textbf{H1}) inequality (\ref{control_phi_1M}) holds true, that is 
\begin{equation}
\left\Vert \psi _{1,M}\right\Vert _{\infty }\leq 4A\text{ }.
\label{control_phi_1M_proof}
\end{equation}%
Also, for $\mathcal{K}_{1,M}$ defined in Section \ref%
{section_complexity_model}, we have%
\begin{equation*}
\mathcal{K}_{1,M}=\sqrt{\frac{1}{D}\sum_{k=1}^{D}%
\var%
\left( \psi _{1,M}%
\cdot%
\varphi _{k}\right) }
\end{equation*}%
for any orthonormal basis $\left( \varphi _{k}\right) _{k=1}^{D}$ of $\left(
M,\left\Vert 
\cdot%
\right\Vert _{2}\right) .$ Moreover, inequality (\ref{upper_K1M}) holds
under (\textbf{H1}) and we have%
\begin{equation}
\mathcal{K}_{1,M}\leq 2\sigma _{\max }+4A\leq 6A\text{ .}
\label{upper_K1M_proof}
\end{equation}%
Assuming (\textbf{H2}), we have from (\ref{lower_K1M})%
\begin{equation}
0<2\sigma _{\min }\leq \mathcal{K}_{1,M}\text{ .}  \label{lower_K1M_proof}
\end{equation}%
Finally, when (\textbf{H3}) holds (it is the case when (\textbf{H4}) holds),
we have by (\ref{def_enveloppe}),%
\begin{equation}
\sup_{s\in M,\text{ }\left\Vert s\right\Vert _{2}\leq 1}\left\Vert
s\right\Vert _{\infty }\leq A_{3,M}\sqrt{D}  \label{control_boule_unite}
\end{equation}%
and so, for any orthonormal basis $\left( \varphi _{k}\right) _{k=1}^{D}$ of 
$\left( M,\left\Vert 
\cdot%
\right\Vert _{2}\right) $, it holds for all $k\in \left\{ 1,...,D\right\} $,
as $P\left( \varphi _{k}^{2}\right) =1$,%
\begin{equation}
\left\Vert \varphi _{k}\right\Vert _{\infty }\leq A_{3,M}\sqrt{D}\text{ .}
\label{control_basis}
\end{equation}

\subsubsection{Proofs of the theorems\label%
{section_proofs_theorem_fixed_mdel}}

The proof of Theorem \ref{principal_reg}\textbf{\ }relies on Lemmas \ref%
{majo_s1c_lemma copy(1)}, \ref{majo_s1>c_lemma} and \ref{mino_S_r_C} stated
in Section \ref{useful_lemmas}, and that give sharp estimates of suprema of
the empirical process on the contrasted functions over slices of interest.

\subsection*{}%
\textbf{Proof of Theorem \ref{principal_reg}. }Let $\alpha >0$ be fixed and
let $\varphi =\left( \varphi _{k}\right) _{k=1}^{D}$ be an orthonormal basis
of $\left( M,\left\Vert 
\cdot%
\right\Vert _{2}\right) $ satisfying (\textbf{H4}). We divide the proof of
Theorem \ref{principal_reg} into four parts, corresponding to the four
Inequalities (\ref{lower_true_risk}), (\ref{upper_true_risk}), (\ref%
{lower_emp_risk}) and (\ref{upper_emp_risk}). The values of $A_{0}$ and $%
A_{\infty }$, respectively defined in (\ref{def_A_0}) and (\ref%
{def_R_n_D_tild}), will then be chosen at the end of the proof.

\bigskip

\subsection*{}%
\textbf{Proof of Inequality (\ref{lower_true_risk}). }Let $r\in \left( 1,2%
\right] $ to be chosen later and $C>0$ such that%
\begin{equation}
rC=\frac{D}{4n}\mathcal{K}_{1,M}^{2}\text{ }.  \label{rC}
\end{equation}%
By (\textbf{H5}) there exists a positive integer $n_{1}$ such that it holds,
for all $n\geq n_{1}$,%
\begin{equation}
\mathbb{P}\left( P\left( Ks_{n}-Ks_{M}\right) \leq C\right) \leq \mathbb{P}%
\left( \left\{ P\left( Ks_{n}-Ks_{M}\right) \leq C\right\} \bigcap \Omega
_{\infty ,\alpha }\right) +n^{-\alpha }  \label{localisation_1}
\end{equation}%
and also 
\begin{align}
& \mathbb{P}\left( \left\{ P\left( Ks_{n}-Ks_{M}\right) \leq C\right\}
\bigcap \Omega _{\infty ,\alpha }\right)  \notag \\
& \leq \mathbb{P}\left( \inf_{s\in \mathcal{F}_{C}}P_{n}\left(
Ks-Ks_{M}\right) \leq \inf_{s\in \mathcal{F}_{>C}}P_{n}\left(
Ks-Ks_{M}\right) \right)  \notag \\
& \leq \mathbb{P}\left( \inf_{s\in \mathcal{F}_{C}}P_{n}\left(
Ks-Ks_{M}\right) \leq \inf_{s\in \mathcal{F}_{\left( C,rC\right]
}}P_{n}\left( Ks-Ks_{M}\right) \right)  \notag \\
& =\mathbb{P}\left( \sup_{s\in \mathcal{F}_{C}}P_{n}\left( Ks_{M}-Ks\right)
\geq \sup_{s\in \mathcal{F}_{\left( C,rC\right] }}P_{n}\left(
Ks_{M}-Ks\right) \right) .  \label{lower_1}
\end{align}%
Now, by (\ref{rC}) and (\ref{lower_K1M_proof}) we have%
\begin{equation*}
\frac{D}{2n}\sigma _{\min }^{2}\leq C\leq \left( 1+A_{4}\nu _{n}\right) ^{2}%
\frac{D}{4n}\mathcal{K}_{1,M}^{2}\text{ }
\end{equation*}%
where $A_{4}$ is defined in Lemma \ref{majo_s1c_lemma copy(1)}. Hence we can
apply Lemma \ref{majo_s1c_lemma copy(1)} with $\alpha =\beta $, $%
A_{l}=\sigma _{\min }^{2}/2$ and $A_{3,M}=r_{M}\left( \varphi \right) $, by
Remark 3. Therefore it holds, for all $n\geq n_{0}\left( A_{\infty
},A_{cons},A_{+},\sigma _{\min },\alpha \right) $,%
\begin{equation}
\mathbb{P}\left[ \sup_{s\in \mathcal{F}_{C}}P_{n}\left( Ks_{M}-Ks\right)
\geq \left( 1+L_{A_{\infty },A,r_{M}\left( \varphi \right) ,\sigma _{\min
},A_{-},\alpha }\times \nu _{n}\right) \sqrt{\frac{CD}{n}}\mathcal{K}_{1,M}-C%
\right] \leq 2n^{-\alpha }\text{ }.  \label{majo_SC_1}
\end{equation}%
Moreover, by using (\ref{lower_K1M_proof}) and (\ref{upper_K1M_proof}) in (%
\ref{rC}) we get%
\begin{equation*}
\frac{D}{n}\sigma _{\min }^{2}\leq rC\leq \frac{D}{n}\left( \sigma _{\max
}+2A\right) ^{2}\text{ }.
\end{equation*}%
We then apply Lemma \ref{mino_S_r_C} with 
\begin{equation*}
\alpha =\beta ,\text{ }A_{l}=\sigma _{\min }^{2},\text{ }A_{u}=\left( \sigma
_{\max }+2A\right) ^{2}
\end{equation*}%
and 
\begin{equation}
\text{ }A_{\infty }\geq 64\sqrt{2}B_{2}A\left( \sigma _{\max }+2A\right)
\sigma _{\min }^{-1}r_{M}\left( \varphi \right) \text{ },
\label{def_A_infini_1}
\end{equation}%
so it holds for all $n\geq n_{0}\left( A_{-},A_{+},A,A_{\infty
},A_{cons},B_{2},r_{M}\left( \varphi \right) ,\sigma _{\max },\sigma _{\min
},\alpha \right) $,%
\begin{equation}
\mathbb{P}\left( \sup_{s\in \mathcal{F}_{\left( C,rC\right] }}P_{n}\left(
Ks_{M}-Ks\right) \leq \left( 1-L_{A_{-},A,A_{\infty },\sigma _{\max },\sigma
_{\min },r_{M}\left( \varphi \right) ,\alpha }\times \nu _{n}\right) \sqrt{%
\frac{rCD}{n}}\mathcal{K}_{1,M}-rC\right) \leq 2n^{-\alpha }\text{ }.
\label{mino_S_r_C_2}
\end{equation}%
Now, from (\ref{majo_SC_1}) and (\ref{mino_S_r_C_2}) we can find a positive
constant $\tilde{A}_{0}$, only depending on $A_{-},A,A_{\infty },$ $\sigma
_{\max },\sigma _{\min },r_{M}\left( \varphi \right) $ and $\alpha $, such
that for all $n\geq n_{0}\left( A_{-},A_{+},A,A_{\infty
},A_{cons},B_{2},r_{M}\left( \varphi \right) ,\sigma _{\max },\sigma _{\min
},\alpha \right) $, there exists an event of probability at least $%
1-4n^{-\alpha }$ on which%
\begin{equation}
\sup_{s\in \mathcal{F}_{C}}P_{n}\left( Ks_{M}-Ks\right) \leq \left( 1+\tilde{%
A}_{0}\nu _{n}\right) \sqrt{\frac{CD}{n}}\mathcal{K}_{1,M}-C
\label{majo_SC_2}
\end{equation}%
and 
\begin{equation}
\sup_{s\in \mathcal{F}_{\left( C,rC\right] }}P_{n}\left( Ks_{M}-Ks\right)
\geq \left( 1-\tilde{A}_{0}\nu _{n}\right) \sqrt{\frac{rCD}{n}}\mathcal{K}%
_{1,M}-rC\text{ }.  \label{mino_S_r_C_3}
\end{equation}%
Hence, from (\ref{majo_SC_2}) and (\ref{mino_S_r_C_3}) we deduce, using (\ref%
{localisation_1}) and (\ref{lower_1}), that if we choose $r\in \left( 1,2%
\right] $ such that%
\begin{equation}
\left( 1+\tilde{A}_{0}\nu _{n}\right) \sqrt{\frac{CD}{n}}\mathcal{K}%
_{1,M}-C<\left( 1-\tilde{A}_{0}\nu _{n}\right) \sqrt{\frac{rCD}{n}}\mathcal{K%
}_{1,M}-rC  \label{equation_K}
\end{equation}%
then, for all $n\geq n_{0}\left( A_{-},A_{+},A,A_{\infty
},A_{cons},B_{2},r_{M}\left( \varphi \right) ,\sigma _{\max },\sigma _{\min
},n_{1},\alpha \right) $ we have 
\begin{equation*}
P\left( Ks_{n}-Ks_{M}\right) \geq C
\end{equation*}%
with probability at least $1-5n^{-\alpha }$. Now, by (\ref{rC}) it holds 
\begin{equation*}
\sqrt{\frac{rCD}{n}}\mathcal{K}_{1,M}=2rC=\frac{1}{2}\frac{D}{n}\mathcal{K}%
_{1,M}^{2}\text{ ,}
\end{equation*}%
and as a consequence Inequality (\ref{equation_K}) is equivalent to%
\begin{equation}
\left( 1-2\tilde{A}_{0}\nu _{n}\right) r-2\left( 1+\tilde{A}_{0}\nu
_{n}\right) \sqrt{r}+1>0\text{ }.  \label{equation_K_2}
\end{equation}%
Moreover, we have by (\ref{nu_n_proof}) and (\textbf{H5}), for all $n\geq
n_{0}\left( A_{+},A_{-},A_{cons},\tilde{A}_{0},\alpha \right) $, 
\begin{equation}
\tilde{A}_{0}\nu _{n}\leq \frac{1}{4}  \label{majo_nu_n_end}
\end{equation}%
and so, for all $n\geq n_{0}\left( A_{+},A_{-},A_{cons},\tilde{A}_{0},\alpha
\right) $, simple computations involving (\ref{majo_nu_n_end}) show that by
taking 
\begin{equation}
r=1+48\sqrt{\tilde{A}_{0}\nu _{n}}  \label{def_r}
\end{equation}%
inequality (\ref{equation_K_2}) is satisfied. Notice that, for all $n\geq
n_{0}\left( A_{+},A_{-},A_{cons},\tilde{A}_{0},\alpha \right) $ we have $0<48%
\sqrt{\tilde{A}_{0}\nu _{n}}<1$, so that $r\in \left( 1,2\right) $. Finally,
we compute $C$ by (\ref{rC}) and (\ref{def_r}), in such a way that for all $%
n\geq n_{0}\left( A_{+},A_{-},A_{cons},\tilde{A}_{0},\alpha \right) $, 
\begin{equation}
C=\frac{rC}{r}=\frac{1}{1+48\sqrt{\tilde{A}_{0}\nu _{n}}}\frac{1}{4}\frac{D}{%
n}\mathcal{K}_{1,M}^{2}\geq \left( 1-48\sqrt{\tilde{A}_{0}\nu _{n}}\right) 
\frac{1}{4}\frac{D}{n}\mathcal{K}_{1,M}^{2}>0  \label{def_A_0_1}
\end{equation}%
which yields the result by noticing that the dependence on $\sigma _{\max }$
can be released in $n_{0}$ and $\tilde{A}_{0}$ since by (\textbf{H1})\ we
have $\sigma _{\max }\leq A$.

\vspace*{0.5in}

\subsection*{}%
\textbf{Proof of Inequality (\ref{upper_true_risk}).} Let $C>0$ and $\delta
\in \left( 0,\frac{1}{2}\right) $ to be chosen later in such a way that%
\begin{equation}
\left( 1-\delta \right) C=\frac{D}{4n}\mathcal{K}_{1,M}^{2}
\label{def_1-delta_C}
\end{equation}%
and 
\begin{equation}
C\geq \frac{1}{4}\left( 1+A_{5}\nu _{n}\right) ^{2}\frac{D}{n}\mathcal{K}%
_{1,M}^{2}\text{ ,}  \label{C}
\end{equation}%
where $A_{5}$ is defined in Lemma \ref{majo_s1>c_lemma}. We have by (\textbf{%
H5}), for all $n\geq n_{1},$%
\begin{equation}
\mathbb{P}\left( P\left( Ks_{n}-Ks_{M}\right) >C\right) \leq \mathbb{P}%
\left( \left\{ P\left( Ks_{n}-Ks_{M}\right) >C\right\} \bigcap \Omega
_{\infty ,\alpha }\right) +n^{-\alpha }  \label{conditionnement_2}
\end{equation}%
and also%
\begin{align}
& \mathbb{P}\left( \left\{ P\left( Ks_{n}-Ks_{M}\right) >C\right\} \bigcap
\Omega _{\infty ,\alpha }\right)  \notag \\
& \leq \mathbb{P}\left( \inf_{s\in \mathcal{F}_{C}}P_{n}\left(
Ks-Ks_{M}\right) \geq \inf_{s\in \mathcal{F}_{>C}}P_{n}\left(
Ks-Ks_{M}\right) \right)  \notag \\
& =\mathbb{P}\left( \sup_{s\in \mathcal{F}_{C}}P_{n}\left( Ks_{M}-Ks\right)
\leq \sup_{s\in \mathcal{F}_{>C}}P_{n}\left( Ks_{M}-Ks\right) \right)  \notag
\\
& \leq \mathbb{P}\left( \sup_{s\in \mathcal{F}_{\left( \frac{C}{2},\left(
1-\delta \right) C\right] }}P_{n}\left( Ks_{M}-Ks\right) \leq \sup_{s\in 
\mathcal{F}_{>C}}P_{n}\left( Ks_{M}-Ks\right) \right) \text{ }.
\label{deconditionnement_2}
\end{align}%
Now by (\ref{C}) we can apply Lemma \ref{majo_s1>c_lemma} with $\alpha
=\beta $ and we obtain, for all $n\geq n_{0}\left( A_{\infty
},A_{cons},A_{+},\alpha \right) $,%
\begin{equation}
\mathbb{P}\left[ \sup_{s\in \mathcal{F}_{>C}}P_{n}\left( Ks_{M}-Ks\right)
\geq \left( 1+A_{5}\nu _{n}\right) \sqrt{\frac{CD}{n}}\mathcal{K}_{1,M}-C%
\right] \leq 2n^{-\alpha }\text{ }  \label{majo_S>C_end}
\end{equation}%
where $A_{5}$ only depends on $A,A_{3,M},A_{\infty },\sigma _{\min },A_{-}$
and $\alpha $. Moreover, we can take $A_{3,M}=r_{M}\left( \varphi \right) $
by Remark 3. Also, by (\ref{def_1-delta_C}), (\ref{lower_K1M_proof}) and (%
\ref{upper_K1M_proof}) we can apply Lemma \ref{mino_S_r_C} with the quantity 
$C$ in Lemma \ref{mino_S_r_C} replaced by $C/2$, $\alpha =\beta $, $%
r=2\left( 1-\delta \right) $, $A_{u}=\left( \sigma _{\max }+2A\right) ^{2}$, 
$A_{l}=\sigma _{\min }^{2}$ and the constant $A_{\infty }$ satisfying%
\begin{equation}
A_{\infty }\geq 64\sqrt{2}B_{2}A\left( \sigma _{\max }+2A\right) \sigma
_{\min }^{-1}r_{M}\left( \varphi \right) \text{ },  \label{def_A_infini_2}
\end{equation}%
and so it holds, for all $n\geq n_{0}\left( A_{-},A_{+},A,A_{\infty
},A_{cons},B_{2},r_{M}\left( \varphi \right) ,\sigma _{\max },\sigma _{\min
},\alpha \right) $,%
\begin{equation}
\mathbb{P}\left( 
\begin{array}{c}
\sup_{s\in \mathcal{F}_{\left( \frac{C}{2},\left( 1-\delta \right) C\right]
}}P_{n}\left( Ks_{M}-Ks\right) \\ 
\leq \left( 1-L_{A_{-},A,A_{\infty },\sigma _{\max },\sigma _{\min
},r_{M}\left( \varphi \right) ,\alpha }\times \nu _{n}\right) \sqrt{\frac{%
\left( 1-\delta \right) CD}{n}}\mathcal{K}_{1,M}-\left( 1-\delta \right) C%
\end{array}%
\right) \leq 2n^{-\alpha }\text{ }.  \label{mino_S_1-delta_C}
\end{equation}%
Hence from (\ref{majo_S>C_end}) and (\ref{mino_S_1-delta_C}), we deduce that
a positive constant $\check{A}_{0}$ exists, only depending on $%
A_{-},A,A_{\infty },\sigma _{\max },\sigma _{\min },r_{M}\left( \varphi
\right) $ and $\alpha $, such that

\noindent for all $n\geq n_{0}\left( A_{-},A_{+},A,A_{\infty
},A_{cons},B_{2},r_{M}\left( \varphi \right) ,\sigma _{\max },\sigma _{\min
},\alpha \right) $ it holds on an event of probability at least $%
1-4n^{-\alpha }$,%
\begin{equation}
\sup_{s\in \mathcal{F}_{\left( \frac{C}{2},\left( 1-\delta \right) C\right]
}}P_{n}\left( Ks_{M}-Ks\right) \geq \left( 1-\check{A}_{0}\nu _{n}\right) 
\sqrt{\frac{\left( 1-\delta \right) CD}{n}}\mathcal{K}_{1,M}-\left( 1-\delta
\right) C  \label{mino_S_1-delta_C_2}
\end{equation}%
and 
\begin{equation}
\sup_{s\in \mathcal{F}_{>C}}P_{n}\left( Ks_{M}-Ks\right) \leq \left( 1+%
\check{A}_{0}\nu _{n}\right) \sqrt{\frac{CD}{n}}\mathcal{K}_{1,M}-C\text{ .}
\label{majo_S>C_2}
\end{equation}%
Now, from (\ref{mino_S_1-delta_C_2}) and (\ref{majo_S>C_2}) we deduce, using
(\ref{conditionnement_2}) and (\ref{deconditionnement_2}), that if we choose 
$\delta \in \left( 0,\frac{1}{2}\right) $ such that (\ref{C}) and%
\begin{equation}
\left( 1+\check{A}_{0}\nu _{n}\right) \sqrt{\frac{CD}{n}}\mathcal{K}%
_{1,M}-C<\left( 1-\check{A}_{0}\nu _{n}\right) \sqrt{\frac{\left( 1-\delta
\right) CD}{n}}\mathcal{K}_{1,M}-\left( 1-\delta \right) C
\label{inequation_delta}
\end{equation}%
are satisfied then, for all $n\geq n_{0}\left( A_{-},A_{+},A,A_{\infty
},A_{cons},B_{2},r_{M}\left( \varphi \right) ,\sigma _{\max },\sigma _{\min
},n_{1},\alpha \right) $, 
\begin{equation*}
P\left( Ks_{n}-Ks_{M}\right) \leq C\text{ ,}
\end{equation*}%
with probability at least $1-5n^{-\alpha }.$ By (\ref{def_1-delta_C}) it
holds%
\begin{equation*}
\sqrt{\frac{\left( 1-\delta \right) CD}{n}}\mathcal{K}_{1,M}=2\left(
1-\delta \right) C=\frac{1}{2}\frac{D}{n}\mathcal{K}_{1,M}^{2}\text{ ,}
\end{equation*}%
and by consequence, inequality (\ref{inequation_delta}) is equivalent to 
\begin{equation}
\left( 1-2\check{A}_{0}\nu _{n}\right) \left( 1-\delta \right) -2\left( 1+%
\check{A}_{0}\nu _{n}\right) \sqrt{1-\delta }+1>0\text{ .}
\label{inequation_delta_2}
\end{equation}%
Moreover, we have by (\ref{nu_n_proof}) and (\textbf{H5}), for all $n\geq
n_{0}\left( A_{+},A_{-},A_{cons},\check{A}_{0},A_{5},\alpha \right) $, 
\begin{equation}
\left( \check{A}_{0}\vee A_{5}\right) \nu _{n}<\frac{1}{72}
\label{majo_nu_n_end_2}
\end{equation}%
and so, for all $n\geq n_{0}\left( A_{+},A_{-},A_{cons},\check{A}_{0},\alpha
\right) $, simple computations involving (\ref{majo_nu_n_end_2}) show that
by taking 
\begin{equation}
\delta =6\left( \sqrt{\check{A}_{0}}\vee \sqrt{A_{5}}\right) \sqrt{\nu _{n}}%
\text{ ,}  \label{def_delta}
\end{equation}%
inequalities (\ref{inequation_delta_2}) and (\ref{C}) are satisfied and $%
\delta \in \left( 0,\frac{1}{2}\right) $. Finally, we can compute $C$ by (%
\ref{def_1-delta_C}) and (\ref{def_delta}), in such a way that for all $%
n\geq n_{0}\left( A_{+},A_{-},A_{cons},\check{A}_{0},\alpha \right) $ 
\begin{equation}
0<C=\frac{\left( 1-\delta \right) C}{\left( 1-\delta \right) }=\frac{1}{%
\left( 1-\delta \right) }\frac{1}{4}\frac{D}{n}\mathcal{K}_{1,M}^{2}\leq
\left( 1+12\left( \sqrt{\check{A}_{0}}\vee \sqrt{A_{5}}\right) \sqrt{\nu _{n}%
}\right) \frac{1}{4}\frac{D}{n}\mathcal{K}_{1,M}^{2}\text{ },
\label{def_A_0_2}
\end{equation}%
which yields the result by noticing that the dependence on $\sigma _{\max }$
can be released from $n_{0}$ and $\check{A}_{0}$ since by (\textbf{H1})\ we
have $\sigma _{\max }\leq A$.\ 

\vspace*{0.5in}

\subsection*{}%
\textbf{Proof of Inequality (\ref{lower_emp_risk}). }Let $C=\frac{D}{8n}%
\mathcal{K}_{1,M}^{2}>0$ and let $r=2$. By (\ref{upper_K1M_proof}) and (\ref%
{lower_K1M_proof}) we have%
\begin{equation*}
\frac{D}{n}\sigma _{\min }^{2}\leq rC=\frac{D}{4n}\mathcal{K}_{1,M}^{2}\leq 
\frac{D}{n}\left( \sigma _{\max }+2A\right) ^{2}
\end{equation*}%
so we can apply Lemma \ref{mino_S_r_C} with $\alpha =\beta $, $A_{l}=\sigma
_{\min }^{2}$ and $A_{u}=\left( \sigma _{\max }+2A\right) ^{2}$. So if 
\begin{equation}
A_{\infty }\geq 64\sqrt{2}B_{2}A\left( \sigma _{\max }+2A\right) \sigma
_{\min }^{-1}r_{M}\left( \varphi \right) \text{ },  \label{def_A_infini_3}
\end{equation}%
it holds, for all $n\geq n_{0}\left( A_{-},A_{+},A,A_{\infty
},A_{cons},B_{2},r_{M}\left( \varphi \right) ,\sigma _{\max },\sigma _{\min
},\alpha \right) $,%
\begin{equation}
\mathbb{P}\left( \sup_{s\in \mathcal{F}_{\left( C,rC\right] }}P_{n}\left(
Ks_{M}-Ks\right) \leq \left( 1-L_{A_{-},A,A_{\infty },\sigma _{\max },\sigma
_{\min },r_{M}\left( \varphi \right) ,\alpha }\times \nu _{n}\right) \sqrt{%
\frac{rCD}{n}}\mathcal{K}_{1,M}-rC\right) \leq 2n^{-\alpha }\text{ }.
\label{def_L}
\end{equation}%
Since $rC=\frac{D}{4n}\mathcal{K}_{1,M}^{2}$, if we set $\hat{A}%
_{0}=2L_{A_{-},A,A_{\infty },\sigma _{\max },\sigma _{\min },r_{M}\left(
\varphi \right) ,\alpha }$ with $L_{A_{-},A,A_{\infty },\sigma _{\max
},\sigma _{\min },r_{M}\left( \varphi \right) ,\alpha }$ the constant in (%
\ref{def_L}), we get%
\begin{equation}
\mathbb{P}\left( \sup_{s\in \mathcal{F}_{\left( C,rC\right] }}P_{n}\left(
Ks_{M}-Ks\right) \leq \left( 1-\hat{A}_{0}\nu _{n}\right) \frac{D}{4n}%
\mathcal{K}_{1,M}^{2}\right) \leq 2n^{-\alpha }\text{ }.  \label{dev_3}
\end{equation}%
Notice that 
\begin{equation*}
P_{n}\left( Ks_{M}-Ks_{n}\right) =\sup_{s\in M}P_{n}\left( Ks_{M}-Ks\right)
\geq \sup_{s\in \mathcal{F}_{\left( C,rC\right] }}P_{n}\left(
Ks_{M}-Ks\right)
\end{equation*}%
so from (\ref{dev_3}) we deduce that%
\begin{equation}
\mathbb{P}\left( P_{n}\left( Ks_{M}-Ks_{n}\right) \geq \left( 1-\hat{A}%
_{0}\nu _{n}\right) \frac{D}{4n}\mathcal{K}_{1,M}^{2}\right) \geq
1-2n^{-\alpha }\text{ }.  \label{def_A_0_3}
\end{equation}

\bigskip

\noindent \textbf{Remark 4} \textit{Notice that in the proof of inequality (%
\ref{lower_emp_risk}), we do not need to assume the consistency of the least
squares estimator }$s_{n}$\textit{\ towards the projection }$s_{M}$\textit{.
Straightforward adaptations of Lemma \ref{mino_S_r_C} allow to take }%
\begin{equation*}
\tilde{\nu}_{n}=\max \left\{ \sqrt{\frac{\ln n}{D}},\text{ }\sqrt{\frac{D\ln
n}{n}}\right\}
\end{equation*}%
\textit{instead of the quantity }$\nu _{n}$\textit{\ defined in (\ref%
{nu_n_proof}). This readily gives the expected bound (\ref%
{lower_emp_risk_general}) of Theorem \ref{principal_reg}.}\vspace*{0.5in}\ 

\subsection*{}%
\textbf{Proof of Inequality (\ref{upper_emp_risk}). }Let 
\begin{equation}
C=\frac{1}{4}\left( 1+A_{5}\nu _{n}\right) ^{2}\frac{D}{n}\mathcal{K}%
_{1,M}^{2}>0  \label{def_C_4}
\end{equation}%
where $A_{5}$\textbf{\ }is defined in Lemma \ref{majo_s1>c_lemma} applied
with $\beta =\alpha $. By (\textbf{H5}) we have 
\begin{equation}
\mathbb{P}\left( P_{n}\left( Ks_{M}-Ks_{n}\right) >C\right) \leq \mathbb{P}%
\left( \left\{ P_{n}\left( Ks_{M}-Ks_{n}\right) >C\right\} \bigcap \Omega
_{\infty ,\alpha }\right) +n^{-\alpha }\text{ .}  \label{deconditionnement_3}
\end{equation}%
Moreover, on $\Omega _{\infty ,\alpha }$, we have%
\begin{align}
P_{n}\left( Ks_{M}-Ks_{n}\right) & =\sup_{s\in B_{\left( M,L_{\infty
}\right) }\left( s_{M},\tilde{R}_{n,D,\alpha }\right) }P_{n}\left(
Ks_{M}-Ks\right)  \notag \\
& =\sup_{s\in \mathcal{F}_{>0}}P_{n}\left( Ks_{M}-Ks\right)  \label{sup_p2}
\end{align}%
and by (\ref{majo_S>C_4}) of Lemma \ref{majo_s1>c_lemma} applied with $%
\alpha =\beta $ it holds, for all $n\geq n_{0}\left( A_{\infty
},A_{cons},A_{+},\alpha \right) $,%
\begin{equation}
\mathbb{P}\left( \sup_{s\in \mathcal{F}_{>0}}P_{n}\left( Ks_{M}-Ks\right)
>C\right) \leq 2n^{-\alpha }\text{ .}  \label{sup_F>0}
\end{equation}%
Finally, using (\ref{sup_p2}) and (\ref{sup_F>0}) in (\ref%
{deconditionnement_3}) we get, for all $n\geq n_{0}\left( A_{\infty
},A_{cons},n_{1},A_{+},\alpha \right) $,%
\begin{equation*}
\mathbb{P}\left( P_{n}\left( Ks_{M}-Ks_{n}\right) >C\right) \leq 3n^{-\alpha
}\text{ .}
\end{equation*}

\vspace*{0.5in}

\noindent\textbf{Conclusion.} To complete the proof of Theorem \ref%
{principal_reg}, just notice that by (\ref{def_A_infini_1}), (\ref%
{def_A_infini_2}) and (\ref{def_A_infini_3}) we can take%
\begin{equation*}
A_{\infty}=64\sqrt{2}B_{2}A\left( \sigma_{\max}+2A\right)
\sigma_{\min}^{-1}r_{M}\left( \varphi\right)
\end{equation*}
and by (\ref{def_A_0_1}), (\ref{def_A_0_2}), (\ref{def_A_0_3}) and (\ref%
{def_C_4}),%
\begin{equation*}
A_{0}=\max\left\{ 48\sqrt{\tilde{A}_{0}},\text{ }12\left( \sqrt{\check {A}%
_{0}}\vee\sqrt{A_{5}}\right) ,\text{ }\sqrt{\hat{A}_{0}},\text{ }\sqrt{A_{5}}%
\right\}
\end{equation*}
is convenient. 
$\blacksquare$%

\subsection*{}%
\textbf{Proof of Theorem \ref{theorem_upper_true_small_models}. }We localize
our analysis in the subset 
\begin{equation*}
B_{\left( M,L_{\infty }\right) }\left( s_{M},R_{n,D,\alpha }\right) =\left\{
s\in M,\left\Vert s-s_{M}\right\Vert _{\infty }\leq R_{n,D,\alpha }\right\}
\subset M\text{ .}
\end{equation*}%
Unlike in the proof of Theorem \ref{principal_reg}, see (\ref{def_R_n_D_tild}%
), we need not to consider the quantity $\tilde{R}_{n,D,\alpha }$, a radius
possibly larger than $R_{n,D,\alpha }$. Indeed, the use of $\tilde{R}%
_{n,D,\alpha }$ rather than $R_{n,D,\alpha }$ in the proof of Theorem \ref%
{principal_reg} is only needed in Lemma \ref{bound_moment_order_1}, where we
derive a sharp lower bound for the mean of the supremum of the empirical
process indexed by the contrasted functions centered by the contrasted
projection over a slice of interest. To prove Theorem \ref%
{theorem_upper_true_small_models}, we just need upper bounds, and Lemma \ref%
{bound_moment_order_1} is avoided as well as the use of $\tilde{R}%
_{n,D,\alpha }$.

\noindent Let us define several slices of excess risk on the model $M$ : for
any $C\geq 0$,%
\begin{align*}
\mathcal{G}_{C}& =\left\{ s\in M,P\left( Ks-Ks_{M}\right) \leq C\right\}
\bigcap B_{\left( M,L_{\infty }\right) }\left( s_{M},R_{n,D,\alpha }\right) 
\text{ ,} \\
\mathcal{G}_{>C}& =\left\{ s\in M,P\left( Ks-Ks_{M}\right) >C\right\}
\bigcap B_{\left( M,L_{\infty }\right) }\left( s_{M},R_{n,D,\alpha }\right) 
\text{ .}
\end{align*}%
We also define, for all $U\geq 0$,%
\begin{equation*}
\mathcal{D}_{U}=\left\{ s\in M,P\left( Ks-Ks_{M}\right) =U\right\} \bigcap
B_{\left( M,L_{\infty }\right) }\left( s_{M},R_{n,D,\alpha }\right) \text{ }.
\end{equation*}

\noindent \textbf{I. Proof of Inequality (\ref{upper_true_small}).} Let $%
C_{1}>0$ to be fixed later, satisfying 
\begin{equation}
C_{1}\geq \frac{D}{n}=:C_{-}>0\text{ .}  \label{uts_0}
\end{equation}%
We have by (\textbf{H5}), for all $n\geq n_{1},$%
\begin{equation}
\mathbb{P}\left( P\left( Ks_{n}-Ks_{M}\right) >C_{1}\right) \leq \mathbb{P}%
\left( \left\{ P\left( Ks_{n}-Ks_{M}\right) >C_{1}\right\} \bigcap \Omega
_{\infty ,\alpha }\right) +n^{-\alpha }  \label{uts_1}
\end{equation}%
and also%
\begin{align}
& \mathbb{P}\left( \left\{ P\left( Ks_{n}-Ks_{M}\right) >C_{1}\right\}
\bigcap \Omega _{\infty ,\alpha }\right)  \notag \\
& \leq \mathbb{P}\left( \inf_{s\in \mathcal{G}_{C_{1}}}P_{n}\left(
Ks-Ks_{M}\right) \geq \inf_{s\in \mathcal{G}_{>C_{1}}}P_{n}\left(
Ks-Ks_{M}\right) \right)  \notag \\
& =\mathbb{P}\left( \sup_{s\in \mathcal{G}_{C_{1}}}P_{n}\left(
Ks_{M}-Ks\right) \leq \sup_{s\in \mathcal{G}_{>C_{1}}}P_{n}\left(
Ks_{M}-Ks\right) \right)  \notag \\
& \leq \mathbb{P}\left( 0\leq \sup_{s\in \mathcal{G}_{>C_{1}}}P_{n}\left(
Ks_{M}-Ks\right) \right) \text{ }.  \label{uts_2}
\end{align}%
Moreover, it holds 
\begin{align}
& \sup_{s\in \mathcal{G}_{>C_{1}}}P_{n}\left( Ks_{M}-Ks\right)  \notag \\
& =\sup_{s\in \mathcal{G}_{>C_{1}}}\left\{ P_{n}\left( \psi _{1,M}%
\cdot%
\left( s_{M}-s\right) -\psi _{2}\circ \left( s-s_{M}\right) \right) \right\}
\notag \\
& =\sup_{s\in \mathcal{G}_{>C_{1}}}\left\{ \left( P_{n}-P\right) \left( \psi
_{1,M}%
\cdot%
\left( s_{M}-s\right) \right) -\left( P_{n}-P\right) \left( \psi _{2}\circ
\left( s-s_{M}\right) \right) -P\left( Ks-Ks_{M}\right) \right\}  \notag \\
& =\sup_{s\in \mathcal{G}_{>C_{1}}}\left\{ \left( P_{n}-P\right) \left( \psi
_{1,M}%
\cdot%
\left( s_{M}-s\right) \right) -P\left( Ks-Ks_{M}\right) -\left(
P_{n}-P\right) \left( \psi _{2}\circ \left( s-s_{M}\right) \right) \right\} 
\notag \\
& =\sup_{U>C_{1}}\sup_{s\in \mathcal{D}_{U}}\left\{ \left( P_{n}-P\right)
\left( \psi _{1,M}%
\cdot%
\left( s_{M}-s\right) \right) -U-\left( P_{n}-P\right) \left( \psi _{2}\circ
\left( s-s_{M}\right) \right) \right\}  \notag \\
& \leq \sup_{U>C_{1}}\left\{ \sqrt{U}\sqrt{\sum_{k=1}^{D}\left(
P_{n}-P\right) ^{2}\left( \psi _{1,M}%
\cdot%
\varphi _{k}\right) }-U+\sup_{s\in \mathcal{G}_{U}}\left\vert \left(
P_{n}-P\right) \left( \psi _{2}\circ \left( s-s_{M}\right) \right)
\right\vert \right\} \text{ .}  \label{uts_3}
\end{align}%
Now, from inequality (\ref{chi_reg_2}) of Lemma \ref{mino_sphere_lemma}
applied with $\beta =\alpha $, we get%
\begin{equation}
\mathbb{P}\left[ \sqrt{\sum_{k=1}^{D}\left( P_{n}-P\right) ^{2}\left( \psi
_{1,M}%
\cdot%
\varphi _{k}\right) }\geq L_{A,A_{3,M},\alpha }\sqrt{\frac{D\vee \ln n}{n}}%
\right] \leq n^{-\alpha }\text{ .}  \label{uts_4}
\end{equation}%
In addition, we handle the empirical process indexed by the second order
terms by straightforward modifications of Lemmas \ref{lower_rademacher
copy(1)} and \ref{majo_uniform_s2c}\ as well as their proofs. It thus holds,
by the same type of arguments as those given in Lemma \ref{lower_rademacher
copy(1)},%
\begin{equation}
\mathbb{E}\left[ \sup_{s\in \mathcal{G}_{C_{1}}}\left\vert \left(
P_{n}-P\right) \left( \psi _{2,M}^{s}%
\cdot%
\left( s-s_{M}\right) \right) \right\vert \right] \leq 8\sqrt{\frac{CD}{n}}%
R_{n,D,\alpha }\text{ }.  \label{upper_mean_2_G_reg}
\end{equation}%
Moreover, using (\ref{upper_mean_2_G_reg}), the same type of arguments as
those leading to inequality (\ref{dev_2}) of Lemma \ref{majo_uniform_s2c},
allow to show that for any $q\geq 1$ and $j\in \mathbb{N}^{\ast }$, for all $%
x>0,$%
\begin{eqnarray}
&&\mathbb{P}\left[ \sup_{s\in \mathcal{G}_{q^{j}C_{-}}}\left\vert \left(
P_{n}-P\right) \left( \psi _{2}\circ \left( s-s_{M}\right) \right)
\right\vert \geq 16\sqrt{\frac{q^{j}C_{-}D}{n}}R_{n,D,\alpha }+\sqrt{\frac{%
2R_{n,D,\alpha }^{2}q^{j}C_{-}x}{n}}+\frac{8}{3}\frac{R_{n,D,\alpha }^{2}x}{n%
}\right]  \notag \\
&\leq &\exp \left( -x\right) \text{ }.  \label{uts_5}
\end{eqnarray}%
Hence, taking $x=\gamma \ln n$ in (\ref{uts_5}) and using the fact that $%
C_{-}=Dn^{-1}\geq n^{-1}$, we get%
\begin{equation}
\mathbb{P}\left[ \sup_{s\in \mathcal{G}_{q^{j}C_{-}}}\left\vert \left(
P_{n}-P\right) \left( \psi _{2}\circ \left( s-s_{M}\right) \right)
\right\vert \geq L_{A_{cons},\gamma }R_{n,D,\alpha }\sqrt{\frac{%
q^{j}C_{-}\left( D\vee \ln n\right) }{n}}\right] \leq n^{-\gamma }\text{ }.
\label{uts_6}
\end{equation}%
Now, by straightforward modifications of the proof of Lemma \ref%
{majo_uniform_s2c}, we get that for all $n\geq n_{0}\left( A_{cons}\right) $,%
\begin{equation}
\mathbb{P}\left[ \forall U>C_{-},\text{ \ }\sup_{s\in \mathcal{G}%
_{U}}\left\vert \left( P_{n}-P\right) \left( \psi _{2}\circ \left(
s-s_{M}\right) \right) \right\vert \leq L_{A_{cons},\alpha }R_{n,D,\alpha }%
\sqrt{\frac{U\left( D\vee \ln n\right) }{n}}\right] \geq 1-n^{-\alpha }\text{
}.  \label{uts_7}
\end{equation}%
Combining (\ref{uts_3}), (\ref{uts_4}) and (\ref{uts_7}), we have on an
event of probability at least $1-2n^{-\alpha }$, for all $n\geq n_{0}\left(
A_{cons}\right) $,%
\begin{eqnarray}
\sup_{s\in \mathcal{G}_{>C_{1}}}P_{n}\left( Ks_{M}-Ks\right) &\leq
&\sup_{U>C_{1}}\left\{ L_{A,A_{3,M},\alpha }\sqrt{\frac{U\left( D\vee \ln
n\right) }{n}}-U+L_{A_{cons},\alpha }R_{n,D,\alpha }\sqrt{\frac{U\left(
D\vee \ln n\right) }{n}}\right\}  \notag \\
&\leq &\sup_{U>C_{1}}\left\{ L_{A,A_{cons},A_{3,M},\alpha }\left(
1+R_{n,D,\alpha }\right) \sqrt{\frac{U\left( D\vee \ln n\right) }{n}}%
-U\right\} \text{ .}  \label{uts_8}
\end{eqnarray}%
Now, as $R_{n,D,\alpha }\leq A_{cons}\left( \ln n\right) ^{-1/2}$, we deduce
from (\ref{uts_8}) that for 
\begin{equation}
C_{1}=L_{A,A_{cons},A_{3,M},\alpha }\frac{D\vee \ln \left( n\right) }{n}%
>C_{-}  \label{def_C1}
\end{equation}%
with $L_{A,A_{cons},A_{3,M},\alpha }$\ large enough, it holds with
probability at least $1-2n^{-\alpha }$ and for all $n\geq n_{0}\left(
A_{cons}\right) $,%
\begin{equation*}
\sup_{s\in \mathcal{G}_{>C_{1}}}P_{n}\left( Ks_{M}-Ks\right) <0\text{ ,}
\end{equation*}%
and so by using (\ref{uts_1}) and (\ref{uts_2}), this yields inequality (\ref%
{upper_true_small}).

\bigskip

\noindent \textbf{II. Proof of Inequality (\ref{upper_emp_small}).} Let $%
C_{2}>0$ to be fixed later, satisfying 
\begin{equation}
C_{2}\geq \frac{D}{n}=C_{-}>0\text{ .}  \label{uts_9}
\end{equation}%
We have by (\textbf{H5}), for all $n\geq n_{1},$%
\begin{equation}
\mathbb{P}\left( P_{n}\left( Ks_{M}-Ks_{n}\right) >C_{2}\right) \leq \mathbb{%
P}\left( \left\{ P_{n}\left( Ks_{M}-Ks_{n}\right) >C_{2}\right\} \bigcap
\Omega _{\infty ,\alpha }\right) +n^{-\alpha }\text{ .}  \label{uts_10}
\end{equation}%
Moreover, we have on $\Omega _{\infty ,\alpha }$,%
\begin{align}
P_{n}\left( Ks_{M}-Ks_{n}\right) & =\sup_{s\in B_{\left( M,L_{\infty
}\right) }\left( s_{M},R_{n,D,\alpha }\right) }P_{n}\left( Ks_{M}-Ks\right) 
\notag \\
& =\max \left\{ \sup_{s\in \mathcal{G}_{C_{1}}}P_{n}\left( Ks_{M}-Ks\right) 
\text{ };\text{ }\sup_{s\in \mathcal{G}_{>C_{1}}}P_{n}\left(
Ks_{M}-Ks\right) \right\} \text{ ,}  \label{uts_11}
\end{align}%
where $C_{1}$ is defined in the first part of the proof dedicated to the
establishment of inequality (\ref{upper_true_small}). Moreover, let us
recall that in the first part of the proof, we have proved that an event of
probability at least $1-2n^{-\alpha }$ exists, that we call $\Omega _{1}$,
such that it holds on this event, for all $n\geq n_{0}\left( A_{cons}\right) 
$,%
\begin{equation}
\sqrt{\sum_{k=1}^{D}\left( P_{n}-P\right) ^{2}\left( \psi _{1,M}%
\cdot%
\varphi _{k}\right) }\leq L_{A,A_{3,M},\alpha }\sqrt{\frac{D\vee \ln n}{n}}%
\text{ ,}  \label{uts_12}
\end{equation}%
\begin{equation}
\forall U>C_{-},\text{ \ }\sup_{s\in \mathcal{G}_{U}}\left\vert \left(
P_{n}-P\right) \left( \psi _{2}\circ \left( s-s_{M}\right) \right)
\right\vert \leq L_{A_{cons},\alpha }R_{n,D,\alpha }\sqrt{\frac{U\left(
D\vee \ln n\right) }{n}}\text{ ,}  \label{uts_13}
\end{equation}%
and%
\begin{equation}
\sup_{s\in \mathcal{G}_{>C_{1}}}P_{n}\left( Ks_{M}-Ks\right) <0\text{ .}
\label{uts_14}
\end{equation}%
By (\ref{uts_11}) and (\ref{uts_14}), we thus have on $\Omega _{\infty
,\alpha }\bigcap \Omega _{1}$, for all $n\geq n_{0}\left( A_{cons}\right) $, 
\begin{equation}
0\leq P_{n}\left( Ks_{M}-Ks_{n}\right) =\sup_{s\in \mathcal{G}%
_{C_{1}}}P_{n}\left( Ks_{M}-Ks\right) \text{ .}  \label{uts_15}
\end{equation}%
In addition, it holds%
\begin{align}
& \sup_{s\in \mathcal{G}_{C_{1}}}P_{n}\left( Ks_{M}-Ks\right)  \notag \\
& =\sup_{s\in \mathcal{G}_{C_{1}}}\left\{ P_{n}\left( \psi _{1,M}%
\cdot%
\left( s_{M}-s\right) -\psi _{2}\circ \left( s-s_{M}\right) \right) \right\}
\notag \\
& =\sup_{s\in \mathcal{G}_{C_{1}}}\left\{ \left( P_{n}-P\right) \left( \psi
_{1,M}%
\cdot%
\left( s_{M}-s\right) \right) -\left( P_{n}-P\right) \left( \psi _{2}\circ
\left( s-s_{M}\right) \right) -P\left( Ks-Ks_{M}\right) \right\}  \notag \\
& \leq \sup_{s\in \mathcal{G}_{C_{1}}}\left\{ \left( P_{n}-P\right) \left(
\psi _{1,M}%
\cdot%
\left( s_{M}-s\right) \right) \right\} +\sup_{s\in \mathcal{G}%
_{C_{1}}}\left\vert \left( P_{n}-P\right) \left( \psi _{2}\circ \left(
s-s_{M}\right) \right) \right\vert \text{ .}  \label{uts_16}
\end{align}%
Now, we have on $\Omega _{1}$, for all $n\geq n_{0}\left( A_{cons}\right) $,%
\begin{eqnarray}
\sup_{s\in \mathcal{G}_{C_{1}}}\left\{ \left( P_{n}-P\right) \left( \psi
_{1,M}%
\cdot%
\left( s_{M}-s\right) \right) \right\} &\leq &\sqrt{C_{1}}\sqrt{%
\sum_{k=1}^{D}\left( P_{n}-P\right) ^{2}\left( \psi _{1,M}%
\cdot%
\varphi _{k}\right) }  \notag \\
&\leq &L_{A,A_{3,M},\alpha }\sqrt{\frac{C_{1}\left( D\vee \ln n\right) }{n}}%
\text{ \ \ \ \ \ by (\ref{uts_12})}  \notag \\
&=&L_{A,A_{cons},A_{3,M},\alpha }\frac{D\vee \ln \left( n\right) }{n}\text{
\ \ \ \ by (\ref{def_C1}) ,}  \label{uts_17}
\end{eqnarray}%
and also, by (\ref{uts_13}) and (\ref{def_C1}), 
\begin{eqnarray}
\sup_{s\in \mathcal{G}_{C_{1}}}\left\vert \left( P_{n}-P\right) \left( \psi
_{2}\circ \left( s-s_{M}\right) \right) \right\vert &\leq
&L_{A_{cons},\alpha }R_{n,D,\alpha }\sqrt{\frac{C_{1}\left( D\vee \ln
n\right) }{n}}  \notag \\
&\leq &L_{A,A_{cons},A_{3,M},\alpha }R_{n,D,\alpha }\frac{D\vee \ln \left(
n\right) }{n}\text{ .}  \label{uts_18}
\end{eqnarray}%
Finally, as $R_{n,D,\alpha }\leq A_{cons}\left( \ln n\right) ^{-1/2}$, we
deduce from (\ref{uts_15}), (\ref{uts_16}), (\ref{uts_17}) and (\ref{uts_18}%
), that it holds on $\Omega _{\infty ,\alpha }\bigcap \Omega _{1}$, for all $%
n\geq n_{0}\left( A_{cons}\right) $, 
\begin{equation*}
P_{n}\left( Ks_{M}-Ks_{n}\right) \leq L_{A,A_{cons},A_{3,M},\alpha }\frac{%
D\vee \ln \left( n\right) }{n}\text{ ,}
\end{equation*}%
and so, this yields to inequality (\ref{upper_emp_small}) by using (\ref%
{uts_10}) and this concludes the proof of Theorem \textbf{\ref%
{theorem_upper_true_small_models}}.\ 
$\blacksquare$%

\subsection{Technical Lemmas\label{useful_lemmas}}

We state here some lemmas needed in the proofs of Section \ref%
{section_proof_fixed_model}. First, in Lemmas \ref{mino_sphere_lemma}, \ref%
{bound_moment_order_1} and \ref{control_lower_infinity}, we derive some
controls, from above and from below, of the empirical process indexed by the
\textquotedblleft linear parts\textquotedblright\ of the contrasted
functions over slices of interest. Secondly, we give upper bounds in Lemmas %
\ref{lower_rademacher copy(1)} and \ref{majo_uniform_s2c} for the empirical
process indexed by the \textquotedblleft quadratic parts\textquotedblright\
of the contrasted functions over slices of interest. And finally, we use all
these results in Lemmas \ref{majo_s1c_lemma copy(1)}, \ref{majo_s1>c_lemma}
and \ref{mino_S_r_C} to derive upper and lower bounds for the empirical
process indexed by the contrasted functions over slices of interest.

\begin{lemma}
\label{mino_sphere_lemma}Assume that (\textbf{H1}), (\textbf{H2}) and (%
\textbf{H3}) hold. Then for any $\beta >0,$ by setting 
\begin{equation*}
\tau _{n}=L_{A,A_{3,M},\sigma _{\min },\beta }\left( \sqrt{\frac{\ln n}{D}}%
\vee \frac{\sqrt{\ln n}}{n^{1/4}}\right) \text{ },
\end{equation*}%
It holds, for any orthonormal basis $\left( \varphi _{k}\right) _{k=1}^{D}$
of $\left( M,\left\Vert 
\cdot%
\right\Vert _{2}\right) $, 
\begin{equation}
\mathbb{P}\left[ \sqrt{\sum_{k=1}^{D}\left( P_{n}-P\right) ^{2}\left( \psi
_{1,M}%
\cdot%
\varphi _{k}\right) }\geq \left( 1+\tau _{n}\right) \sqrt{\frac{D}{n}}%
\mathcal{K}_{1,M}\right] \leq n^{-\beta }\text{ }.  \label{chi_reg_1}
\end{equation}%
If (\textbf{H1}) and (\textbf{H3}) hold, then for any $\beta >0,$ it holds%
\begin{equation}
\mathbb{P}\left[ \sqrt{\sum_{k=1}^{D}\left( P_{n}-P\right) ^{2}\left( \psi
_{1,M}%
\cdot%
\varphi _{k}\right) }\geq L_{A,A_{3,M},\beta }\sqrt{\frac{D\vee \ln n}{n}}%
\right] \leq n^{-\beta }\text{ }.  \label{chi_reg_2}
\end{equation}
\end{lemma}

\subsection*{}%
\textbf{Proof. }By Cauchy-Schwarz inequality we have%
\begin{equation*}
\chi _{M}:=\sqrt{\sum_{k=1}^{D}\left( P_{n}-P\right) ^{2}\left( \psi _{1,M}%
\cdot%
\varphi _{k}\right) }=\sup_{s\in M\text{, }\left\Vert s\right\Vert _{2}\leq
1}\left\{ \left\vert \left( P_{n}-P\right) \left( \psi _{1,M}%
\cdot%
s\right) \right\vert \right\} \text{ }.
\end{equation*}%
Hence, we get by Bousquet's inequality (\ref{bousquet_2}) applied with $%
\mathcal{F=}\left\{ \psi _{1,M}%
\cdot%
s\text{ };\text{ }s\in M,\text{ }\left\Vert s\right\Vert _{2}\leq 1\right\} $%
, for all $x>0$, $\delta >0,$%
\begin{equation}
\mathbb{P}\left[ \chi _{M}\geq \sqrt{2\sigma ^{2}\frac{x}{n}}+\left(
1+\delta \right) \mathbb{E}\left[ \chi _{M}\right] +\left( \frac{1}{3}+\frac{%
1}{\delta }\right) \frac{bx}{n}\right] \leq \exp \left( -x\right)
\label{mino_sup_sphere_1}
\end{equation}%
where 
\begin{equation*}
\sigma ^{2}\leq \sup_{s\in M,\text{ }\left\Vert s\right\Vert _{2}\leq 1}P%
\left[ \left( \psi _{1,M}%
\cdot%
s\right) ^{2}\right] \leq \left\Vert \psi _{1,M}\right\Vert _{\infty
}^{2}\leq 16A^{2}\text{ \ \ \ \ \ \ by (\ref{control_phi_1M_proof})}
\end{equation*}%
and 
\begin{equation*}
b\leq \sup_{s\in M,\text{ }\left\Vert s\right\Vert _{2}\leq 1}\left\Vert
\psi _{1,M}%
\cdot%
s-P\left( \psi _{1,M}%
\cdot%
s\right) \right\Vert _{\infty }\leq 4A\sqrt{D}A_{3,M}\text{ \ \ \ \ \ \ by (%
\ref{center_proof}), (\ref{control_phi_1M_proof})\ and (\ref%
{control_boule_unite}).}
\end{equation*}%
Moreover,%
\begin{equation*}
\mathbb{E}\left[ \chi _{M}\right] \leq \sqrt{\mathbb{E}\left[ \chi _{M}^{2}%
\right] }=\sqrt{\frac{D}{n}}\mathcal{K}_{1,M}\text{ }.
\end{equation*}%
So, from (\ref{mino_sup_sphere_1}) it follows that, for all $x>0$, $\delta
>0,$%
\begin{equation}
\mathbb{P}\left[ \chi _{M}\geq \sqrt{32A^{2}\frac{x}{n}}+\left( 1+\delta
\right) \sqrt{\frac{D}{n}}\mathcal{K}_{1,M}+\left( \frac{1}{3}+\frac{1}{%
\delta }\right) \frac{4A\sqrt{D}A_{3,M}x}{n}\right] \leq \exp \left(
-x\right) \text{ }.  \label{dev_Z_1}
\end{equation}%
Hence, taking $x=\beta \ln n$, $\delta =\frac{\sqrt{\ln n}}{n^{1/4}}$ in (%
\ref{dev_Z_1}), we derive by (\ref{lower_K1M_proof}) that a positive
constant $L_{A,A_{3,M},\sigma _{\min },\beta }$ exists such that 
\begin{equation*}
\mathbb{P}\left[ \chi _{M}\geq \left( 1+L_{A,A_{3,M},\sigma _{\min },\beta
}\left( \sqrt{\frac{\ln n}{D}}\vee \frac{\sqrt{\ln n}}{n^{1/4}}\right)
\right) \sqrt{\frac{D}{n}}\mathcal{K}_{1,M}\right] \leq n^{-\beta }\text{ },
\end{equation*}%
which yields inequality (\ref{chi_reg_1}). By (\ref{upper_K1M_proof}) we
have $\mathcal{K}_{1,M}\leq 6A$, and by taking again $x=\beta \ln n$ and $%
\delta =\frac{\sqrt{\ln n}}{n^{1/4}}$ in (\ref{dev_Z_1}), simple
computations give 
\begin{equation*}
\mathbb{P}\left[ \sqrt{\sum_{k=1}^{D}\left( P_{n}-P\right) ^{2}\left( \psi
_{1,M}%
\cdot%
\varphi _{k}\right) }\geq L_{A,A_{3,M},\beta }\left( \sqrt{\frac{D}{n}}\vee 
\sqrt{\frac{\ln n}{n}}\vee \sqrt{\frac{D\ln n}{n^{3/2}}}\right) \right] \leq
n^{-\beta }\text{ },
\end{equation*}%
and by consequence, (\ref{chi_reg_2}) follows. 
$\blacksquare$%

\noindent In the next lemma, we state sharp lower bounds for the mean of the
supremum of the empirical process on the linear parts of contrasted
functions of $M$ belonging to a slice of excess risk. This is done for a
model of reasonable dimension.

\begin{lemma}
\label{bound_moment_order_1}Let $r>1$ and $C>0$. Assume that (\textbf{H1}), (%
\textbf{H2}), (\textbf{H4}) and (\ref{def_A_cons}) hold and let $\varphi
=\left( \varphi _{k}\right) _{k=1}^{D}$ be an orthonormal basis of $\left(
M,\left\Vert 
\cdot%
\right\Vert _{2}\right) $ satisfying (\textbf{H4}). If positive constants $%
A_{-},A_{+},A_{l},A_{u}$ exist such that 
\begin{equation*}
A_{+}\frac{n}{\left( \ln n\right) ^{2}}\geq D\geq A_{-}\left( \ln n\right)
^{2}\text{ \ \ and \ \ }A_{l}\frac{D}{n}\leq rC\leq A_{u}\frac{D}{n}\text{ },
\end{equation*}%
and if the constant $A_{\infty }$ defined in (\ref{def_R_n_D_tild})
satisfies 
\begin{equation}
A_{\infty }\geq 64B_{2}A\sqrt{2A_{u}}\sigma _{\min }^{-1}r_{M}\left( \varphi
\right) \text{ },  \label{lower_A_infini}
\end{equation}%
then a positive constant $L_{A,A_{l},A_{u},\sigma _{\min }}$ exists such
that, for all $n\geq n_{0}\left( A_{-},A_{+},A_{u},A_{l},A,B_{2},r_{M}\left(
\varphi \right) ,\sigma _{\min }\right) $,%
\begin{equation}
\mathbb{E}\left[ \sup_{s\in \mathcal{F}_{\left( C,rC\right] }}\left(
P_{n}-P\right) \left( \psi _{1,M}%
\cdot%
\left( s_{M}-s\right) \right) \right] \geq \left( 1-\frac{%
L_{A,A_{l},A_{u},\sigma _{\min }}}{\sqrt{D}}\right) \sqrt{\frac{rCD}{n}}%
\mathcal{K}_{1,M}\text{ }.  \label{M1KC_mino}
\end{equation}
\end{lemma}

\noindent Our argument leading to Lemma \ref{bound_moment_order_1} shows
that we have to assume that the constant $A_{\infty}$ introduced in (\ref%
{def_R_n_D_tild}) is large enough. In order to prove Lemma \ref%
{bound_moment_order_1} the following result is needed.

\begin{lemma}
\label{control_lower_infinity}Let $r>1,$ $\beta >0$ and $C\geq 0$. Assume
that (\textbf{H1}), (\textbf{H2}), (\textbf{H4}) and (\ref{def_A_cons}) hold
and let $\varphi =\left( \varphi _{k}\right) _{k=1}^{D}$ be an orthonormal
basis of $\left( M,\left\Vert 
\cdot%
\right\Vert _{2}\right) $ satisfying (\textbf{H4}). If positive constants $%
A_{+},$ $A_{-}$ and $A_{u}$ exist such that 
\begin{equation*}
A_{+}\frac{n}{\left( \ln n\right) ^{2}}\geq D\geq A_{-}\left( \ln n\right)
^{2},\text{ \ }rC\leq A_{u}\frac{D}{n}\text{ },
\end{equation*}%
and if%
\begin{equation*}
A_{\infty }\geq 32B_{2}A\sqrt{2A_{u}\beta }\sigma _{\min }^{-1}r_{M}\left(
\varphi \right)
\end{equation*}%
then for all $n\geq n_{0}\left( A_{-},A_{+},A,B_{2},r_{M}\left( \varphi
\right) ,\sigma _{\min },\beta \right) $, it holds%
\begin{equation*}
\mathbb{P}\left[ \max_{k\in \left\{ 1,...,D\right\} }\left\vert \frac{\sqrt{%
rC}\left( P_{n}-P\right) \left( \psi _{1,M}%
\cdot%
\varphi _{k}\right) }{\sqrt{\sum_{j=1}^{D}\left( P_{n}-P\right) ^{2}\left(
\psi _{1,M}%
\cdot%
\varphi _{j}\right) }}\right\vert \geq \frac{\tilde{R}_{n,D,\alpha }}{%
r_{M}\left( \varphi \right) \sqrt{D}}\right] \leq \frac{2D+1}{n^{\beta }}%
\text{ }.
\end{equation*}
\end{lemma}

\subsection*{}%
\textbf{Proof of Lemma \ref{control_lower_infinity}.} By Cauchy-Schwarz
inequality, we get 
\begin{equation*}
\chi _{M}=\sqrt{\sum_{k=1}^{D}\left( P_{n}-P\right) ^{2}\left( \psi _{1,M}%
\cdot%
\varphi _{k}\right) }=\sup_{s\in S_{M}}\left\vert \left( P_{n}-P\right)
\left( \psi _{1,M}%
\cdot%
s\right) \right\vert \text{ },
\end{equation*}%
where $S_{M}$ is the unit sphere of $M$, that is%
\begin{equation*}
S_{M}=\left\{ s\in M,\text{ }s=\sum_{k=1}^{D}\beta _{k}\varphi _{k}\text{
and }\sqrt{\sum_{k=1}^{D}\beta _{k}^{2}}=1\right\} \text{ }.
\end{equation*}%
Thus we can apply Klein-Rio's inequality (\ref{klein_rio_2}) to $\chi _{M}$
by taking $\mathcal{F=}S_{M}$ and use the fact that%
\begin{align}
\sup_{s\in S_{M}}\left\Vert \psi _{1,M}%
\cdot%
s-P\left( \psi _{1,M}%
\cdot%
s\right) \right\Vert _{\infty }& \leq 4A\sqrt{D}r_{M}\left( \varphi \right) 
\text{ \ \ \ by (\ref{center_proof}), (\ref{control_phi_1M_proof})\ and (%
\textbf{H4}).}  \label{ineg_1_bis_lemma-31} \\
\sup_{s\in S_{M}}%
\var%
\left( \psi _{1,M}%
\cdot%
s\right) & =\sup_{s\in S_{M}}P\left( \psi _{1,M}%
\cdot%
s\right) ^{2}\leq 16A^{2}\text{ \ \ \ by (\ref{center_proof}), (\ref%
{control_phi_1M_proof})}  \notag
\end{align}%
and also, by using (\ref{ineg_1_bis_lemma-31}) in Inequality (\ref{emp_hoff}%
) applied to $\chi _{M}$, we get that%
\begin{align*}
\mathbb{E}\left[ \chi _{M}\right] & \geq B_{2}^{-1}\sqrt{\mathbb{E}\left[
\chi _{M}^{2}\right] }-\frac{4A\sqrt{D}r_{M}\left( \varphi \right) }{n} \\
& =B_{2}^{-1}\sqrt{\frac{D}{n}}\mathcal{K}_{1,M}-\frac{4A\sqrt{D}r_{M}\left(
\varphi \right) }{n}\text{ }.
\end{align*}%
We thus obtain by (\ref{klein_rio_2}), for all $\varepsilon ,x>0,$%
\begin{equation}
\mathbb{P}\left( \chi _{M}\leq \left( 1-\varepsilon \right) B_{2}^{-1}\sqrt{%
\frac{D}{n}}\mathcal{K}_{1,M}-\sqrt{32A^{2}\frac{x}{n}}-\left( 1-\varepsilon
+\left( 1+\frac{1}{\varepsilon }\right) x\right) \frac{4A\sqrt{D}r_{M}\left(
\varphi \right) }{n}\right) \leq \exp \left( -x\right) \text{ }.
\label{chi_1_n_concen}
\end{equation}%
So, by taking $\varepsilon =\frac{1}{2}$ and $x=\beta \ln n$ in (\ref%
{chi_1_n_concen}), and by observing that $D\geq A_{-}\left( \ln n\right)
^{2} $ and $\mathcal{K}_{1,M}\geq 2\sigma _{\min }$, we conclude that, for
all $n\geq n_{0}\left( A_{-},A,B_{2},r_{M}\left( \varphi \right) ,\sigma
_{\min },\beta \right) $,%
\begin{equation}
\mathbb{P}\left[ \chi _{M}\leq \frac{B_{2}^{-1}}{8}\sqrt{\frac{D}{n}}%
\mathcal{K}_{1,M}\right] \leq n^{-\beta }\text{ }.  \label{chi_1_n}
\end{equation}%
Furthermore, combining Bernstein's inequality (\ref{bernstein_ineq_reg}),
with the observation that we have, for every $k\in \left\{ 1,...,D\right\} $,%
\begin{align*}
\left\Vert \psi _{1,M}%
\cdot%
\varphi _{k}\right\Vert _{\infty }& \leq 4A\sqrt{D}r_{M}\left( \varphi
\right) \text{ \ \ \ \ \ \ \ \ \ \ \ \ \ by (\ref{control_phi_1M_proof}) \
and \ (\textbf{H4})} \\
P\left( \psi _{1,M}%
\cdot%
\varphi _{k}\right) ^{2}& \leq \left\Vert \psi _{1,M}\right\Vert _{\infty
}^{2}\leq 16A^{2}\text{ \ \ \ \ \ by (\ref{control_phi_1M_proof})}
\end{align*}%
we get that, for every $x>0$ and every $k\in \left\{ 1,...,D\right\} $,%
\begin{equation*}
\mathbb{P}\left[ \left\vert \left( P_{n}-P\right) \left( \psi _{1,M}%
\cdot%
\varphi _{k}\right) \right\vert \geq \sqrt{32A^{2}\frac{x}{n}}+\frac{4A\sqrt{%
D}r_{M}\left( \varphi \right) }{3}\frac{x}{n}\right] \leq 2\exp \left(
-x\right)
\end{equation*}%
and so 
\begin{equation}
\mathbb{P}\left[ \max_{k\in \left\{ 1,...,D\right\} }\left\vert \left(
P_{n}-P\right) \left( \psi _{1,M}%
\cdot%
\varphi _{k}\right) \right\vert \geq \sqrt{32A^{2}\frac{x}{n}}+\frac{4A\sqrt{%
D}r_{M}\left( \varphi \right) }{3}\frac{x}{n}\right] \leq 2D\exp \left(
-x\right) \text{ }.  \label{max_1_n_x}
\end{equation}%
Hence, taking $x=\beta \ln n$ in (\ref{max_1_n_x}), it comes 
\begin{equation}
\mathbb{P}\left[ \max_{k\in \left\{ 1,...,D\right\} }\left\vert \left(
P_{n}-P\right) \left( \psi _{1,M}%
\cdot%
\varphi _{k}\right) \right\vert \geq \sqrt{\frac{32A^{2}\beta \ln n}{n}}+%
\frac{4A\sqrt{D}r_{M}\left( \varphi \right) \beta \ln n}{3n}\right] \leq 
\frac{2D}{n^{\beta }}\text{ },  \label{max_1_n}
\end{equation}%
then, by using (\ref{chi_1_n}) and (\ref{max_1_n}), we get for all $n\geq
n_{0}\left( A_{-},A,B_{2},r_{M}\left( \varphi \right) ,\sigma _{\min },\beta
\right) $,%
\begin{equation*}
\mathbb{P}\left[ \max_{k\in \left\{ 1,...,D\right\} }\left\vert \frac{\sqrt{%
rC}\left( P_{n}-P\right) \left( \psi _{1,M}%
\cdot%
\varphi _{k}\right) }{\chi _{M}}\right\vert \geq \frac{8B_{2}\sqrt{rC}}{%
\sqrt{\frac{D}{n}}\mathcal{K}_{1,M}}\left( \sqrt{\frac{32A^{2}\beta \ln n}{n}%
}+\frac{4A\sqrt{D}r_{M}\left( \varphi \right) \beta \ln n}{3n}\right) \right]
\leq \frac{2D+1}{n^{\beta }}\text{ }.
\end{equation*}%
Finally, as $A_{+}\frac{n}{\left( \ln n\right) ^{2}}\geq D$ we have, for all 
$n\geq n_{0}\left( A,A_{+},r_{M}\left( \varphi \right) ,\beta \right) $, 
\begin{equation*}
\frac{4A\sqrt{D}r_{M}\left( \varphi \right) \beta \ln n}{3n}\leq \sqrt{\frac{%
32A^{2}\beta \ln n}{n}}
\end{equation*}%
and we can check that, since $rC\leq A_{u}\frac{D}{n}$ and $\mathcal{K}%
_{1,M}\geq 2\sigma _{\min }$, if%
\begin{equation*}
A_{\infty }\geq 32B_{2}\sqrt{2A_{u}A^{2}\beta }\sigma _{\min
}^{-1}r_{M}\left( \varphi \right)
\end{equation*}%
then, for all $n\geq n_{0}\left( A_{-},A_{+},A,B_{2},r_{M}\left( \varphi
\right) ,\sigma _{\min },\beta \right) $,%
\begin{equation*}
\mathbb{P}\left[ \max_{k\in \left\{ 1,...,D\right\} }\left\vert \frac{\sqrt{%
rC}\left( P_{n}-P\right) \left( \psi _{1,M}%
\cdot%
\varphi _{k}\right) }{\chi _{M}}\right\vert \geq \frac{A_{\infty }}{%
r_{M}\left( \varphi \right) }\sqrt{\frac{\ln n}{n}}\right] \leq \frac{2D+1}{%
n^{\beta }}
\end{equation*}%
which readily gives the result. 
$\blacksquare$%

\noindent We are now ready to prove the lower bound (\ref{M1KC_mino}) for
the expected value of the largest increment of the empirical process over $%
\mathcal{F}_{\left( C,rC\right] }.$

\subsection*{}%
\textbf{Proof of Lemma \ref{bound_moment_order_1}.} Let us begin with the
lower bound of 
\begin{equation*}
\mathbb{E}^{\frac{1}{2}}\left( \sup_{s\in \mathcal{F}_{\left( C,rC\right]
}}\left( P_{n}-P\right) \left( \psi _{1,M}%
\cdot%
\left( s_{M}-s\right) \right) \right) ^{2}\text{ },
\end{equation*}%
a result that will be need further in the proof. Introduce for all $k\in
\left\{ 1,...,D\right\} $,%
\begin{equation*}
\beta _{k,n}=\frac{\sqrt{rC}\left( P_{n}-P\right) \left( \psi _{1,M}%
\cdot%
\varphi _{k}\right) }{\sqrt{\sum_{j=1}^{D}\left( P_{n}-P\right) ^{2}\left(
\psi _{1,M}%
\cdot%
\varphi _{j}\right) }}\text{ },
\end{equation*}%
and observe that the excess risk on $M$ of $\left( \sum_{k=1}^{D}\beta
_{k,n}\varphi _{k}+s_{M}\right) \in M$ is equal to $rC$. We also set 
\begin{equation*}
\tilde{\Omega}=\left\{ \max_{k\in \left\{ 1,...,D\right\} }\left\vert \beta
_{k,n}\right\vert \leq \frac{\tilde{R}_{n,D,\alpha }}{r_{M}\left( \varphi
\right) \sqrt{D}}\right\} \text{ }.
\end{equation*}%
By Lemma \ref{control_lower_infinity} we have that for all $\beta >0$, if $%
A_{\infty }\geq 32B_{2}\sqrt{2A_{u}A^{2}\beta }\sigma _{\min
}^{-1}r_{M}\left( \varphi \right) $ then,

\noindent for all $n\geq n_{0}\left( A_{-},A_{+},A,B_{2},r_{M}\left( \varphi
\right) ,\sigma _{\min },\beta \right) $,%
\begin{equation}
\mathbb{P}\left( \tilde{\Omega}\right) \geq 1-\frac{2D+1}{n^{\beta }}\text{ }%
.  \label{proba_omega_tilde_reg_fixed}
\end{equation}%
Moreover, by (\textbf{H4}), we get on the event $\tilde{\Omega}$,%
\begin{equation*}
\left\Vert \sum_{k=1}^{D}\beta _{k,n}\varphi _{k}\right\Vert _{\infty }\leq 
\tilde{R}_{n,D,\alpha }\text{ },
\end{equation*}%
and so, on $\tilde{\Omega}$, 
\begin{equation}
\left( s_{M}+\sum_{k=1}^{D}\beta _{k,n}\varphi _{k}\right) \in \mathcal{F}%
_{\left( C,rC\right] }\text{ }.  \label{appart_FRC}
\end{equation}%
As a consequence, by (\ref{appart_FRC}) it holds%
\begin{align}
& \mathbb{E}^{\frac{1}{2}}\left( \sup_{s\in \mathcal{F}_{\left( C,rC\right]
}}\left( P_{n}-P\right) \left( \psi _{1,M}%
\cdot%
\left( s_{M}-s\right) \right) \right) ^{2}  \notag \\
& \geq \mathbb{E}^{\frac{1}{2}}\left[ \left( \left( P_{n}-P\right) \left(
\psi _{1,M}%
\cdot%
\left( \sum_{k=1}^{D}\beta _{k,n}\varphi _{k}\right) \right) \right) ^{2}%
\mathbf{1}_{\tilde{\Omega}}\right]  \notag \\
& =\sqrt{rC}\sqrt{\mathbb{E}\left[ \left( \sum_{k=1}^{D}\left(
P_{n}-P\right) ^{2}\left( \psi _{1,M}%
\cdot%
\varphi _{k}\right) \right) \mathbf{1}_{\tilde{\Omega}}\right] }\text{ }.
\label{minoration_mean_infinity}
\end{align}%
Furthermore, since by (\ref{center_proof}) $P\left( \psi _{1,M}%
\cdot%
\varphi _{k}\right) =0$ and by (\textbf{H4}) $\left\Vert \varphi
_{k}\right\Vert _{\infty }\leq \sqrt{D}r_{M}\left( \varphi \right) $ for all 
$k\in \left\{ 1,...,D\right\} ,$ we have%
\begin{align*}
\left\vert \sum_{k=1}^{D}\left( P_{n}-P\right) ^{2}\left( \psi _{1,M}%
\cdot%
\varphi _{k}\right) \right\vert & \leq D\max_{k=1,...,D}\left\vert \left(
P_{n}-P\right) ^{2}\left( \psi _{1,M}%
\cdot%
\varphi _{k}\right) \right\vert \\
& =D\max_{k=1,...,D}\left\vert P_{n}^{2}\left( \psi _{1,M}%
\cdot%
\varphi _{k}\right) \right\vert \\
& \leq D\max_{k=1,...,D}\left\Vert \psi _{1,M}%
\cdot%
\varphi _{k}\right\Vert _{\infty }^{2} \\
& \leq 16A^{2}D^{2}r_{M}^{2}\left( \varphi \right)
\end{align*}%
and it ensures 
\begin{equation}
\mathbb{E}\left[ \left( \sum_{k=1}^{D}\left( P_{n}-P\right) ^{2}\left( \psi
_{1,M}%
\cdot%
\varphi _{k}\right) \right) 1_{\tilde{\Omega}}\right] \geq \mathbb{E}\left[
\left( \sum_{k=1}^{D}\left( P_{n}-P\right) ^{2}\left( \psi _{1,M}%
\cdot%
\varphi _{k}\right) \right) \right] -16A^{2}D^{2}r_{M}^{2}\left( \varphi
\right) \mathbb{P}\left[ \left( \tilde{\Omega}\right) ^{c}\right] \text{ }.
\label{mino_ordre_2_omega}
\end{equation}%
Comparing inequality (\ref{mino_ordre_2_omega}) with (\ref%
{minoration_mean_infinity}) and using (\ref{proba_omega_tilde_reg_fixed}),
we obtain the following lower bound for all $n\geq n_{0}\left(
A_{-},A_{+},A,B_{2},r_{M}\left( \varphi \right) ,\sigma _{\min },\beta
\right) $,%
\begin{align}
\mathbb{E}^{\frac{1}{2}}\left( \sup_{s\in \mathcal{F}_{\left( C,rC\right]
}}\left( P_{n}-P\right) \left( \psi _{1,M}%
\cdot%
\left( s_{M}-s\right) \right) \right) ^{2}& \geq \sqrt{rC}\sqrt{\mathbb{E}%
\left[ \left( \sum_{k=1}^{D}\left( P_{n}-P\right) ^{2}\left( \psi _{1,M}%
\cdot%
\varphi _{k}\right) \right) \right] }  \notag \\
& -4Ar_{M}\left( \varphi \right) D\sqrt{rC}\sqrt{\mathbb{P}\left[ \left( 
\tilde{\Omega}\right) ^{c}\right] }  \notag \\
& \geq \sqrt{\frac{rCD}{n}}\mathcal{K}_{1,M}-4Ar_{M}\left( \varphi \right) D%
\sqrt{rC}\sqrt{\frac{2D+1}{n^{\beta }}}\text{ }.  \label{mino_ordre_2_beta}
\end{align}%
We take $\beta =4$, and we must have%
\begin{equation*}
A_{\infty }\geq 64AB_{2}\sqrt{2A_{u}}\sigma _{\min }^{-1}r_{M}\left( \varphi
\right) \text{ .}
\end{equation*}%
Since $D\leq A_{+}n\left( \ln n\right) ^{-2}$ and $\mathcal{K}_{1,M}\geq
2\sigma _{\min }$ under (\textbf{H2}), we get, for all $n\geq n_{0}\left(
A,A_{+},r_{M}\left( \varphi \right) ,\sigma _{\min }\right) $,%
\begin{equation}
4Ar_{M}\left( \varphi \right) D\sqrt{rC}\sqrt{\frac{2D+1}{n^{\beta }}}\leq 
\frac{1}{\sqrt{D}}\times \sqrt{\frac{rCD}{n}}\mathcal{K}_{1,M}
\label{compa_reg}
\end{equation}%
and so, by combining (\ref{mino_ordre_2_beta}) and (\ref{compa_reg}), for
all $n\geq n_{0}\left( A_{-},A_{+},A,B_{2},r_{M}\left( \varphi \right)
,\sigma _{\min }\right) $, it holds%
\begin{equation}
\mathbb{E}^{\frac{1}{2}}\left( \sup_{s\in \mathcal{F}_{\left( C,rC\right]
}}\left( P_{n}-P\right) \left( \psi _{1,M}%
\cdot%
\left( s_{M}-s\right) \right) \right) ^{2}\geq \left( 1-\frac{1}{\sqrt{D}}%
\right) \sqrt{\frac{rCD}{n}}\mathcal{K}_{1,M}\text{ }.
\label{mino_moment_2_proof}
\end{equation}%
Now, as $D\geq A_{-}\left( \ln n\right) ^{2}$ we have for all $n\geq
n_{0}\left( A_{-}\right) $, $D^{-1/2}\leq 1/2$. Moreover, we have $\mathcal{K%
}_{1,M}\geq 2\sigma _{\min }$ by (\textbf{H2}) and $rC\geq A_{l}Dn^{-1}$, so
we finally deduce from (\ref{mino_moment_2_proof}) that, for all $n\geq
n_{0}\left( A_{-},A_{+},A,B_{2},A_{l},r_{M}\left( \varphi \right) ,\sigma
_{\min }\right) $,%
\begin{equation}
\mathbb{E}^{\frac{1}{2}}\left( \sup_{s\in \mathcal{F}_{\left( C,rC\right]
}}\left( P_{n}-P\right) \left( \psi _{1,M}%
\cdot%
\left( s_{M}-s\right) \right) \right) ^{2}\geq \sigma _{\min }\sqrt{A_{l}}%
\frac{D}{n}\text{ }.  \label{mino_moment_2_cor_28}
\end{equation}

\noindent We turn now to the lower bound of $\mathbb{E}\left[ \sup_{s\in 
\mathcal{F}_{\left( C,rC\right] }}\left( P_{n}-P\right) \left( \psi _{1,M}%
\cdot%
\left( s_{M}-s\right) \right) \right] $. First observe that $s\in \mathcal{F}%
_{\left( C,rC\right] }$ implies that $\left( 2s_{M}-s\right) \in \mathcal{F}%
_{\left( C,rC\right] }$, so that%
\begin{equation}
\mathbb{E}\left[ \sup_{s\in \mathcal{F}_{\left( C,rC\right] }}\left(
P_{n}-P\right) \left( \psi _{1,M}%
\cdot%
\left( s_{M}-s\right) \right) \right] =\mathbb{E}\left[ \sup_{s\in \mathcal{F%
}_{\left( C,rC\right] }}\left\vert \left( P_{n}-P\right) \left( \psi _{1,M}%
\cdot%
\left( s_{M}-s\right) \right) \right\vert \right] .  \label{abs_values_proof}
\end{equation}%
In the next step, we apply Corollary \ref{cor_alt_hoff_2}. More precisely,
using notations of Corollary \ref{cor_alt_hoff_2}, we set 
\begin{equation*}
\mathcal{F=}\left\{ \psi _{1,M}%
\cdot%
\left( s_{M}-s\right) \text{ };\text{ }s\in \mathcal{F}_{\left( C,rC\right]
}\right\}
\end{equation*}%
and%
\begin{equation*}
Z=\sup_{s\in \mathcal{F}_{\left( C,rC\right] }}\left\vert \left(
P_{n}-P\right) \left( \psi _{1,M}%
\cdot%
\left( s_{M}-s\right) \right) \right\vert \text{ }.
\end{equation*}%
Now, since for all $n\geq n_{0}\left( A_{+},A_{-},A_{\infty
},A_{cons}\right) $ we have $\tilde{R}_{n,D,\alpha }\leq 1$, we get by (\ref%
{center_proof}) and (\ref{control_phi_1M_proof}), for all $n\geq n_{0}\left(
A_{+},A_{-},A_{\infty },A_{cons}\right) $,%
\begin{equation*}
\sup_{f\in \mathcal{F}}\left\Vert f-Pf\right\Vert _{\infty }=\sup_{s\in 
\mathcal{F}_{\left( C,rC\right] }}\left\Vert \psi _{1,M}%
\cdot%
\left( s_{M}-s\right) \right\Vert _{\infty }\leq 4A\tilde{R}_{n,D,\alpha
}\leq 4A
\end{equation*}%
we set $b=4A$. Since we assume that $rC\leq A_{u}\frac{D}{n}$, it moreover
holds by (\ref{control_phi_1M_proof}), 
\begin{equation*}
\sup_{f\in \mathcal{F}}%
\var%
\left( f\right) \leq \sup_{s\in \mathcal{F}_{\left( C,rC\right] }}P\left(
\psi _{1,M}%
\cdot%
\left( s_{M}-s\right) \right) ^{2}\leq 16A^{2}rC\leq 16A^{2}A_{u}\frac{D}{n}
\end{equation*}%
and so we set $\sigma ^{2}=16A^{2}A_{u}\frac{D}{n}$. Now, by (\ref%
{mino_moment_2_cor_28}) we have, for all $n\geq n_{0}\left(
A_{-},A_{+},A,B_{2},A_{l},r_{M}\left( \varphi \right) ,\sigma _{\min
}\right) $, 
\begin{equation}
\sqrt{\mathbb{E}\left[ Z^{2}\right] }\geq \sigma _{\min }\sqrt{A_{l}}\frac{D%
}{n}\text{ }.  \label{lower_moment_order_2}
\end{equation}%
Hence, a positive constant $L_{A,A_{l},A_{u},\sigma _{\min }}$ ( $\max
\left( 4A\sqrt{A_{u}}A_{l}^{-1/2}\sigma _{\min }^{-1}\text{ };\text{ }2\sqrt{%
A}A_{l}^{-1/4}\sigma _{\min }^{-1/2}\right) $ holds) exists such that, by
setting 
\begin{equation*}
\varkappa _{n}=\frac{L_{A,A_{l},A_{u},\sigma _{\min }}}{\sqrt{D}}
\end{equation*}%
we get, using (\ref{lower_moment_order_2}), that, for all $n\geq n_{0}\left(
A_{-},A_{+},A_{l},A_{u},A,B_{2},r_{M}\left( \varphi \right) ,A_{cons},\sigma
_{\min }\right) $,%
\begin{equation*}
\varkappa _{n}^{2}\mathbb{E}\left[ Z^{2}\right] \geq \frac{\sigma ^{2}}{n}%
\text{ ,}
\end{equation*}%
\begin{equation*}
\varkappa _{n}^{2}\sqrt{\mathbb{E}\left[ Z^{2}\right] }\geq \frac{b}{n}\text{
.}
\end{equation*}%
Furthermore, since $D\geq A_{-}\left( \ln n\right) ^{2}$, we have for all $%
n\geq n_{0}\left( A_{-},A,A_{u},A_{l},\sigma _{\min }\right) $,%
\begin{equation*}
\varkappa _{n}\in \left( 0,1\right) \text{ }.
\end{equation*}%
So, using (\ref{abs_values_proof}) and Corollary \ref{cor_alt_hoff_2}, it
holds for all $n\geq n_{0}\left( A_{-},A_{+},A_{l},A_{u},A,B_{2},r_{M}\left(
\varphi \right) ,\sigma _{\min }\right) $,%
\begin{align}
& \mathbb{E}\left[ \sup_{s\in \mathcal{F}_{\left( C,rC\right] }}\left(
P_{n}-P\right) \left( \psi _{1,M}%
\cdot%
\left( s_{M}-s\right) \right) \right]  \notag \\
& \geq \left( 1-\frac{L_{A,A_{l},A_{u},\sigma _{\min }}}{\sqrt{D}}\right) 
\mathbb{E}^{\frac{1}{2}}\left( \sup_{s\in \mathcal{F}_{\left( C,rC\right]
}}\left( P_{n}-P\right) \left( \psi _{1,M}%
\cdot%
\left( s_{M}-s\right) \right) \right) ^{2}\text{ }.  \label{mino_m1kc_proof}
\end{align}%
Finally, by comparing (\ref{mino_moment_2_proof}) and (\ref{mino_m1kc_proof}%
), we deduce that for all $n\geq n_{0}\left(
A_{-},A_{+},A_{l},A_{u},A,B_{2},r_{M}\left( \varphi \right) ,\sigma _{\min
}\right) $,%
\begin{equation*}
\mathbb{E}\left[ \sup_{s\in \mathcal{F}_{\left( C,rC\right] }}\left(
P_{n}-P\right) \left( \psi _{1,M}%
\cdot%
\left( s_{M}-s\right) \right) \right] \geq \left( 1-\frac{%
L_{A,A_{l},A_{u},\sigma _{\min }}}{\sqrt{D}}\right) \sqrt{\frac{rCD}{n}}%
\mathcal{K}_{1,M}
\end{equation*}%
and so (\ref{M1KC_mino}) is proved.\ 
$\blacksquare$%

\noindent Let us now turn to the control of second order terms appearing in
the expansion of the least squares contrast, see (\ref{dev_contrast}). Let
us define 
\begin{equation*}
\Omega _{C}\left( x\right) =\sup_{s\in \mathcal{F}_{\left( C,rC\right]
}}\left\{ \frac{\left\vert \psi _{2}\left( \left( s-s_{M}\right) \left(
x\right) \right) -\psi _{2}\left( \left( t-s_{M}\right) \left( x\right)
\right) \right\vert }{\left\vert s\left( x\right) -t\left( x\right)
\right\vert }\text{ };\text{ }\left( s,t\right) \in \mathcal{F}_{C}\text{ },%
\text{ }s\left( x\right) \neq t\left( x\right) \right\} .
\end{equation*}%
After straightforward computations using that $\psi _{2}\left( t\right)
=t^{2}$ for all $t\in \mathbb{R}$ and assuming (\textbf{H3}), we get that,
for all $x\in \mathcal{X}$,%
\begin{align}
\Omega _{C}\left( x\right) & =2\sup_{s\in \mathcal{F}_{C}}\left\{ \left\vert
s\left( x\right) -s_{M}\left( x\right) \right\vert \right\}
\label{omega_c_z} \\
& \leq 2\left( \tilde{R}_{n,D,\alpha }\wedge \sqrt{CD}A_{3,M}\right) \text{ }%
.  \label{majo_omega_c_z}
\end{align}

\begin{lemma}
\label{lower_rademacher copy(1)}Let $C\geq 0$. Under (\textbf{H3}), it holds%
\begin{equation*}
\mathbb{E}\left[ \sup_{s\in \mathcal{F}_{C}}\left\vert \left( P_{n}-P\right)
\left( \psi _{2}\circ \left( s-s_{M}\right) \right) \right\vert \right] \leq
8\sqrt{\frac{CD}{n}}\left( \tilde{R}_{n,D,\alpha }\wedge \sqrt{CD}%
A_{3,M}\right) \text{ }.
\end{equation*}
\end{lemma}

\subsection*{}%
\textbf{Proof. }We define the Rademacher process $\mathcal{R}_{n}$ on a
class $\mathcal{F}$ of measurable functions from $\mathcal{X}$ to $\mathbb{R}
$, to be%
\begin{equation*}
\mathcal{R}_{n}\left( f\right) =\frac{1}{n}\sum_{i=1}^{n}\varepsilon
_{i}f\left( X_{i}\right) \text{ , \ \ }f\in \mathcal{F}
\end{equation*}%
where $\varepsilon _{i}$ are independent Rademacher random variables also
independent from the $X_{i}$. By the usual symmetrization argument we have%
\begin{equation*}
\mathbb{E}\left[ \sup_{s\in \mathcal{F}_{C}}\left\vert \left( P_{n}-P\right)
\left( \psi _{2}\circ \left( s-s_{M}\right) \right) \right\vert \right] \leq
2\mathbb{E}\left[ \sup_{s\in \mathcal{F}_{C}}\left\vert \mathcal{R}%
_{n}\left( \psi _{2}\circ \left( s-s_{M}\right) \right) \right\vert \right] .
\end{equation*}%
Taking the expectation with respect to the Rademacher variables, we get%
\begin{align}
& \mathbb{E}_{\varepsilon }\left[ \underset{s\in \mathcal{F}_{C}}{\sup }%
\left\vert \mathcal{R}_{n}\left( \psi _{2}\circ \left( s-s_{M}\right)
\right) \right\vert \right]  \notag \\
& =\mathbb{E}_{\varepsilon }\left[ \underset{s\in \mathcal{F}_{C}}{\sup }%
\left\vert \mathcal{R}_{n}\left( \left( s-s_{M}\right) ^{2}\right)
\right\vert \right]  \notag \\
& \leq \left( \underset{1\leq i\leq n}{\max }\Omega _{C}\left( X_{i}\right)
\right) \mathbb{E}_{\varepsilon }\left[ \sup_{s\in \mathcal{F}%
_{C}}\left\vert \frac{1}{n}\sum_{i=1}^{n}\varepsilon _{i}\varphi _{i}\left(
\left( s-s_{M}\right) \left( X_{i}\right) \right) \right\vert \right]
\label{rad_1}
\end{align}%
where the functions $\varphi _{i}:\mathbb{R\longrightarrow R}$ are defined
by 
\begin{equation*}
\varphi _{i}\left( t\right) =\left\{ 
\begin{tabular}{l}
$\left( \Omega _{C}\left( X_{i}\right) \right) ^{-1}t^{2}$ \ \ \ \ for $%
\left\vert t\right\vert \leq \sup_{s\in \mathcal{F}_{C}}\left\{ \left\vert
s\left( X_{i}\right) -s_{M}\left( X_{i}\right) \right\vert \right\} =\frac{%
\Omega _{C}\left( X_{i}\right) }{2}$ \\ 
$\frac{1}{4}\Omega _{C}\left( X_{i}\right) $\ \ \ \ \ \ \ \ \ \ \ otherwise%
\end{tabular}%
\ \right.
\end{equation*}%
Then by (\ref{omega_c_z}) we deduce that $\varphi _{i}$ is a contraction
mapping with $\varphi _{i}\left( 0\right) =0$. We thus apply Theorem \ref%
{comparison_theorem} to get%
\begin{align}
\mathbb{E}_{\varepsilon }\left[ \sup_{s\in \mathcal{F}_{C}}\left\vert \frac{1%
}{n}\sum_{i=1}^{n}\varepsilon _{i}\varphi _{i}\left( \left( s-s_{M}\right)
\left( X_{i}\right) \right) \right\vert \right] & \leq 2\mathbb{E}%
_{\varepsilon }\left[ \sup_{s\in \mathcal{F}_{C}}\left\vert \frac{1}{n}%
\sum_{i=1}^{n}\varepsilon _{i}\left( s-s_{M}\right) \left( X_{i}\right)
\right\vert \right]  \notag \\
& =2\mathbb{E}_{\varepsilon }\left[ \sup_{s\in \mathcal{F}_{C}}\left\vert 
\mathcal{R}_{n}\left( s-s_{M}\right) \right\vert \right]  \label{rad_2}
\end{align}%
and so we derive successively the following upper bounds in mean,%
\begin{align*}
& \mathbb{E}\left[ \sup_{s\in \mathcal{F}_{C}}\left\vert \mathcal{R}%
_{n}\left( \psi _{2}\circ \left( s-s_{M}\right) \right) \right\vert \right] =%
\mathbb{E}\left[ \mathbb{E}_{\varepsilon }\left[ \sup_{s\in \mathcal{F}%
_{C}}\left\vert \mathcal{R}_{n}\left( \psi _{2}\circ \left( s-s_{M}\right)
\right) \right\vert \right] \right] \\
& \leq \mathbb{E}\left[ \left( \underset{1\leq i\leq n}{\max }\Omega
_{C}\left( X_{i}\right) \right) \mathbb{E}_{\varepsilon }\left[ \sup_{s\in 
\mathcal{F}_{C}}\left\vert \frac{1}{n}\sum_{i=1}^{n}\varepsilon _{i}\varphi
_{i}\left( \left( s-s_{M}\right) \left( X_{i}\right) \right) \right\vert %
\right] \right] \text{ \ \ \ by (\ref{rad_1})} \\
& \leq 2\mathbb{E}\left[ \left( \underset{1\leq i\leq n}{\max }\Omega
_{C}\left( X_{i}\right) \right) \mathbb{E}_{\varepsilon }\left[ \sup_{s\in 
\mathcal{F}_{C}}\left\vert \mathcal{R}_{n}\left( s-s_{M}\right) \right\vert %
\right] \right] \text{ \ \ \ \ \ \ \ \ \ \ \ \ \ \ \ \ \ \ \ \ \ \ \ \ by (%
\ref{rad_2})} \\
& =2\mathbb{E}\left[ \left( \underset{1\leq i\leq n}{\max }\Omega _{C}\left(
X_{i}\right) \right) \sup_{s\in \mathcal{F}_{C}}\left\vert \mathcal{R}%
_{n}\left( s-s_{M}\right) \right\vert \right] \\
& \leq 2\sqrt{\mathbb{E}\left[ \underset{1\leq i\leq n}{\max }\Omega
_{C}^{2}\left( X_{i}\right) \right] }\sqrt{\mathbb{E}\left[ \left(
\sup_{s\in \mathcal{F}_{C}}\left\vert \mathcal{R}_{n}\left( s-s_{M}\right)
\right\vert \right) ^{2}\right] }
\end{align*}%
We consider now an orthonormal basis of $\left( M,\left\Vert 
\cdot%
\right\Vert _{2}\right) $ and denote it by $\left( \varphi _{k}\right)
_{k=1}^{D}$. Whence%
\begin{align*}
& \sqrt{\mathbb{E}\left[ \left( \sup_{s\in \mathcal{F}_{C}}\left\vert 
\mathcal{R}_{n}\left( s-s_{M}\right) \right\vert \right) ^{2}\right] } \\
& \leq \sqrt{\mathbb{E}\left[ \left( \sup \left\{ \left\vert
\sum_{k=1}^{D}a_{k}\mathcal{R}_{n}\left( \varphi _{k}\right) \right\vert ;%
\underset{k=1}{\overset{D}{\sum }}a_{k}^{2}\leq C\right\} \right) ^{2}\right]
} \\
& =\sqrt{C}\sqrt{\mathbb{E}\left[ \sum_{k=1}^{D}\left( \mathcal{R}_{n}\left(
\varphi _{k}\right) \right) ^{2}\right] }=\sqrt{\frac{CD}{n}}\text{ },
\end{align*}%
to complete the proof, it remains to observe that, by (\ref{majo_omega_c_z}),%
\begin{equation*}
\sqrt{\mathbb{E}\left[ \underset{1\leq i\leq n}{\max }\Omega _{C}^{2}\left(
X_{i}\right) \right] }\leq 2\left( \tilde{R}_{n,D,\alpha }\wedge \sqrt{CD}%
A_{3,M}\right) \text{ }.
\end{equation*}

\ \ \ \ \ \ \ \ \ \ \ \ \ \ \ \ \ \ \ \ \ \ \ \ \ \ \ \ \ \ \ \ \ \ \ \ \ \
\ \ \ \ \ \ \ \ \ \ \ \ \ \ \ \ \ \ \ \ \ \ \ \ \ \ \ \ \ \ \ \ \ \ \ \ \ \
\ \ \ \ \ \ \ \ \ \ \ \ \ \ \ \ \ \ \ \ \ \ \ \ \ \ \ \ \ \ \ \ \ \ \ \ \ \
\ \ \ \ \ \ \ \ \ \ \ \ \ \ \ \ \ \ \ \ \ \ \ \ \ \ \ \ \ \ \ \ 
$\blacksquare$%

\noindent In the following Lemma, we provide uniform upper bounds for the
supremum of the empirical process of second order terms in the contrast
expansion when the considered slices are not too small.

\begin{lemma}
\label{majo_uniform_s2c}Let $A_{+},A_{-},A_{l},\beta ,C_{-}>0$, and assume (%
\textbf{H3}) and (\ref{def_A_cons}). If $C_{-}\geq A_{l}\frac{D}{n}$ and $%
A_{+}n\left( \ln n\right) ^{-2}\geq D\geq A_{-}\left( \ln n\right) ^{2}$,
then a positive constant $L_{A_{-},A_{l},\beta }$ exists such that, for all $%
n\geq n_{0}\left( A_{\infty },A_{cons},A_{+},A_{l}\right) $,%
\begin{equation*}
\mathbb{P}\left[ \forall C>C_{-},\text{ \ }\sup_{s\in \mathcal{F}%
_{C}}\left\vert \left( P_{n}-P\right) \left( \psi _{2}\circ \left(
s-s_{M}\right) \right) \right\vert \leq L_{A_{-},A_{l},\beta }\sqrt{\frac{CD%
}{n}}\tilde{R}_{n,D,\alpha }\right] \geq 1-n^{-\beta }\text{ }.
\end{equation*}
\end{lemma}

\subsection*{}%
\textbf{Proof.}\ First notice that, as $A_{+}n\left( \ln n\right) ^{-2}\geq
D $, we have by (\ref{def_A_cons}), 
\begin{equation*}
\tilde{R}_{n,D,\alpha }\leq \frac{\max \left\{ A_{cons}\text{ };\text{ }%
A_{\infty }\sqrt{A_{+}}\right\} }{\sqrt{\ln n}}\text{ }.
\end{equation*}%
By consequence, for all $n\geq n_{0}\left( A_{\infty },A_{cons},A_{+}\right) 
$,%
\begin{equation}
\tilde{R}_{n,D,\alpha }\leq 1\text{ }.  \label{control_R_n_tild}
\end{equation}%
Now, since $\cup _{C>C_{-}}\mathcal{F}_{C}\subset B_{\left( M,L_{\infty
}\right) }\left( s_{M},\tilde{R}_{n,D,\alpha }\right) $ where 
\begin{equation*}
B_{\left( M,L_{\infty }\right) }\left( s_{M},\tilde{R}_{n,D,\alpha }\right)
=\left\{ s\in M,\left\Vert s-s_{M}\right\Vert _{\infty }\leq \tilde{R}%
_{n,D,\alpha }\right\} \text{ },
\end{equation*}%
we have by (\ref{control_R_n_tild}), for all $s\in \cup _{C>C_{-}}\mathcal{F}%
_{C}$ and for all $n\geq n_{0}\left( A_{\infty },A_{cons},A_{+}\right) $,%
\begin{align*}
P\left( Ks-Ks_{M}\right) & =P\left[ \left( s-s_{M}\right) ^{2}\right] \\
& \leq \left\Vert s-s_{M}\right\Vert _{\infty }^{2} \\
& \leq \tilde{R}_{n,D,\alpha }^{2}\leq 1\text{.}
\end{align*}%
We thus have, for all $n\geq n_{0}\left( A_{\infty },A_{cons},A_{+}\right) $%
, 
\begin{equation*}
\bigcup\limits_{C>C_{-}}\mathcal{F}_{C}=\bigcup\limits_{C_{-}\wedge 1<C\leq
1}\mathcal{F}_{C}
\end{equation*}%
and by monotonicity of the collection $\mathcal{F}_{C}$, for some $q>1$ and $%
J=\left\lfloor \frac{\left\vert \ln \left( C_{-}\wedge 1\right) \right\vert 
}{\ln q}\right\rfloor +1$, it holds%
\begin{equation*}
\bigcup\limits_{C_{-}\wedge 1<C\leq 1}\mathcal{F}_{C}\subset
\bigcup\limits_{j=0}^{J}\mathcal{F}_{q^{j}C_{-}}.
\end{equation*}%
Simple computations show that, since $D\geq 1$ and $C_{-}\geq A_{l}\frac{D}{n%
}\geq \frac{A_{l}}{n}$, one can find a constant $L_{A_{l},q}$ such that%
\begin{equation*}
J\leq L_{A_{l},q}\ln n.
\end{equation*}%
Moreover, by monotonicity of $C\longmapsto \sup_{s\in \mathcal{F}%
_{C}}\left\vert \left( P_{n}-P\right) \left( \psi _{2}\circ \left(
s-s_{M}\right) \right) \right\vert $, we have uniformly in $C\in \left(
q^{j-1}C_{-},q^{j}C_{-}\right] ,$%
\begin{equation*}
\sup_{s\in \mathcal{F}_{C}}\left\vert \left( P_{n}-P\right) \left( \psi
_{2}\circ \left( s-s_{M}\right) \right) \right\vert \leq \sup_{s\in \mathcal{%
F}_{q^{j+1}C_{-}}}\left\vert \left( P_{n}-P\right) \left( \psi _{2}\circ
\left( s-s_{M}\right) \right) \right\vert \text{ }.
\end{equation*}%
Hence, taking the convention $\sup_{s\in \mathcal{\emptyset }}\left\vert
\left( P_{n}-P\right) \left( \psi _{2}\circ \left( s-s_{M}\right) \right)
\right\vert =0$, we get for all $n\geq n_{0}\left( A_{\infty
},A_{cons},A_{+}\right) $ and any $L>0$, 
\begin{align*}
& \mathbb{P}\left[ \forall C>C_{-},\text{ \ }\sup_{s\in \mathcal{F}%
_{C}}\left\vert \left( P_{n}-P\right) \left( \psi _{2}\circ \left(
s-s_{M}\right) \right) \right\vert \leq L\sqrt{\frac{CD}{n}}\tilde{R}%
_{n,D,\alpha }\right] \\
& \geq \mathbb{P}\left[ \forall j\in \left\{ 1,...,J\right\} ,\text{ \ }%
\sup_{s\in \mathcal{F}_{q^{j}C_{-}}}\left\vert \left( P_{n}-P\right) \left(
\psi _{2}\circ \left( s-s_{M}\right) \right) \right\vert \leq L\sqrt{\frac{%
q^{j}C_{-}D}{n}}\tilde{R}_{n,D,\alpha }\right] \text{ }.
\end{align*}%
Now, for any $L>0$,%
\begin{align}
& \mathbb{P}\left[ \forall j\in \left\{ 1,...,J\right\} ,\text{ \ }%
\sup_{s\in \mathcal{F}_{q^{j}C_{-}}}\left\vert \left( P_{n}-P\right) \left(
\psi _{2}\circ \left( s-s_{M}\right) \right) \right\vert \leq L\sqrt{\frac{%
q^{j}C_{-}D}{n}}\tilde{R}_{n,D,\alpha }\right]  \notag \\
& =1-\mathbb{P}\left[ \exists j\in \left\{ 1,...,J\right\} ,\text{ \ }%
\sup_{s\in \mathcal{F}_{q^{j}C_{-}}}\left\vert \left( P_{n}-P\right) \left(
\psi _{2}\circ \left( s-s_{M}\right) \right) \right\vert >L\sqrt{\frac{%
q^{j}C_{-}D}{n}}\tilde{R}_{n,D,\alpha }\right]  \notag \\
& \geq 1-\sum_{j=1}^{J}\mathbb{P}\left[ \sup_{s\in \mathcal{F}%
_{q^{j}C_{-}}}\left\vert \left( P_{n}-P\right) \left( \psi _{2}\circ \left(
s-s_{M}\right) \right) \right\vert >L\sqrt{\frac{q^{j}C_{-}D}{n}}\tilde{R}%
_{n,D,\alpha }\right] \text{ }.  \label{total_probability}
\end{align}%
Given $j\in \left\{ 1,...,J\right\} ,$ Lemma \ref{lower_rademacher copy(1)}
yields%
\begin{equation*}
\mathbb{E}\left[ \sup_{s\in \mathcal{F}_{q^{j}C_{-}}}\left\vert \left(
P_{n}-P\right) \left( \psi _{2}\circ \left( s-s_{M}\right) \right)
\right\vert \right] \leq 8\sqrt{\frac{q^{j}C_{-}D}{n}}\tilde{R}_{n,D,\alpha }%
\text{ },
\end{equation*}%
and next, we apply Bousquet's inequality (\ref{bousquet_2}) to handle the
deviations around the mean. We have%
\begin{align*}
& \sup_{s\in \mathcal{F}_{q^{j}C_{-}}}\left\Vert \psi _{2}\circ \left(
s-s_{M}\right) -P\left( \psi _{2}\circ \left( s-s_{M}\right) \right)
\right\Vert _{\infty } \\
& \leq 2\sup_{s\in \mathcal{F}_{q^{j}C_{-}}}\left\Vert \left( s-s_{M}\right)
^{2}\right\Vert _{\infty }\leq 2\tilde{R}_{n,D,\alpha }^{2}
\end{align*}%
and, for all $s\in \mathcal{F}_{q^{j}C_{-}}$,%
\begin{align*}
& 
\var%
\left( \psi _{2}\circ \left( s-s_{M}\right) \right) \\
& \leq P\left[ \left( s-s_{M}\right) ^{4}\right] \\
& \leq \left\Vert s-s_{M}\right\Vert _{\infty }^{2}P\left[ \left(
s-s_{M}\right) ^{2}\right] \\
& \leq \tilde{R}_{n,D,\alpha }^{2}q^{j}C_{-}\text{ }.
\end{align*}%
It follows that, for $\varepsilon =1$ and all $x>0,$%
\begin{equation}
\mathbb{P}\left[ \sup_{s\in \mathcal{F}_{q^{j}C_{-}}}\left\vert \left(
P_{n}-P\right) \left( \psi _{2}\circ \left( s-s_{M}\right) \right)
\right\vert \geq 16\sqrt{\frac{q^{j}C_{-}D}{n}}\tilde{R}_{n,D,\alpha }+\sqrt{%
\frac{2\tilde{R}_{n,D,\alpha }^{2}q^{j}C_{-}x}{n}}+\frac{8}{3}\frac{\tilde{R}%
_{n,D,\alpha }^{2}x}{n}\right] \leq \exp \left( -x\right) \text{ }.
\label{dev_2}
\end{equation}%
By consequence, as $D\geq A_{-}\left( \ln n\right) ^{2}$ and as $\tilde{R}%
_{n,D,\alpha }\leq 1$ for all $n\geq n_{0}\left( A_{\infty
},A_{cons},A_{+}\right) $, taking $x=\gamma \ln n$ in (\ref{dev_2}) for some 
$\gamma >0$, easy computations show that a positive constant $%
L_{A_{-},A_{l},\gamma }$ independent of $j$ exists such that for all $n\geq
n_{0}\left( A_{\infty },A_{cons},A_{+}\right) $,%
\begin{equation*}
\mathbb{P}\left[ \sup_{s\in \mathcal{F}_{q^{j}C_{-}}}\left\vert \left(
P_{n}-P\right) \left( \psi _{2}\circ \left( s-s_{M}\right) \right)
\right\vert \geq L_{A_{-},A_{l},\gamma }\sqrt{\frac{q^{j}C_{-}D}{n}}\tilde{R}%
_{n,D,\alpha }\right] \leq \frac{1}{n^{\gamma }}\text{ }.
\end{equation*}%
Hence, using (\ref{total_probability}), we get for all $n\geq n_{0}\left(
A_{\infty },A_{cons},A_{+}\right) $,%
\begin{align*}
& \mathbb{P}\left[ \forall C>C_{-},\text{ \ }\sup_{s\in \mathcal{F}%
_{C}}\left\vert \left( P_{n}-P\right) \left( \psi _{2}\circ \left(
s-s_{M}\right) \right) \right\vert \leq L_{A_{-},A_{l},\gamma }\sqrt{\frac{CD%
}{n}}\tilde{R}_{n,D,\alpha }\right] \\
& \geq 1-\frac{J}{n^{\gamma }}\text{ }.
\end{align*}%
And finally, as $J$ $\leq L_{A_{l},q}\ln n$, taking $\gamma =\beta +1$ and $%
q=2$ gives the result for all $n\geq n_{0}\left( A_{\infty
},A_{cons},A_{+},A_{l}\right) $.

\ \ \ \ \ \ \ \ \ \ \ \ \ \ \ \ \ \ \ \ \ \ \ \ \ \ \ \ \ \ \ \ \ \ \ \ \ \
\ \ \ \ \ \ \ \ \ \ \ \ \ \ \ \ \ \ \ \ \ \ \ \ \ \ \ \ \ \ \ \ \ \ \ \ \ \
\ \ \ \ \ \ \ \ \ \ \ \ \ \ \ \ \ \ \ \ \ \ \ \ \ \ \ \ \ \ \ \ \ \ \ \ \ \
\ \ \ \ \ \ \ \ \ \ \ \ \ \ \ \ \ \ \ \ \ \ \ \ \ \ \ \ \ \ \ \ 
$\blacksquare$%

\noindent Having controlled the residual empirical process driven by the
remainder terms in the expansion of the contrast, and having proved sharp
bounds for the expectation of the increments of the main empirical process
on the slices, it remains to combine the above lemmas in order to establish
the probability estimates controlling the empirical excess risk on the
slices.

\begin{lemma}
\label{majo_s1c_lemma copy(1)}Let $\beta ,A_{-},A_{+},A_{l},C>0$. Assume
that (\textbf{H1}), (\textbf{H2}), (\textbf{H3}) and (\ref{def_A_cons})
hold. A positive constant $A_{4}$ exists, only depending on $%
A,A_{3,M},\sigma _{\min },\beta $, such that, if 
\begin{equation*}
A_{l}\frac{D}{n}\leq C\leq \frac{1}{4}\left( 1+A_{4}\nu _{n}\right) ^{2}%
\frac{D}{n}\mathcal{K}_{1,M}^{2}\text{ \ \ \ and \ \ \ \ }A_{+}\frac{n}{%
\left( \ln n\right) ^{2}}\geq D\geq A_{-}\left( \ln n\right) ^{2}
\end{equation*}%
where $\nu _{n}=\max \left\{ \sqrt{\frac{\ln n}{D}},\text{ }\sqrt{\frac{D\ln
n}{n}},\text{ }R_{n,D,\alpha }\right\} $ is defined in (\ref{nu_n_proof}),
then for all $n\geq n_{0}\left( A_{\infty },A_{cons},A_{+},A_{l}\right) $,%
\begin{equation*}
\mathbb{P}\left[ \sup_{s\in \mathcal{F}_{C}}P_{n}\left( Ks_{M}-Ks\right)
\geq \left( 1+L_{A_{\infty },A,A_{3,M},\sigma _{\min },A_{-},A_{l},\beta
}\times \nu _{n}\right) \sqrt{\frac{CD}{n}}\mathcal{K}_{1,M}-C\right] \leq
2n^{-\beta }\text{ }.
\end{equation*}
\end{lemma}

\subsection*{}%
\textbf{Proof. }Start with%
\begin{align}
\sup_{s\in \mathcal{F}_{C}}P_{n}\left( Ks_{M}-Ks\right) & =\sup_{s\in 
\mathcal{F}_{C}}\left\{ P_{n}\left( \psi _{1,M}%
\cdot%
\left( s_{M}-s\right) -\psi _{2}\circ \left( s-s_{M}\right) \right) \right\}
\notag \\
& =\sup_{s\in \mathcal{F}_{C}}\left\{ \left( P_{n}-P\right) \left( \psi
_{1,M}%
\cdot%
\left( s_{M}-s\right) \right) -\left( P_{n}-P\right) \left( \psi _{2}\circ
\left( s-s_{M}\right) \right) -P\left( Ks-Ks_{M}\right) \right\}  \notag \\
& \leq \sup_{s\in \mathcal{F}_{C}}\left\{ \left( P_{n}-P\right) \left( \psi
_{1,M}%
\cdot%
\left( s_{M}-s\right) \right) -P\left( Ks-Ks_{M}\right) \right\}  \notag \\
& +\sup_{s\in \mathcal{F}_{C}}\left\vert \left( P_{n}-P\right) \left( \psi
_{2}\circ \left( s-s_{M}\right) \right) \right\vert .
\label{decomposition_sup_fc_bis}
\end{align}%
Next, recall that by definition,%
\begin{equation*}
D_{L}=\left\{ s\in B_{\left( M,L_{\infty }\right) }\left( s_{M},\tilde{R}%
_{n,D,\alpha }\right) ,\text{ }P\left( Ks-Ks_{M}\right) =L\right\} ,
\end{equation*}%
so we have%
\begin{align*}
& \sup_{s\in \mathcal{F}_{C}}\left\{ \left( P_{n}-P\right) \left( \psi _{1,M}%
\cdot%
\left( s_{M}-s\right) \right) -P\left( Ks-Ks_{M}\right) \right\} \\
& =\sup_{0\leq L\leq C}\sup_{s\in D_{L}}\left\{ \left( P_{n}-P\right) \left(
\psi _{1,M}%
\cdot%
\left( s_{M}-s\right) \right) -L\right\} \\
& \leq \sup_{0\leq L\leq C}\left\{ \sqrt{L}\sqrt{\sum_{k=1}^{D}\left(
P_{n}-P\right) ^{2}\left( \psi _{1,M}%
\cdot%
\varphi _{k}\right) }-L\right\}
\end{align*}%
where the last bound follows from Cauchy-Schwarz inequality. Hence, we
deduce from Lemma \ref{mino_sphere_lemma} that%
\begin{equation}
\mathbb{P}\left[ \sup_{s\in \mathcal{F}_{C}}\left\{ \left( P_{n}-P\right)
\left( \psi _{1,M}%
\cdot%
\left( s_{M}-s\right) \right) -P\left( Ks-Ks_{M}\right) \right\} \geq
\sup_{0\leq L\leq C}\left\{ \sqrt{L}\left( 1+\tau _{n}\right) \sqrt{\frac{D}{%
n}}\mathcal{K}_{1,M}-L\right\} \right] \leq n^{-\beta }\text{ ,}
\label{majo_s1c_1}
\end{equation}%
where 
\begin{align}
\tau _{n}& =L_{A,A_{3,M},\sigma _{\min },\beta }\left( \sqrt{\frac{\ln n}{D}}%
\vee \frac{\sqrt{\ln n}}{n^{1/4}}\right)  \notag \\
& \leq L_{A,A_{3,M},\sigma _{\min },\beta }\left( \sqrt{\frac{\ln n}{D}}\vee 
\sqrt{\frac{D\ln n}{n}}\right)  \notag \\
& \leq L_{A,A_{3,M},\sigma _{\min },\beta }\times \nu _{n}\text{ }.
\label{majo_nu_n}
\end{align}%
So, injecting (\ref{majo_nu_n}) in (\ref{majo_s1c_1}) we have%
\begin{equation*}
\mathbb{P}\left[ 
\begin{array}{c}
\sup_{s\in \mathcal{F}_{C}}\left\{ \left( P_{n}-P\right) \left( \psi _{1,M}%
\cdot%
\left( s_{M}-s\right) \right) -P\left( Ks-Ks_{M}\right) \right\} \\ 
\geq \sup_{0\leq L\leq C}\left\{ \sqrt{L}\left( 1+L_{A,A_{3,M},\sigma _{\min
},\beta }\times \nu _{n}\right) \sqrt{\frac{D}{n}}\mathcal{K}_{1,M}-L\right\}%
\end{array}%
\right] \leq n^{-\beta }
\end{equation*}%
and since we assume $C\leq \frac{1}{4}\left( 1+L_{A,A_{3,M},\sigma _{\min
},\beta }\times \nu _{n}\right) ^{2}\frac{D}{n}\mathcal{K}_{1,M}^{2}$ we see
that 
\begin{equation*}
\sup_{0\leq L\leq C}\left\{ \sqrt{L}\left( 1+L_{A,A_{3,M},\sigma _{\min
},\beta }\nu _{n}\right) \sqrt{\frac{D}{n}}\mathcal{K}_{1,M}-L\right\} =%
\sqrt{C}\left( 1+L_{A,A_{3,M},\sigma _{\min },\beta }\times \nu _{n}\right) 
\sqrt{\frac{D}{n}}\mathcal{K}_{1,M}-C
\end{equation*}%
and therefore%
\begin{equation}
\mathbb{P}\left[ \sup_{s\in \mathcal{F}_{C}}\left\{ \left( P_{n}-P\right)
\left( \psi _{1,M}%
\cdot%
\left( s_{M}-s\right) \right) -P\left( Ks-Ks_{M}\right) \right\} \geq \left(
1+L_{A,A_{3,M},\sigma _{\min },\beta }\nu _{n}\right) \sqrt{\frac{CD}{n}}%
\mathcal{K}_{1,M}-C\right] \leq n^{-\beta }\text{ }.  \label{majo_s1c}
\end{equation}%
Moreover, as $C\geq A_{l}\frac{D}{n}$, we derive from Lemma \ref%
{majo_uniform_s2c} that it holds, for all $n\geq n_{0}\left( A_{\infty
},A_{cons},A_{+},A_{l}\right) $,%
\begin{equation}
\mathbb{P}\left[ \sup_{s\in \mathcal{F}_{C}}\left\vert \left( P_{n}-P\right)
\left( \psi _{2}\circ \left( s-s_{M}\right) \right) \right\vert \geq
L_{A_{-},A_{l},\beta }\sqrt{\frac{CD}{n}}\tilde{R}_{n,D,\alpha }\right] \leq
n^{-\beta }\text{ }.  \label{majo_s2c}
\end{equation}%
Finally, noticing that%
\begin{align*}
\tilde{R}_{n,D,\alpha }& =\max \left\{ R_{n,D,\alpha },A_{\infty }\sqrt{%
\frac{D\ln n}{n}}\right\} \\
& \leq L_{A_{\infty },\sigma _{\min }}\max \left\{ R_{n,D,\alpha },\sqrt{%
\frac{D\ln n}{n}}\right\} \times \mathcal{K}_{1,M}\text{ \ \ \ \ \ \ by (\ref%
{lower_K1M_proof})} \\
& \leq L_{A_{\infty },\sigma _{\min }}\times \nu _{n}\times \mathcal{K}_{1,M}%
\text{ ,}
\end{align*}%
we deduce from (\ref{majo_s2c}) that, for all $n\geq n_{0}\left( A_{\infty
},A_{cons},A_{+},A_{l}\right) $,%
\begin{equation}
\mathbb{P}\left[ \sup_{s\in \mathcal{F}_{C}}\left\vert \left( P_{n}-P\right)
\left( \psi _{2}\circ \left( s-s_{M}\right) \right) \right\vert \geq
L_{A_{\infty },\sigma _{\min },A_{-},A_{l},\beta }\times \nu _{n}\sqrt{\frac{%
CD}{n}}\mathcal{K}_{1,M}\right] \leq n^{-\beta }  \label{majo_s2c_2}
\end{equation}%
and the conclusion follows by making use of (\ref{majo_s1c}) and (\ref%
{majo_s2c_2}) in inequality (\ref{decomposition_sup_fc_bis}).\ 
$\blacksquare$%

\noindent The second deviation bound for the empirical excess risk we need
to establish on the upper slice is proved in a similar way.

\begin{lemma}
\label{majo_s1>c_lemma}Let $\beta ,A_{-},A_{+},C\geq 0$. Assume that (%
\textbf{H1}), (\textbf{H2}), (\textbf{H3}) and (\ref{def_A_cons}) hold. A
positive constant $A_{5}$, depending on $A,A_{3,M},A_{\infty },\sigma _{\min
},A_{-}$ and $\beta $, exists such that, if it holds 
\begin{equation*}
C\geq \frac{1}{4}\left( 1+A_{5}\nu _{n}\right) ^{2}\frac{D}{n}\mathcal{K}%
_{1,M}^{2}\text{ \ \ \ and \ \ \ \ }A_{+}\frac{n}{\left( \ln n\right) ^{2}}%
\geq D\geq A_{-}\left( \ln n\right) ^{2}
\end{equation*}%
where $\nu _{n}=\max \left\{ \sqrt{\frac{\ln n}{D}},\text{ }\sqrt{\frac{D\ln
n}{n}},\text{ }R_{n,D,\alpha }\right\} $ is defined in (\ref{nu_n_proof}),
then for all $n\geq n_{0}\left( A_{\infty },A_{cons},A_{+}\right) $,%
\begin{equation*}
\mathbb{P}\left[ \sup_{s\in \mathcal{F}_{>C}}P_{n}\left( Ks_{M}-Ks\right)
\geq \left( 1+A_{5}\nu _{n}\right) \sqrt{\frac{CD}{n}}\mathcal{K}_{1,M}-C%
\right] \leq 2n^{-\beta }\text{ }.
\end{equation*}%
Moreover, when we only assume $C\geq 0$, we have for all $n\geq n_{0}\left(
A_{\infty },A_{cons},A_{+}\right) $,%
\begin{equation}
\mathbb{P}\left[ \sup_{s\in \mathcal{F}_{>C}}P_{n}\left( Ks_{M}-Ks\right)
\geq \frac{1}{4}\left( 1+A_{5}\nu _{n}\right) ^{2}\frac{D}{n}\mathcal{K}%
_{1,M}^{2}\right] \leq 2n^{-\beta }\text{ }.  \label{majo_S>C_4}
\end{equation}
\end{lemma}

\subsection*{}%
\textbf{Proof. }First observe that%
\begin{align}
\sup_{s\in \mathcal{F}_{>C}}P_{n}\left( Ks_{M}-Ks\right) & =\sup_{s\in 
\mathcal{F}_{>C}}\left\{ P_{n}\left( \psi _{1,M}%
\cdot%
\left( s_{M}-s\right) -\psi _{2}\circ \left( s-s_{M}\right) \right) \right\}
\notag \\
& =\sup_{s\in \mathcal{F}_{>C}}\left\{ \left( P_{n}-P\right) \left( \psi
_{1,M}%
\cdot%
\left( s_{M}-s\right) \right) -\left( P_{n}-P\right) \left( \psi _{2}\circ
\left( s-s_{M}\right) \right) -P\left( Ks-Ks_{M}\right) \right\}  \notag \\
& =\sup_{s\in \mathcal{F}_{>C}}\left\{ \left( P_{n}-P\right) \left( \psi
_{1,M}%
\cdot%
\left( s_{M}-s\right) \right) -P\left( Ks-Ks_{M}\right) -\left(
P_{n}-P\right) \left( \psi _{2}\circ \left( s-s_{M}\right) \right) \right\} 
\notag \\
& =\sup_{L>C}\sup_{s\in D_{L}}\left\{ \left( P_{n}-P\right) \left( \psi
_{1,M}%
\cdot%
\left( s_{M}-s\right) \right) -L-\left( P_{n}-P\right) \left( \psi _{2}\circ
\left( s-s_{M}\right) \right) \right\}  \notag \\
& \leq \sup_{L>C}\left\{ \sqrt{L}\sqrt{\sum_{k=1}^{D}\left( P_{n}-P\right)
^{2}\left( \psi _{1,M}%
\cdot%
\varphi _{k}\right) }-L+\sup_{s\in \mathcal{F}_{L}}\left\vert \left(
P_{n}-P\right) \left( \psi _{2}\circ \left( s-s_{M}\right) \right)
\right\vert \right\}  \label{control_s>C}
\end{align}%
where the last bound follows from Cauchy-Schwarz inequality. Now, the end of
the proof is similar to that of Lemma \ref{majo_s1c_lemma copy(1)} and
follows from the same kind of computations. Indeed, from Lemma \ref%
{mino_sphere_lemma} we deduce that%
\begin{equation}
\mathbb{P}\left[ \sqrt{\sum_{k=1}^{D}\left( P_{n}-P\right) ^{2}\left( \psi
_{1,M}%
\cdot%
\varphi _{k}\right) }\geq \left( 1+L_{A,A_{3,M},\sigma _{\min },\beta
}\times \nu _{n}\right) \sqrt{\frac{D}{n}}\mathcal{K}_{1,M}\right] \leq
n^{-\beta }  \label{majo_s1c_2}
\end{equation}%
and, since 
\begin{equation*}
C\geq \frac{1}{4}\frac{D}{n}\mathcal{K}_{1,M}^{2}\geq \sigma _{\min }^{2}%
\frac{D}{n}\text{ ,}
\end{equation*}%
we apply Lemma \ref{majo_uniform_s2c} with $A_{l}=\sigma _{\min }^{2}$, and
deduce that, for all $n\geq n_{0}\left( A_{\infty },A_{cons},A_{+}\right) $,%
\begin{equation}
\mathbb{P}\left[ \forall L>C,\text{ \ }\sup_{s\in \mathcal{F}_{L}}\left\vert
\left( P_{n}-P\right) \left( \psi _{2,M}^{s}%
\cdot%
\left( s-s_{M}\right) \right) \right\vert \geq L_{A_{\infty },\sigma _{\min
},A_{-},\beta }\times \nu _{n}\sqrt{\frac{LD}{n}}\mathcal{K}_{1,M}\right]
\leq n^{-\beta }\text{ .}  \label{majo_unif_proof}
\end{equation}%
Now using (\ref{majo_s1c_2}) and (\ref{majo_unif_proof}) in (\ref%
{control_s>C}) we obtain, for all $n\geq n_{0}\left( A_{\infty
},A_{cons},A_{+}\right) $,%
\begin{equation}
\mathbb{P}\left[ \sup_{s\in \mathcal{F}_{>C}}P_{n}\left( Ks_{M}-Ks\right)
\geq \sup_{L>C}\left\{ \left( 1+L_{A,A_{3,M},A_{\infty },\sigma _{\min
},A_{-},\beta }\times \nu _{n}\right) \sqrt{\frac{LD}{n}}\mathcal{K}%
_{1,M}-L\right\} \right] \leq 2n^{-\beta }  \label{def_A_4}
\end{equation}%
and we set $A_{5}=L_{A,A_{3,M},A_{\infty },\sigma _{\min },A_{-},\beta }$
where $L_{A,A_{3,M},A_{\infty },\sigma _{\min },A_{-},\beta }$ is the
constant in (\ref{def_A_4}). For $C\geq \frac{1}{4}\left( 1+A_{5}\nu
_{n}\right) ^{2}\frac{D}{n}\mathcal{K}_{1,M}^{2}$ we get 
\begin{equation*}
\sup_{L>C}\left\{ \sqrt{L}\left( 1+A_{5}\nu _{n}\right) \sqrt{\frac{D}{n}}%
\mathcal{K}_{1,M}-L\right\} =\left( 1+A_{5}\nu _{n}\right) \sqrt{\frac{CD}{n}%
}\mathcal{K}_{1,M}-C
\end{equation*}%
and by consequence,%
\begin{equation*}
\mathbb{P}\left[ \sup_{s\in \mathcal{F}_{>C}}P_{n}\left( Ks_{M}-Ks\right)
\geq \left( 1+A_{5}\nu _{n}\right) \sqrt{\frac{CD}{n}}\mathcal{K}_{1,M}-C%
\right] \leq 2n^{-\beta }\text{ ,}
\end{equation*}%
which gives the first part of the lemma. The second part comes from (\ref%
{def_A_4}) and the fact that, for any value of $C\geq 0$,%
\begin{equation*}
\sup_{L>C}\left\{ \sqrt{L}\left( 1+A_{5}\nu _{n}\right) \sqrt{\frac{D}{n}}%
\mathcal{K}_{1,M}-L\right\} \leq \left( 1+A_{5}\nu _{n}\right) ^{2}\frac{D}{%
4n}\mathcal{K}_{1,M}^{2}\text{ .}
\end{equation*}

\ \ \ \ \ \ \ \ \ \ \ \ \ \ \ \ \ \ \ \ \ \ \ \ \ \ \ \ \ \ \ \ \ \ \ \ \ \
\ \ \ \ \ \ \ \ \ \ \ \ \ \ \ \ \ \ \ \ \ \ \ \ \ \ \ \ \ \ \ \ \ \ \ \ \ \
\ \ \ \ \ \ \ \ \ \ \ \ \ \ \ \ \ \ \ \ \ \ \ \ \ \ \ \ \ \ \ \ \ \ \ \ \ \
\ \ \ \ \ \ \ \ \ \ \ \ \ \ \ \ \ \ \ \ \ \ \ \ \ \ \ \ \ \ \ \ 
$\blacksquare$%

\begin{lemma}
\label{mino_S_r_C}Let $r>1$ and $C,\beta >0$. Assume that (\textbf{H1}), (%
\textbf{H2}), (\textbf{H4}) and (\ref{def_A_cons}) hold and let $\varphi
=\left( \varphi _{k}\right) _{k=1}^{D}$ be an orthonormal basis of $\left(
M,\left\Vert 
\cdot%
\right\Vert _{2}\right) $ satisfying (\textbf{H4}). If positive constants $%
A_{-},A_{+},A_{l},A_{u}$ exist such that 
\begin{equation*}
A_{+}\frac{n}{\left( \ln n\right) ^{2}}\geq D\geq A_{-}\left( \ln n\right)
^{2}\text{ \ \ and \ \ }A_{l}\frac{D}{n}\leq rC\leq A_{u}\frac{D}{n}\text{ },
\end{equation*}%
and if the constant $A_{\infty }$ defined in (\ref{def_R_n_D_tild})
satisfies 
\begin{equation*}
A_{\infty }\geq 64B_{2}A\sqrt{2A_{u}}\sigma _{\min }^{-1}r_{M}\left( \varphi
\right) ,
\end{equation*}%
then a positive constant $L_{A_{-},A_{l},A_{u},A,A_{\infty },\sigma _{\min
},r_{M}\left( \varphi \right) ,\beta }$ exists such that,

\noindent for all $n\geq n_{0}\left( A_{-},A_{+},A_{u},A_{l},A,A_{\infty
},A_{cons},B_{2},r_{M}\left( \varphi \right) ,\sigma _{\min }\right) $,%
\begin{equation*}
\mathbb{P}\left( \sup_{s\in \mathcal{F}_{\left( C,rC\right] }}P_{n}\left(
Ks_{M}-Ks\right) \leq \left( 1-L_{A_{-},A_{l},A_{u},A,A_{\infty },\sigma
_{\min },r_{M}\left( \varphi \right) ,\beta }\times \nu _{n}\right) \sqrt{%
\frac{rCD}{n}}\mathcal{K}_{1,M}-rC\right) \leq 2n^{-\beta }\text{ },
\end{equation*}%
where $\nu _{n}=\max \left\{ \sqrt{\frac{\ln n}{D}},\text{ }\sqrt{\frac{D\ln
n}{n}},\text{ }R_{n,D,\alpha }\right\} $ is defined in (\ref{nu_n_proof}).
\end{lemma}

\subsection*{}%
\textbf{Proof. }Start with 
\begin{align}
& \sup_{s\in \mathcal{F}_{\left( C,rC\right] }}P_{n}\left( Ks_{M}-Ks\right) 
\notag \\
& =\sup_{s\in \mathcal{F}_{\left( C,rC\right] }}\left\{ \left(
P_{n}-P\right) \left( Ks_{M}-Ks\right) +P\left( Ks_{M}-Ks\right) \right\} 
\notag \\
& \geq \sup_{s\in \mathcal{F}_{\left( C,rC\right] }}\left( P_{n}-P\right)
\left( \psi _{1,M}%
\cdot%
\left( s_{M}-s\right) \right) -\sup_{s\in \mathcal{F}_{\left( C,rC\right]
}}\left( P_{n}-P\right) \left( \psi _{2}\circ \left( s-s_{M}\right) \right)
-\sup_{s\in \mathcal{F}_{\left( C,rC\right] }}P\left( Ks-Ks_{M}\right) 
\notag \\
& \geq \sup_{s\in \mathcal{F}_{\left( C,rC\right] }}\left( P_{n}-P\right)
\left( \psi _{1,M}%
\cdot%
\left( s_{M}-s\right) \right) -\sup_{s\in \mathcal{F}_{rC}}\left(
P_{n}-P\right) \left( \psi _{2}\circ \left( s-s_{M}\right) \right) -rC
\label{mino_skc}
\end{align}%
and set%
\begin{align*}
S_{1,r,C}& =\sup_{s\in \mathcal{F}_{\left( C,rC\right] }}\left(
P_{n}-P\right) \left( \psi _{1,M}%
\cdot%
\left( s_{M}-s\right) \right) \\
M_{1,r,C}& =\mathbb{E}\left[ \sup_{s\in \mathcal{F}_{\left( C,rC\right]
}}\left( P_{n}-P\right) \left( \psi _{1,M}%
\cdot%
\left( s_{M}-s\right) \right) \right] \\
b_{1,r,C}& =\sup_{s\in \mathcal{F}_{\left( C,rC\right] }}\left\Vert \psi
_{1,M}%
\cdot%
\left( s_{M}-s\right) -P\left( \psi _{1,M}%
\cdot%
\left( s_{M}-s\right) \right) \right\Vert _{\infty } \\
\sigma _{1,r,C}^{2}& =\sup_{s\in \mathcal{F}_{\left( C,rC\right] }}%
\var%
\left( \psi _{1,M}%
\cdot%
\left( s_{M}-s\right) \right) .
\end{align*}%
By Klein-Rio's Inequality (\ref{klein_rio_2}), we get, for all $\delta ,x>0$,%
\begin{equation}
\mathbb{P}\left( S_{1,r,C}\leq \left( 1-\delta \right) M_{1,r,C}-\sqrt{\frac{%
2\sigma _{1,r,C}^{2}x}{n}}-\left( 1+\frac{1}{\delta }\right) \frac{b_{1,r,C}x%
}{n}\right) \leq \exp \left( -x\right) \text{ }.  \label{dev_s1kc}
\end{equation}%
Then, notice that all conditions of Lemma \ref{bound_moment_order_1} are
satisfied, and that it gives by (\ref{M1KC_mino}),

\noindent for all $n\geq n_{0}\left(
A_{-},A_{+},A_{u},A_{l},A,B_{2},r_{M}\left( \varphi \right) ,\sigma _{\min
}\right) $,%
\begin{equation}
M_{1,r,C}\geq \left( 1-\frac{L_{A,A_{l},A_{u},\sigma _{\min }}}{\sqrt{D}}%
\right) \sqrt{\frac{rCD}{n}}\mathcal{K}_{1,M}\text{ }.
\label{mino_mean_proof}
\end{equation}%
In addition, observe that%
\begin{equation}
\sigma _{1,r,C}^{2}\leq \sup_{s\in \mathcal{F}_{\left( C,rC\right] }}P\left(
\psi _{1,M}^{2}%
\cdot%
\left( s_{M}-s\right) ^{2}\right) \leq 16A^{2}rC\text{ \ \ by (\ref%
{control_phi_1M_proof})}  \label{sigma1rC}
\end{equation}%
and%
\begin{equation}
b_{1,r,C}=\sup_{s\in \mathcal{F}_{\left( C,rC\right] }}\left\Vert \psi _{1,M}%
\cdot%
\left( s_{M}-s\right) \right\Vert _{\infty }\leq 4Ar_{M}\left( \varphi
\right) \sqrt{rCD}\text{ \ \ \ by (\ref{control_phi_1M_proof}) and (\textbf{%
H4})}  \label{b1rC}
\end{equation}%
Hence, using (\ref{mino_mean_proof}), (\ref{sigma1rC}) and (\ref{b1rC}) in
inequality (\ref{dev_s1kc}), we get for all $x>0$ and

\noindent for all $n\geq n_{0}\left(
A_{-},A_{+},A_{u},A_{l},A,B_{2},r_{M}\left( \varphi \right) ,\sigma _{\min
}\right) $,%
\begin{align*}
& \mathbb{P}\left( S_{1,r,C}\leq \left( 1-\delta \right) \left( 1-\frac{%
L_{A,A_{l},A_{u},\sigma _{\min }}}{\sqrt{D}}\right) \sqrt{\frac{rCD}{n}}%
\mathcal{K}_{1,M}-\sqrt{\frac{32A^{2}rCx}{n}}-\left( 1+\frac{1}{\delta }%
\right) \frac{4Ar_{M}\left( \varphi \right) \sqrt{rCD}x}{n}\right) \\
& \leq \exp \left( -x\right) \text{ }.
\end{align*}%
Now, taking $x=\beta \ln n$, $\delta =\frac{\sqrt{\ln n}}{n^{1/4}}$ and
using (\ref{lower_K1M_proof}), we deduce by simple computations that for all 
$n\geq n_{0}\left( A_{-},A_{+},A_{u},A_{l},A,B_{2},r_{M}\left( \varphi
\right) ,\sigma _{\min }\right) $,%
\begin{equation}
\mathbb{P}\left( S_{1,r,C}\leq \left( 1-L_{A,A_{l},A_{u},\sigma _{\min
},r_{M}\left( \varphi \right) ,\beta }\times \left( \sqrt{\frac{\ln n}{D}}%
\vee \frac{\sqrt{\ln n}}{n^{1/4}}\right) \right) \sqrt{\frac{rCD}{n}}%
\mathcal{K}_{1,M}\right) \leq n^{-\beta }  \label{mino_s1kc}
\end{equation}%
and as%
\begin{equation*}
\sqrt{\frac{\ln n}{D}}\vee \frac{\sqrt{\ln n}}{n^{1/4}}\leq \sqrt{\frac{\ln n%
}{D}}\vee \sqrt{\frac{D\ln n}{n}}\leq \nu _{n}
\end{equation*}%
(\ref{mino_s1kc}) gives, for all $n\geq n_{0}\left(
A_{-},A_{+},A_{u},A_{l},A,B_{2},r_{M}\left( \varphi \right) ,\sigma _{\min
}\right) $,%
\begin{equation}
\mathbb{P}\left( S_{1,r,C}\leq \left( 1-L_{A,A_{l},A_{u},\sigma _{\min
},r_{M}\left( \varphi \right) ,\beta }\times \nu _{n}\right) \sqrt{\frac{rCD%
}{n}}\mathcal{K}_{1,M}\right) \leq n^{-\beta }\text{ }.  \label{mino_s1kc_2}
\end{equation}%
Moreover, from Lemma \ref{majo_uniform_s2c} we deduce that, for all $n\geq
n_{0}\left( A_{\infty },A_{cons},A_{+},A_{l}\right) $,%
\begin{equation}
\mathbb{P}\left[ \sup_{s\in \mathcal{F}_{rC}}\left\vert \left(
P_{n}-P\right) \left( \psi _{2}\circ \left( s-s_{M}\right) \right)
\right\vert \geq L_{A_{-},A_{l},\beta }\sqrt{\frac{rCD}{n}}\tilde{R}%
_{n,D,\alpha }\right] \leq n^{-\beta }  \label{S2C}
\end{equation}%
and noticing that%
\begin{align*}
\tilde{R}_{n,D,\alpha }& =\max \left\{ R_{n,D,\alpha }\text{ };\text{ }%
A_{\infty }\sqrt{\frac{D\ln n}{n}}\right\} \\
& \leq L_{A_{\infty },\sigma _{\min }}\max \left\{ R_{n,D,\alpha }\text{ };%
\text{ }\sqrt{\frac{D\ln n}{n}}\right\} \times \mathcal{K}_{1,M}\text{ \ \ \
\ \ \ by (\ref{lower_K1M_proof})} \\
& \leq L_{A_{\infty },\sigma _{\min }}\times \nu _{n}\times \mathcal{K}_{1,M}%
\text{ ,}
\end{align*}%
we deduce from (\ref{S2C}) that for all $n\geq n_{0}\left( A_{\infty
},A_{cons},A_{+},A_{l}\right) $,%
\begin{equation}
\mathbb{P}\left[ \sup_{s\in \mathcal{F}_{rC}}\left\vert \left(
P_{n}-P\right) \left( \psi _{2}\circ \left( s-s_{M}\right) \right)
\right\vert \geq L_{A_{-},A_{l},A_{\infty },\sigma _{\min },\beta }\times
\nu _{n}\times \sqrt{\frac{rCD}{n}}\mathcal{K}_{1,M}\right] \leq n^{-\beta }%
\text{ }.  \label{S2C_2}
\end{equation}%
Finally, using (\ref{mino_s1kc_2}) and (\ref{S2C_2}) in (\ref{mino_skc}) we
get that,

\noindent for all $n\geq n_{0}\left( A_{-},A_{+},A_{u},A_{l},A,A_{\infty
},A_{cons},B_{2},r_{M}\left( \varphi \right) ,\sigma _{\min }\right) $,%
\begin{equation*}
\mathbb{P}\left( \sup_{s\in \mathcal{F}_{\left( C,rC\right] }}P_{n}\left(
Ks_{M}-Ks\right) \leq \left( 1-L_{A_{-},A_{l},A_{u},A,A_{\infty },\sigma
_{\min },r_{M}\left( \varphi \right) ,\beta }\times \nu _{n}\right) \sqrt{%
\frac{rCD}{n}}\mathcal{K}_{1,M}-rC\right) \leq 2n^{-\beta }\text{ },
\end{equation*}%
which concludes the proof. 
$\blacksquare$%

\subsection{Probabilistic Tools\label{section_probabilistic_tools}}

We recall here the main probabilistic results that are instrumental in our
proofs.

\noindent Let us begin with the $L_{p}$-version of Hoffmann-J%
\o{}%
rgensen's inequality, that can be found for example in \cite{LedouxTal:91},
Proposition 6.10, p.157.

\begin{theorem}
\label{hoffmann_jorgensen}For any independent mean zero random variables $%
Y_{j},$ $j=1,...,n$ taking values in a Banach space $\left( \mathcal{B}%
,\left\Vert .\right\Vert \right) $ and satisfying $\mathbb{E}\left[
\left\Vert Y_{j}\right\Vert ^{p}\right] <+\infty$ for some $p\geq1,$ we have%
\begin{equation*}
\mathbb{E}^{1/p}\left\Vert \sum_{j=1}^{n}Y_{j}\right\Vert ^{p}\leq
B_{p}\left( \mathbb{E}\left\Vert \sum_{j=1}^{n}Y_{j}\right\Vert +\mathbb{E}%
^{1/p}\left( \max_{1\leq j\leq n}\left\Vert Y_{j}\right\Vert \right)
^{p}\right)
\end{equation*}
where $B_{p}$ is a universal constant depending only on $p$.
\end{theorem}

\noindent We will use this theorem for $p=2$ in order to control suprema of
empirical processes. In order to be more specific, let $\mathcal{F}$ be a
class of measurable functions from a measurable space $\mathcal{Z}$ to $%
\mathbb{R}$ and $\left( X_{1},...,X_{n}\right) $ be independent variables of
common law $P$ taking values in $\mathcal{Z}$. We then denote by $\mathcal{B}%
=l^{\infty}\left( \mathcal{F}\right) $ the space of uniformly bounded
functions on $\mathcal{F}$ and, for any $b\in\mathcal{B}$, we set $%
\left\Vert b\right\Vert =\sup_{f\in\mathcal{F}}\left\vert b\left( f\right)
\right\vert $. Thus $\left( \mathcal{B},\left\Vert .\right\Vert \right) $ is
a Banach space. Indeed we shall apply Theorem \ref{hoffmann_jorgensen} to
the independent random variables, with mean zero and taking values in $%
\mathcal{B}$, defined by 
\begin{equation*}
Y_{j}=\left\{ f\left( X_{j}\right) -Pf,\text{ }f\in\mathcal{F}\right\} \text{
}.
\end{equation*}
More precisely, we will use the following result, which is a straightforward
application of Theorem \ref{hoffmann_jorgensen}. Denote by 
\begin{equation*}
P_{n}=\frac{1}{n}\sum_{i=1}^{n}\delta_{X_{i}}
\end{equation*}
the empirical measure associated to the sample $\left(
X_{1},...,X_{n}\right) $ and by 
\begin{equation*}
\left\Vert P_{n}-P\right\Vert _{\mathcal{F}}=\sup_{f\in\mathcal{F}%
}\left\vert \left( P_{n}-P\right) \left( f\right) \right\vert
\end{equation*}
the supremum of the empirical process over $\mathcal{F}$.

\begin{corollary}
If $\mathcal{F}$ is a class of measurable functions from a measurable space $%
\mathcal{Z}$ to $\mathbb{R}$ satisfying 
\begin{equation*}
\sup_{z\in\mathcal{Z}}\sup_{f\in\mathcal{F}}\left\vert f\left( z\right)
-Pf\right\vert =\sup_{f\in\mathcal{F}}\left\Vert f-Pf\right\Vert _{\infty
}<+\infty
\end{equation*}
and $\left( X_{1},...,X_{n}\right) $ are $n$ i.i.d. random variables taking
values in $\mathcal{Z}$, then an absolute constant $B_{2}$ exists such that,%
\begin{equation}
\mathbb{E}^{1/2}\left[ \left\Vert P_{n}-P\right\Vert _{\mathcal{F}}^{2}%
\right] \leq B_{2}\left( \mathbb{E}\left[ \left\Vert P_{n}-P\right\Vert _{%
\mathcal{F}}\right] +\frac{\sup_{f\in\mathcal{F}}\left\Vert f-Pf\right\Vert
_{\infty}}{n}\right) \text{ }.  \label{emp_hoff}
\end{equation}
\end{corollary}

\noindent Another tool we need is a comparison theorem for Rademacher
processes, see Theorem 4.12 of \cite{LedouxTal:91}. A function $\varphi :%
\mathbb{R\rightarrow R}$ is called a contraction if $\left\vert
\varphi\left( u\right) -\varphi\left( v\right) \right\vert \leq\left\vert
u-v\right\vert $ for all $u,v\in\mathbb{R}$. Moreover, for a subset $T$ $%
\subset$ $\mathbb{R}^{n}$ we set 
\begin{equation*}
\left\Vert h\left( t\right) \right\Vert _{T}=\left\Vert h\right\Vert
_{T}=\sup_{t\in T}\left\vert h\left( t\right) \right\vert \text{ }.
\end{equation*}

\begin{theorem}
\label{comparison_theorem}Let $\left( \varepsilon_{1},...,\varepsilon
_{n}\right) $ be $n$ i.i.d. Rademacher variables and $F:\mathbb{R}%
_{+}\longrightarrow\mathbb{R}_{+}$ be a convex and increasing function.
Furthermore, let $\varphi_{i}:\mathbb{R\longrightarrow R},$ $i\leq n,$ be
contractions such that $\varphi_{i}\left( 0\right) =0$. Then, for any
bounded subset $T$ $\subset$ $\mathbb{R}^{n},$%
\begin{equation*}
\mathbb{E}F\left( \left\Vert \sum_{i}\varepsilon_{i}\varphi_{i}\left(
t_{i}\right) \right\Vert _{T}\right) \leq2\mathbb{E}F\left( \left\Vert
\sum_{i}\varepsilon_{i}t_{i}\right\Vert _{T}\right) .
\end{equation*}
\end{theorem}

\noindent The next tool is the well known Bernstein's inequality, that can
be found for example in \cite{Massart:07}, Proposition 2.9.

\begin{theorem}
\label{bernstein}(Bernstein's inequality) Let $\left( X_{1},...,X_{n}\right) 
$ be independent real valued random variables and define%
\begin{equation*}
S=\frac{1}{n}\sum_{i=1}^{n}\left( X_{i}-\mathbb{E}\left[ X_{i}\right]
\right) .
\end{equation*}%
Assuming that%
\begin{equation*}
v=\frac{1}{n}\sum_{i=1}^{n}\mathbb{E}\left[ X_{i}^{2}\right] <\infty
\end{equation*}%
and%
\begin{equation*}
\left\vert X_{i}\right\vert \leq b\text{ \ }a.s.
\end{equation*}%
we have, for every $x>0$,%
\begin{equation}
\mathbb{P}\left[ \left\vert S\right\vert \geq \sqrt{2v\frac{x}{n}}+\frac{bx}{%
3n}\right] \leq 2\exp \left( -x\right) .  \label{bernstein_ineq_reg}
\end{equation}
\end{theorem}

\noindent We turn now to concentration inequalities for the empirical
process around its mean. Bousquet's inequality \cite{Bousquet:02} provides
optimal constants for the deviations at the right. Klein-Rio's inequality 
\cite{Klein_Rio:05} gives sharp constants for the deviations at the left,
that slightly improves Klein's inequality \cite{Klein:02}.

\begin{theorem}
\label{theorem_concentrations}Let $\left( \xi _{1},...,\xi _{n}\right) $ be $%
n$ i.i.d. random variables having common law $P$ and taking values in a
measurable space $\mathcal{Z}$. If $\mathcal{F}$ is a class of measurable
functions from $\mathcal{Z}$ to $\mathbb{R}$ satisfying%
\begin{equation*}
\text{\ }\left\vert f\left( \xi _{i}\right) -Pf\right\vert \leq b\text{ \ \ }%
a.s.,\text{ for all }f\in \mathcal{F},\text{ }i\leq n,
\end{equation*}%
then, by setting 
\begin{equation*}
\sigma _{\mathcal{F}}^{2}=\sup_{f\in \mathcal{F}}\left\{ P\left(
f^{2}\right) -\left( Pf\right) ^{2}\right\} ,
\end{equation*}%
we have,\ for all $x\geq 0$,

\textbf{Bousquet's inequality }:%
\begin{equation}
\mathbb{P}\left[ \left\Vert P_{n}-P\right\Vert _{\mathcal{F}}-\mathbb{E}%
\left[ \left\Vert P_{n}-P\right\Vert _{\mathcal{F}}\right] \geq \sqrt{%
2\left( \sigma_{\mathcal{F}}^{2}+2b\mathbb{E}\left[ \left\Vert
P_{n}-P\right\Vert _{\mathcal{F}}\right] \right) \frac{x}{n}}+\frac{bx}{3n}%
\right] \leq\exp\left( -x\right)  \label{bousquet}
\end{equation}
and we can deduce that, for all $\varepsilon,x>0$, it holds%
\begin{equation}
\mathbb{P}\left[ \left\Vert P_{n}-P\right\Vert _{\mathcal{F}}-\mathbb{E}%
\left[ \left\Vert P_{n}-P\right\Vert _{\mathcal{F}}\right] \geq\sqrt {%
2\sigma_{\mathcal{F}}^{2}\frac{x}{n}}+\varepsilon\mathbb{E}\left[ \left\Vert
P_{n}-P\right\Vert _{\mathcal{F}}\right] +\left( \frac{1}{\varepsilon}+\frac{%
1}{3}\right) \frac{bx}{n}\right] \leq\exp\left( -x\right) .
\label{bousquet_2}
\end{equation}

\textbf{Klein-Rio's inequality :}%
\begin{equation}
\mathbb{P}\left[ \mathbb{E}\left[ \left\Vert P_{n}-P\right\Vert _{\mathcal{F}%
}\right] -\left\Vert P_{n}-P\right\Vert _{\mathcal{F}}\geq \sqrt{2\left(
\sigma_{\mathcal{F}}^{2}+2b\mathbb{E}\left[ \left\Vert P_{n}-P\right\Vert _{%
\mathcal{F}}\right] \right) \frac{x}{n}}+\frac{bx}{n}\right] \leq\exp\left(
-x\right)  \label{klein_rio}
\end{equation}
and again, we can deduce that, for all $\varepsilon,x>0$, it holds%
\begin{equation}
\mathbb{P}\left[ \mathbb{E}\left[ \left\Vert P_{n}-P\right\Vert _{\mathcal{F}%
}\right] -\left\Vert P_{n}-P\right\Vert _{\mathcal{F}}\geq \sqrt{2\sigma_{%
\mathcal{F}}^{2}\frac{x}{n}}+\varepsilon\mathbb{E}\left[ \left\Vert
P_{n}-P\right\Vert _{\mathcal{F}}\right] +\left( \frac {1}{\varepsilon}%
+1\right) \frac{bx}{n}\right] \leq\exp\left( -x\right) .  \label{klein_rio_2}
\end{equation}
\end{theorem}

\noindent The following result is due to Ledoux \cite{Ledoux:96}. We will
use it along the proofs through Corollary \ref{cor_alt_hoff_2} which is
stated below. From now on, we set for short $Z=\left\Vert P_{n}-P\right\Vert
_{\mathcal{F}}$.

\begin{theorem}
\label{Theorem_Ledoux}Let $\left( \xi _{1},...,\xi _{n}\right) $ be
independent random with values in some measurable space $\left( \mathcal{Z},%
\mathcal{T}\right) $ and $\mathcal{F}$ be some countable class of
real-valued measurable functions from $\mathcal{Z}$. Let $\left( \xi
_{1}^{\prime },...,\xi _{n}^{\prime }\right) $ be independent from $\left(
\xi _{1},...,\xi _{n}\right) $ and with the same distribution. Setting%
\begin{equation*}
v=\mathbb{E}\left[ \sup_{f\in \mathcal{F}}\frac{1}{n}\sum_{i=1}^{n}\left(
f\left( \xi _{i}\right) -f\left( \xi _{i}^{\prime }\right) \right) ^{2}%
\right]
\end{equation*}%
then%
\begin{equation*}
\mathbb{E}\left[ Z^{2}\right] -\mathbb{E}\left[ Z\right] ^{2}\leq \frac{v}{n}%
\text{ .}
\end{equation*}
\end{theorem}

\begin{corollary}
\label{cor_alt_hoff_2}Under notations of Theorem \ref{theorem_concentrations}%
, if some $\varkappa _{n}\in \left( 0,1\right) $ exists such that%
\begin{equation*}
\varkappa _{n}^{2}\mathbb{E}\left[ Z^{2}\right] \geq \frac{\sigma ^{2}}{n}
\end{equation*}%
and%
\begin{equation*}
\varkappa _{n}^{2}\sqrt{\mathbb{E}\left[ Z^{2}\right] }\geq \frac{b}{n}
\end{equation*}%
then we have, for a numerical constant $A_{1,-}$, 
\begin{equation*}
\left( 1-\varkappa _{n}A_{1,-}\right) \sqrt{\mathbb{E}\left[ Z^{2}\right] }%
\leq \mathbb{E}\left[ Z\right] \text{ }.
\end{equation*}
\end{corollary}

\subsection*{}%
\textbf{Proof of Corollary \ref{cor_alt_hoff_2}.} Just use Theorem \ref%
{Theorem_Ledoux}, noticing the fact that%
\begin{equation*}
\sqrt{\mathbb{E}\left[ Z^{2}\right] }-\mathbb{E}\left[ Z\right] \leq \sqrt{%
\mathbb{V}\left( Z\right) }
\end{equation*}%
and that, with notations of Theorem \ref{Theorem_Ledoux},%
\begin{equation*}
v\leq 2\sigma ^{2}+32b\mathbb{E}\left[ Z\right] \text{ .}
\end{equation*}%
The result then follows from straightforward calculations.\ 
$\blacksquare$%

\pagebreak

\noindent {\LARGE Acknowledgements}\bigskip

The author thanks gratefully the editor David Ruppert, the associate editor
and the two anonymous referees for their comments and suggestions, that
greatly improved the paper. I am also grateful to Gilles Celeux for having
helped me to edit the English in the text.

\bibliographystyle{plain}
\bibliography{acompat,Slope_heuristics_regression_13}

\end{document}